\documentclass[a4paper]{amsart}

\usepackage{amsmath,amssymb}
\usepackage[dvipsnames]{xcolor}
\usepackage{pictexwd}
\usepackage{graphicx}

\newtheorem{theorem}{Theorem}[section]
\newtheorem{proposition}[theorem]{Proposition}
\newtheorem{corollary}[theorem]{Corollary}
\newtheorem{lemma}[theorem]{Lemma}
\newtheorem{algorithm}[theorem]{Algorithm}

\newtheorem{preremark}[theorem]{Remark}
\newtheorem{predefinition}[theorem]{Definition}
\newtheorem{preexample}[theorem]{Example}
\newtheorem{prenotation}[theorem]{Notation}
\newtheorem{preconjecture}[theorem]{Conjecture}

\newenvironment{remark}{\begin{preremark}\rm}{\end{preremark}}
\newenvironment{definition}{\begin{predefinition}\rm}
{\end{predefinition}}
\newenvironment{example}{\begin{preexample}\rm}{\end{preexample}}

\newenvironment{notation}{\begin{prenotation}\rm}{\end{prenotation}}

\allowdisplaybreaks
\def\OO{{\mathcal{O}}}
\def\QQ{\mathbb{Q}}

\def\AA{\mathbb{A}}

\newcommand{\M}{{\mathfrak{M}}}

\newcommand{\m}{{\mathfrak{m}}}

\let\epsilon=\varepsilon

\def\phi{{\varphi}}
\let\Psi=\varPsi
\let\Phi=\varPhi
\let\theta=\vartheta
\let\rho=\varrho

\def\LT{\mathop{\rm LT}\nolimits}

\def\ND{\mathop{\rm ND}\nolimits}
\def\AR{\mathop{\rm AR}\nolimits}
\def\GL{\mathop{\rm GL}\nolimits}

\def\Mat{\mathop{\rm Mat}\nolimits}

\def\Supp{\mathop{\rm Supp}\nolimits}
\def\Spec{\mathop{\rm Spec}\nolimits}
\def\Hilb{\mathop{\rm Hilb}\nolimits}
\def\Cot{\mathop{\rm Cot}\nolimits}

\def\GFan{\mathop{\rm GFan}\nolimits}
\def\LTGFan{\mathop{\rm LTGFan}\nolimits}

\def\Syz{\mathop{\rm Syz}\nolimits}
\def\Ker{\mathop{\rm Ker}\nolimits}

\def\wt{\mathop{\rm wt}\nolimits}

\def\lin{{\mathop{\rm lin}\nolimits}}
\def\ord{{\mathop{\rm ord}\nolimits}}

\newcommand{\Lin}{\mathop{\rm Lin}\nolimits}
\newcommand{\edim}{\mathop{\rm edim}\nolimits}

\newcommand{\BO}{\mathbb{B}_{\mathcal{O}}}

\let\To=\longrightarrow

\def\tr{^{\,\rm tr}}
\def\tfrac #1#2{{\textstyle\frac{#1}{#2}}}
\def\tsum_#1^#2{{\textstyle\sum\limits_{#1}^{#2}}}
\def\tbinom #1#2{{\textstyle\binom{#1}{#2}}}

\definecolor{red}{rgb}{1.0, 0.0, 0.0}

\def\cocoa{\mbox{\rm
  C\kern-.13em o\kern-.07 em C\kern-.13em o\kern-.15em A}}
\def\apcocoa{\mbox{\rm
A\kern-0.13em p\kern -0.07em C\kern-.13em o\kern-.07 em C\kern-.13em
o\kern-.15em A}}

\begin{document}

\title{Optimal Re-Embeddings of Border Basis Schemes}

\author{Martin Kreuzer}
\address{Fakult\"at f\"ur Informatik und Mathematik, Universit\"at
Passau, D-94030 Passau, Germany}
\email{Martin.Kreuzer@uni-passau.de}

\author{Le Ngoc Long}
\address{Fakult\"at f\"ur Informatik und Mathematik, Universit\"at
Passau, D-94030 Passau, Germany and Department of Mathematics,
University of Education - Hue University, 34 Le Loi Street, Hue City, Vietnam}
\email{lelong@hueuni.edu.vn}

\author{Lorenzo Robbiano}
\address{Dipartimento di Matematica, Universit\`a di Genova,
Via Dodecaneso 35,
I-16146 Genova, Italy}
\email{lorobbiano@gmail.com}

\date{\today}

\begin{abstract}
Border basis schemes are open subschemes of Hilbert schemes
parametrizing 0-dimensional subschemes of~$\mathbb{P}^n$ of given length.
They yield open coverings and are easy to describe 
and to compute with. Our topic is to find re-embeddings of border basis schemes 
into affine spaces of minimal dimension.

Given $P = K[x_1,\dots,x_n]$ over a field~$K$, an ideal
$I\subseteq \langle X \rangle = \langle x_1,\dots, x_n\rangle$, and a tuple~$Z$ 
of indeterminates, in previous papers the authors developed techniques for computing 
{\it $Z$-sep\-a\-ra\-ting re-embeddings} of the ideal~$I$, i.e., of isomorphisms 
$\Phi: P/I \rightarrow K[X{\setminus} Z] / (I\cap K[X{\setminus }Z])$.
Here these general techniques are developed further and improved 
by constructing a new algorithm for checking candidate tuples~$Z$ and by using
the Gr\"obner fan of the linear part of~$I$ advantageously. 

Then we apply this to the ideals defining border basis schemes $\mathbb{B}_{\mathcal{O}}$, 
where~$\mathcal{O}$ is an order ideal of terms, and to their natural 
generating polynomials. The fact that these
ideals are homogeneous w.r.t.\ the arrow grading allows us to look for suitable 
tuples~$Z$ more systematically. Using the equivalence of indeterminates
modulo the square of the maximal ideal, we compute the Gr\"obner fan of
the linear part of the ideal quickly and determine which indeterminates
should be in~$Z$ when we are looking for optimal re-embeddings.
Specific applications include re-embeddings of border basis schemes where
$\mathcal{O}\subseteq K[x,y]$ and where~$\mathcal{O}$ consists 
of all terms up to some degree.
\end{abstract}



\maketitle

\noindent{\bf IMPORTANT NOTE:} {\it This preprint will not be published. It has been
split into several parts which will be extended and published separately. When the last
part is finished, this preprint will be withdrawn.}
\bigskip

%
%

\section*{Introduction}

One of the characteristic features of Algebraic Geometry is the existence
of moduli spaces. A systematic construction of so-called Hilbert schemes
was given by A.\ Grothendieck in~\cite{Gro}.
One particularly famous and well-researched type of these moduli spaces are 
the Hilbert schemes $\Hilb_\mu(\mathbb{P}^n_K)$ parametrizing
0-dimensional subschemes of length~$\mu$ in a projective space $\mathbb{P}^n_K$ 
over a field~$K$. In Computational Algebraic Geometry, Grothendieck's construction
leads to considerable difficulties, as the defining ideals of these moduli
spaces involve too many indeterminates to allow acceptably efficient calculations
in all but the most elementary cases.

One tool which has emerged in the last 20 years to study
$\Hilb_\mu(\mathbb{P}^n_K)$ is the theory of border basis schemes.
These schemes define an open affine cover of the Hilbert scheme. For an order ideal
$\OO = \{t_1,\dots,t_\mu\}$ in $P=K[x_1,\dots,x_n]$, i.e., for
a set of terms (power products of indeterminates) which is closed under forming divisors,
the $\OO$-border basis scheme $\BO$ parametrizes all 0-dimensional schemes
in~$\mathbb{A}^n_K$ such that the residue classes of the terms in~$\OO$ form
a $K$-vector space basis of their coordinate ring. The various $\OO$-border basis
schemes for order ideals of length~$\mu$ cover $\Hilb_\mu(\mathbb{P}^n_K)$ and have
easily computable defining ideals (see~\cite{BC}, \cite{Hai}, \cite{Hui1}, \cite{Hui2}, \cite{Hui3},
\cite{KR3}, \cite{KR4}). This makes it possible to use them for explicit calculations
to check if they are non-singular (see \cite{Hui2}, \cite{KSL}), to  verify when they are 
affine cells (see \cite{Hui1}), or to compute subschemes defining particular loci
(see \cite{KLR1}, \cite{KR3}, \cite{KR4}, \cite{Rob}). But one impediment remains:
although the equations defining border basis schemes are sparse, there are still
many of them and they involve a significant number of indeterminates, making it hard
to carry out some of the necessary calculations in practice.

This leads us to the topic of this paper: we want to re-embed border basis schemes
into lower-dimensional affine spaces and to efficiently calculate those re-embeddings.
The border basis scheme $\BO$ is an affine scheme whose coordinate ring is 
of the form $K[C]/I(\BO)$, where $C=\{c_{11}, \dots, c_{\mu\nu}\}$ is a set of distinct 
indeterminates (see~Definition~\ref{def-BO}).
One obstruction to the idea of re-embedding $\BO$ is given by the cotangent space at the monomial 
point of the border basis scheme. The monomial point corresponds to the monomial ideal in~$P$ which
is generated by all terms not in~$\OO$. On the $\OO$-border basis scheme, it corresponds
to the origin defined by the maximal ideal~$\m$ generated by the residue classes of the
indeterminates in~$C$. 

A key observation is that the
vector space dimension $\dim_K(\m/\m^2)$ yields a lower bound for the
possible embedding dimensions, since it is an invariant under isomorphisms, 
and if this bound is achieved, we say that we have found
an {\it optimal} re-embedding of~$\BO$ (see Definition~\ref{def-optimal}).

The usual way re-embeddings are computed in Computational Algebraic Geometry utilizes
the calculation of elimination ideals via Gr\"obner bases. 
However, due to the large number of indeterminates,
this method is not feasible in most situations involving border basis schemes. In two
previous papers (see~\cite{KLR2}, \cite{KLR3}) the authors introduced and developed
a more efficient method based on separating tuples of indeterminates.
This method can be applied to any finitely generated $K$-al\-ge\-bra $P/I$,
where~$I$ is an ideal contained in the maximal ideal 
$\M = \langle x_1,\dots,x_n \rangle$ of~$P$. 
For a tuple $Z=(z_1,\dots,z_s)$ of indeterminates
in~$P$, we say that~$I$ is {\it $Z$-separating} if there exist a tuple of non-zero 
polynomials $(f_1,\dots,f_s)$ of~$I$ and a term ordering~$\sigma$ such that
$\LT_\sigma(f_i)=z_i$ for $i=1,\dots,s$. Given a tuple~$Z$ such that~$I$ is $Z$-separating, 
we can eliminate the indeterminates in~$Z$ and get a {\it $Z$-separating
re-embedding} $\Phi:\; P/I \longrightarrow K[Y] / (I\cap K[Y])$ where $Y= \{x_1,\dots,x_n\}
\setminus \{z_1,\dots,z_s\}$. Moreover, we described some methods to search for such 
separating tuples based on Linear Programming Feasibility (LPF) solvers and the computation
of restricted Gr\"obner fans (see~\cite{KLR3}).

The first topic in this paper, treated in Sections~1 and~2, is to improve those methods 
further in the general framework of finitely generated $K$-algebras, in particular 
if we are looking for optimal re-embeddings of~$I$.
In Section~1 we show how to check whether a given tuple of distinct indeterminates~$Z$
yields a $Z$-separating re-embedding for a given ideal~$I$. Besides the general method,
which involves a Gr\"obner basis computation (see Proposition~\ref{prop-CheckZ}), 
we provide an algorithm which works particularly well for the type of ideals considered
later in the paper, namely defining ideals of border basis schemes 
(see Algorithm~\ref{alg-checkZsepI}). Since this algorithm uses only linear algebra
techniques, it is suitable for finding good re-embeddings even in the case of large examples.

In Subsection~2.A we introduce the idea of using $\GFan(\langle\Lin_\M(I)\rangle)$.
This is the Gr\"obner fan of the {\it linear part} $\Lin_{\M}(I)$ of~$I$, 
i.e., the ideal generated by the 
homogeneous linear parts of the elements of~$I$. It is used to find tuples~$Z$ which may
yield separating, possibly optimal re-embeddings of~$I$ (see Algorithm~\ref{alg-Emb}).
Using Remark~\ref{rem-onlysufficient} and Example~\ref{ex-isotoK[x]} we show that 
our methods may fail to detect {\it all}\/ optimal embeddings.
In Subsection~2.B we look at the Gr\"obner fan of an ideal~$I_L$ generated by homogeneous
linear polynomials more closely. First we reduce the task to check whether a given
tuple of indeterminates~$Z$ consists of leading terms of a reduced Gr\"obner basis
of~$I_L$ to a linear algebra computation (see Proposition~\ref{prop-charLinZ}). Then
we reduce the calculation of the Gr\"obner fan of~$I_L$
to the task of finding the maximal minors of a matrix (see Theorem~\ref{thm-bijection} 
and Remark~\ref{rem-compGFan}). We also note that this task can be tackled via any 
method for computing the bases of a linear matroid (see Remark~\ref{rem-compGFan}.b).

\medskip
The actual application of these methods to finding optimal re-em\-bed\-dings of
border basis schemes commences in Section~3. When we are given an order ideal of terms 
$\OO = \{ t_1,\dots,t_\mu\}$ in a polynomial ring $P=K[x_1,\dots,x_n]$ over a field~$K$, 
we let $\partial\OO = \{b_1,\dots,b_\nu\}$ be its {\it border}, and $C = \{ c_{ij} \mid
i=1,\dots,\mu;\; j=1,\dots,\nu\}$ the set of indeterminates which appear as coefficients
in the {\it generic $\OO$-border prebasis} $G = \{ g_1,\dots,g_\nu \}$,
where $g_j = b_j -\sum_{i=1}^\mu c_{ij} t_i$ for $j=1,\dots,\nu$. 
Then the defining ideal $I(\BO)$ of the $\OO$-border basis scheme $\BO$ is generated 
in~$K[C]$ by two types of polynomials:
\begin{enumerate} 
\item[(1)] For border terms $b_j = x_\ell b_{j'}$, called {\it next-door neighbours},
there is a tuple $\ND(j, j')$ of $\mu$ polynomials which expresses the 
fact that the syzygy $b_j - x_\ell b_{j'} = 0$ of border terms can be lifted to a syzygy
of~$G$.

\item[(2)] For border terms $b_j = x_\ell t_m$ and $b_{j'} = x_k t_m$, called
{\it across-the-rim neighbours}, there is a tuple $\AR(j, j')$ of $\mu$ polynomials
which expresses the fact that the syzygy  $x_k b_j - x_\ell b_{j'} = 0$ 
of border terms can be lifted to a syzygy of~$G$.
\end{enumerate}
Altogether, both types of polynomials form the set of {\it natural generators} of~$I(\BO)$.
In Proposition~\ref{prop-lin_cij} we recall the detailed description of the linear parts
of the natural generators of~$I(\BO)$ given in~\cite{KSL}.

Then Section~4 introduces a further important structure governing the vanishing ideal
of a border basis scheme, namely the {\it arrow degree}. It is the $\mathbb{Z}^n$-grading 
on~$K[C]$ given by $\deg_W(c_{ij}) = \log(b_j) -\log(t_i)$ for $i=1,\dots,\mu$ and $j=1,\dots,\nu$.
The natural generators, and thus the entire ideal $I(\BO)$, are homogeneous with respect to
the arrow grading. Unfortunately, the arrow grading is in general not positive, and in fact
not even non-negative, so that we cannot avail ourselves of the many nice properties of positive
or non-negative gradings. However, the arrow degree allows us to visualize the linear parts
of the natural generators of~$I(\BO)$ nicely (see Remark~\ref{rem-LinParts}) and to describe
the terms in the quadratic parts of the natural generators in great detail
(see Proposition~\ref{prop-QuadrParts}). 

Another fundamental ingredient is provided in Section~5. 
If we let $\m = \langle \bar{c}_{ij} \rangle$
denote the homogeneous maximal ideal in $K[C]/ I(\BO)$, the cotangent space of~$\BO$ at the monomial
point is given by $\Cot_{\m}(B_{\OO}) = \m / \m^2$. As mentioned before, 
the dimension of this vector space is a lower bound for the embedding dimension of~$\BO$. 
It also yields a useful
equivalence relation on the indeterminates~$c_{ij}$. We say that $c_{ij}$ and $c_{i' j'}$
are {\it cotangent equivalent} if $\bar{c}_{ij} + \m^2 = \bar{c}_{i' j'} + \m^2$. 
Then there are three main types of cotangent equivalence classes.
The indeterminates~$c_{ij}$ such that $\bar{c}_{ij}\in\m^2$ are called {\it trivial} and form
the {\it trivial equivalence class}.  A non-trivial indeterminate~$c_{ij}$ is called {\it basic} 
if its cotangent equivalence class (which is then also called basic) consists only of~$c_{ij}$. 
Finally, a non-trivial indeterminate~$c_{ij}$ is called {\it proper} if its equivalence class 
(which is then also called proper) contains at least two elements.
An interesting property of this equivalence relation is that it provides a 
beautiful description of $\GFan(\langle\Lin_\M(I(\BO))\rangle)$ and a 
method to compute it (see Theorem~\ref{thm-shapeofSsigma}) which is more
efficient than the general method described in Subsection 2.A.

Moreover, an important application of these equivalence classes to the task of 
finding optimal $Z$-separating re-embeddings of~$I(\BO)$ is that, 
to find a suitable tuple~$Z$, we have to put into $Y = C {\setminus} Z$ 
all basic indeterminates~$c_{ij}$ and exactly one from each proper equivalence
class (see Theorem~\ref{thm-classifycij}). Since the cotangent equivalence
classes can be calculated very quickly (see~\cite{KSL}), we obtain a rather fast
method of checking for optimal $Z$-separating re-embeddings of~$I(\BO)$ in 
Algorithm~\ref{alg-compEmb}.

The final three sections of the paper are devoted to examining special cases
in which we can (or cannot) find optimal re-embeddings of~$\BO$. Section~6 starts
with the case of {\it planar} border basis schemes, i.e., the case $n=2$.
Most of the results in this section are inspired by, improve upon, or make explicit 
the excellent and thorough study of this setting in~\cite{Hui1}. 
A key concept here is the notion of an {\it $x_\ell$-exposed} indeterminate~$c_{ij}$.
When we form a product $x_\ell g_j$ as part of lifting a neighbour syzygy or an across-the-rim
syzygy, where $g_j$ is an element of the generic $\OO$-border prebasis, some coefficients
of the result will be coefficients of a border term, and the indeterminates  appearing
in them are called the $x_\ell$-exposed indeterminates. By varying~$\ell$ and~$j$, 
we obtain the set $C^{\rm ex}$ of all exposed indeterminates in~$C$ 
(see Definition~\ref{def-exposedIndets} and Proposition~\ref{prop:CharExposed}). 

The main result of the section is that all non-exposed indeterminates can be eliminated,
i.e., that the exposed indeterminates generate $K[C] / I(\BO)$. In~\cite{Hui1}, Prop.~7.1.2, 
this is shown via the Hilbert-Burch Theorem. Instead, here we construct an explicit
Weight Assignment Algorithm~\ref{alg-WeightAssign} which allows us to find, 
for $Z = C\setminus C^{\rm ex}$, a $Z$-separating tuple of natural generators of~$I(\BO)$
and the corresponding elimination ordering~$\sigma$. Thus all non-exposed indeterminates
$c_{ij}$ can be eliminated (see Theorem~\ref{thm-elimNex}). Consequently, we can also improve
Algorithm~\ref{alg-compEmb} by restricting the choice of indeterminates for $Y=C{\setminus} Z$
in each proper cotangent equivalent class to exposed indeterminates only 
(see Proposition~\ref{prop-E0plusNonexp}). 
However, as for instance Example~\ref{ex-Section5-continued} shows, eliminating the 
non-exposed indeterminates is in general not sufficient to achieve optimal embeddings.
Furthermore, Remark~\ref{rem:intnonexp} shows that all exposed indeterminates are rim 
indeterminates (see Definition~\ref{def-interiorrim}), but in 
Example~\ref{ex-Section5-continued} some rim indeterminates are not exposed.
An explicit application  of the Weight Assignment Algorithm~\ref{alg-WeightAssign} 
is to the family of order ideals called {\it planar boxes} (see~Remark~\ref{rem-twoboxes}).

Another family of examples in which we can compute optimal separating re-embeddings 
of~$I(\BO)$ is examined in Section~7, where we treat {\it simplicial} order ideals~$\OO$ 
of type~$d$, i.e., order ideals consisting of all terms of degree $\le d$ for some $d>0$.
Here we prove that one can optimally re-embed $I(\BO)$ by using the tuple~$Z$
of all interior indeterminates~$c_{ij}$, i.e., all indeterminates~$c_{ij}$ such that
no multiple $x_k t_i$ is in~$\partial\OO$ (see Proposition~\ref{prop-simplicialZ}). 
As an immediate consequence, we see that
border basis schemes of simplicial order ideals are affine cells if $n=2$ and singular
schemes if $n\ge 3$ (see Corollary~\ref{cor-affand singular}).

The last section comprises a small collection of interesting further examples
of re-embeddings of border basis schemes. In Example~\ref{ex-theLshapeCont}
we look at the order ideal $\OO = \{ 1,\, y,\, x,\, y^2,\, x^2\}$ in two indeterminates
which we call the {\it L-shape}. Its border basis scheme~$\BO$ is shown to be a
10-dimensional smooth variety which can be embedded as a complete intersection 
of codimension~2 in $\mathbb{A}^{12}$. This embedding is the best $Z$-separating re-embedding,
and we conjecture that~12 is in fact the embedding dimension of~$\BO$ in this case.
In particular, this would imply that~$\BO$ is not an affine cell.

Two similar looking further examples exhibit very different behaviors. For the order ideal
$\OO = \{ 1,\, z,\, y,\, x,\, z^2 \}$ we get a singular border basis scheme~$\BO$
which has an optimal $Z$-separating re-embedding (see Example~\ref{ex-13z^2}),
while for $\OO = \{ 1,\, z,\, y,\, x,\, xy\}$ we get a smooth border basis scheme
for which we cannot find an optimal $Z$-separating re-embedding.

All algorithms mentioned in this paper were implemented in the computer 
algebra system CoCoA (see~\cite{CoCoA}) which was also used to perform the calculations underlying 
most examples. Their use was essential in the discovery of properties and features which eventually 
evolved into theorems or disproved previous conjectures. A large folder of supporting functions 
and explicit examples calculated using CoCoA is available from the authors upon request.
The general notation and definitions in this paper follow~\cite{KR1} and~\cite{KR2},
and for specific terminology relating to border basis schemes we adhere to our earlier
papers~\cite{KLR1}, \cite{KR3}, \cite{KR4},  and~\cite{KSL}.

\bigskip\bigbreak
%
%

\section{$Z$-Separating Re-embeddings of Affine Schemes}
\label{Z-Separating Re-embeddings of Affine Schemes}

In this paper we let $K$ be an arbitrary field, let $P=K[x_1,\dots,x_n]$,
and let $\M=\langle x_1, \dots, x_n \rangle$. The tuple formed by the indeterminates 
of~$P$ is denoted by $X =(x_1, \dots, x_n)$. Moreover, let $1\le s\le n$, 
let $z_1,\dots, z_s$ be pairwise distinct indeterminates in~$X$, 
and let $Z =(z_1, \dots, z_s)$.
Denote the remaining indeterminates by 
$\{y_1,\dots,y_{n-s}\} = \{x_1,\dots,x_n\} \setminus \{z_1,\dots,z_s\}$,
and let $Y =(y_1,\dots, y_{n-s})$. Committing a slight abuse of notation, 
we shall also write $Y = X \setminus Z$.
The monoid of terms in~$P$ is denoted by $\mathbb{T}^n = \{ x_1^{\alpha_1} \cdots
x_n^{\alpha_n} \mid \alpha_i \ge 0\}$. 
Given a term ordering $\sigma$ on~$\mathbb{T}^n$,
its restriction to $\mathbb{T}(y_1, \dots, y_{n-s})$ is denoted by~$\sigma_Y$.

\medskip
\subsection*{1.A.\ \textit{Z}-Separating Tuples and \textit{Z}-Separating 
Re-embeddings}$\mathstrut$

In~\cite{KLR2} and~\cite{KLR3}, the authors examined re-embeddings of affine algebras,
i.e., isomorphisms with presentations requiring fewer $K$-algebra
generators. In particular, it was shown that the following situation leads to
such re-embeddings.

\begin{definition}\label{def-Zseparating}
Let $I$ be an ideal in~$P$ which is contained in~$\M$, and let $Z=(z_1,\dots,z_s)$
be a tuple of distinct indeterminates in~$X$. We say that the ideal ~$I$ 
is {\bf $Z$-separating} if there exist a term ordering~$\sigma$ on~$\mathbb{T}^n$ 
and $f_1,\dots,f_s \in I \setminus \{0\}$ such that
$\LT_\sigma(f_i)= z_i$ for $i=1,\dots,s$. In this situation~$\sigma$ is called a 
{\bf $Z$-separating term ordering} for~$I$, and the tuple $(f_1,\dots,f_s )$ is called
a \bf $Z$-separating tuple.
\end{definition}

Recall that, in the setting of the definition, the reduced $\sigma$-Gr\"obner basis
of the ideal $\langle f_1,\dots,f_s\rangle$ is {\bf coherently $Z$-separating}, 
i.e., it is of the form $\{g_1,\dots,g_s\}$ where we have $\LT_\sigma(g_i)=z_i$ 
for $i=1,\dots,s$, and where~$z_i$ divides no other term in~$g_i$ 
or in~$g_j$ for $i,j\in\{1,\dots,s\}$ and $j\ne i$ (see~\cite{KLR3}, Prop.~2.2).
Furthermore, recall that $\Lin_{\M}(f)$ denotes the homogeneous component of degree~1 
of a polynomial $f\in\M$ and $\Lin_{\M}(I) = \langle \Lin_{\M}(f) \mid f\in I \rangle_K$ 
is the {\bf $\M$-linear part} of~$I$. Moreover, if we denote the image of~$\M$ in~$P/I$
by~$\m$, it was shown in~\cite{KLR2}, Prop.~1.8.b, that
$P_1 / \Lin_{\M}(I) \cong \m/\m^2$, where $\Cot_{\m}(B_\OO) = \m / \m^2$ 
is the {\bf cotangent space} of~$P/I$ at the origin, and in\cite{KLR2}, Prop.~1.9
that $\Lin_{\M}(I) $ is easy to compute, since it is equal 
to $\langle \Lin_{\M}(f_1),...,\Lin_{\M}(f_s)\rangle_K$, where $\{f_1, \dots, f_s\}$ 
is any set of  generators of~$I$.

In this setting we obtain the following kind of re-embedding of~$I$.

\begin{proposition}\label{prop-cohSepPoly}
Let~$I$ be an ideal in~$P$ which is contained in~$\M$. Assume that there 
exists a tuple of distinct indeterminates $Z=(z_1,\dots,z_s)$ in~$X$ such that~$I$
is $Z$-separating. Let~$\sigma$ be a $Z$-separating term ordering for~$I$,
let $Y=X \setminus Z$, and let~$S_\sigma$ be the minimal set of terms generating
$\LT_\sigma(\langle \Lin_{\M}(I)\rangle )$.
\begin{enumerate}
\item[(a)] The reduced $\sigma$-Gr\"obner basis of~$I$ is of the form 
$G = \{ z_1 - h_1,\allowbreak \dots,\, z_s - h_s, \allowbreak g_1,\, \dots,\, g_r\}$
with $h_i, g_j \in K[Y]$. 

\item[(b)] The $K$-algebra homomorphism $\Phi: P/I \longrightarrow K[Y] / 
(I \cap K[Y])$ given by $\Phi(\bar{x}_i) = \bar{x}_i$ for 
$x_i\in Y$ and $\Phi(\bar{x}_i) = \bar{h}_j$
for $x_i = z_j \in Z$ is an isomorphism of $K$-algebras,  called the
{\bf $Z$-separating re-embedding} of~$I$.

\item[(c)] We have $\langle Z\rangle \subseteq  \LT_\sigma(\langle \Lin_{\M}(I) \rangle)$,
and hence $Y \supseteq X \setminus S_\sigma$.  

\item[(d)] The residue classes of the indeterminates in $Y = X \setminus Z$
generate the cotangent space $P_1/\Lin_{\M}(I)\cong \m / \m^2$ of $P/I$ at the origin. 
\end{enumerate}
\end{proposition}

\begin{proof}
Claims~(a) and~(b) follow from~\cite{KLR2}, Theorem~2.13. Let us prove~(c).
By assumption, there exist 
polynomials $f_1, \dots, f_s\in I$ such that they form the reduced $\sigma$-Gr\"obner 
basis of $\langle f_1,\dots,f_s\rangle$ and such that $z_i = \LT_\sigma(f_i)$ 
for $i=1,\dots, s$. Then we have $z_i = \LT_\sigma(\Lin_\M(f_i))$ for $i=1,\dots, s$, 
and thus $\langle Z\rangle \subseteq \LT_\sigma(\langle\Lin_{\M}(I)\rangle)$. 
Claim~(d) follows from~(c), since the set of the residue classes 
of the elements in $X \setminus S_\sigma$ is a $K$-basis of $P_1/\Lin_{\M}(I)$.
\end{proof}

Next we recall a definition which plays a fundamental role in this paper
(see~\cite{KLR2}, Definition 3.3).

\begin{definition}\label{def-optimal}
Let $P=K[x_1,\dots,x_n]$ be a polynomial ring over a field~$K$, and let~$I$ 
be a proper ideal in~$P$.
\begin{enumerate}
\item[(a)] A $K$-algebra isomorphism $\Psi:\; P/I \longrightarrow P'/I'$,
where $P'$ is a polynomial ring over~$K$ and~$I'$ is an ideal in~$P'$, is
called a {\bf re-embedding} of~$I$.

\item[(b)] A re-embedding $\Psi:\; P/I \longrightarrow P'/I'$ of~$I$ is called
{\bf optimal} if every \hbox{$K$-al}\-ge\-bra isomorphism $P/I \longrightarrow P''/I''$ 
with a polynomial ring~$P''$ over~$K$ and an ideal~$I''$ in~$P''$ 
satisfies the inequality $\dim(P'')\ge \dim(P')$.

\end{enumerate}
\end{definition}

When we are only looking for $Z$-separating re-embeddings of~$I$ which are optimal,
the above proposition yields the following characterization.

\begin{corollary}\label{cor-cohSepPoly}
In the setting of the proposition, assume that ${s=\#Z}$ is equal to the $K$-vector space 
dimension of $\Lin_{\M}(I)$. Then the following claims hold.
\begin{enumerate}
\item[(a)] The map~$\Phi$ is an optimal re-embedding of~$I$.

\item[(b)] We have $\langle Z\rangle = \LT_\sigma(\langle \Lin_{\M}(I) \rangle)$,
and hence $Y = X \setminus S_\sigma$. 

\item[(c)] The residue classes of the indeterminates in $Y = X \setminus Z$
form a $K$-vector space basis of~the cotangent space 
$P_1/\Lin_{\M}(I)\cong \m / \m^2$ of $P/I$ at the origin. 

\end{enumerate}
\end{corollary}

\begin{proof}
Claim~(a) is a consequence of~\cite{KLR2}, Cor.~4.2. Let us prove~(b). 
By claim~(c) of the proposition, we have 
$\langle Z\rangle \subseteq  \LT_\sigma(\langle \Lin_{\M}(I) \rangle)$.
Since we have $s = \dim_K(\langle Z\rangle_K) = \dim_K(\Lin_{\M}(I))$, we deduce 
that~$Z$ minimally generates $\LT_\sigma(\langle \Lin_{\M}(I) \rangle)$ 
which concludes the proof. 

Finally, let us prove (c). By~(b), we have $Z = S_\sigma$. As mentioned in
the proof of claim~(d) of the proposition, the set of the residue classes 
of the elements in $X \setminus S_\sigma$ is a $K$-basis 
of $P_1/\Lin_{\M}(I)$. This implies the claim.
\end{proof}

An important situation in which we obtain an optimal re-embedding is described 
by the following corollary.

\begin{corollary}\label{cor-affinescheme}
In the setting of the proposition, the following conditions are equivalent.
\begin{enumerate}
\item[(a)] The map~$\Phi:\; P/I \longrightarrow K[Y] / (I \cap K[Y])$
is an optimal re-embedding, and the ring $P_{\M}/ I_{\M}$ is a regular local ring.

\item[(b)] We have $I \cap K[Y] = \{0\}$.

\end{enumerate}
In particular, if these conditions are satisfied, then~$\Phi$ is an isomorphism
with a polynomial ring, the scheme $\Spec(P/I)$ is isomorphic to an affine space,
and the set $G = \{z_1-h_1,\dots, z_s-h_s\}$ is a minimal set of generators of~$I$.  
\end{corollary}

\begin{proof} 
The claims  follow from~\cite{KLR3}, Proposition~6.7.
\end{proof}

Sometimes an optimal separating re-embedding
cannot be achieved using a given set of generators of~$I$, 
as the following easy example shows.

\begin{example}\label{ex-easynongens}
Let $P = \QQ[x,y]$, and let $I=\langle x,\ y\rangle$.
Using  $\{x,\, y\}$ as a set of generators of~$I$, we obtain the 
re-embedding $P/I\cong \QQ$ which cannot be achieved using the set $\{x+y,\ x-y\}$.
\end{example}

If the conditions of Corollary~\ref{cor-affinescheme} are satisfied, 
we have the following additional information.

\begin{proposition}\label{prop-minsysgens}
Let $F=(f_1,\dots,f_r)$ be a tuple of non-zero polynomials in~$\M$,
and let $I=\langle f_1,\dots,f_r \rangle$. Suppose that $P_{\M}/I_{\M}$ is
a regular local ring and that, for $s= \dim_K (\Lin_{\M}(I))$, we have a tuple~$Z$
of~$s$ indeterminates in~$X$ such that there exist $i_1,\dots,i_s\in \{1,\dots,r\}$
for which the tuple $F_Z = (f_{i_1}, \dots, f_{i_s})$ is $Z$-separating.

\begin{enumerate}
\item[(a)]
The tuple $F_Z$ is a Gr\"obner basis of~$I$ with respect to any elimination ordering
for~$Z$ and a minimal system of generators of~$I$. 

\item[(b)]
The tuple $F_Z$ represents $P/I$ as a complete intersection.

\item[(c)] The reduced Gr\"obner bases of all elimination orderings for~$Z$ coincide.
\end{enumerate}
\end{proposition}

\begin{proof} First we prove~(a). Since~$F_Z$ is separating, we have $\langle Z\rangle = 
\LT_\sigma(\langle F_Z\rangle)$ for every elimination ordering~$\sigma$
for~$Z$. The theorem on the computation of elimination modules (cf.~\cite{KR1},
Theorem 3.4.5) implies that every $\sigma$-Gr\"obner basis of~$I$ consists
of polynomials with leading terms in $\langle Z\rangle$ and polynomials in~$K[Y]$
generating $I \cap K[Y]$.
By Corollary~\ref{cor-affinescheme}, we have $I\cap K[Y] = \{0\}$, whence
it follows that~$F_Z$ is in fact a minimal $\sigma$-Gr\"obner basis of~$I$. 
Consequently, $F_Z$ is a system of generators of~$I$. 

To prove (b) we note that the fact that the leading terms of $f_{i_1},\dots,f_{i_s}$ 
form a regular sequence implies that~$F_Z$ is actually a minimal system of generators of~$I$ 
and a regular sequence.

Finally, to prove claim~(c), we use the definitions and results of~\cite{KLR3}, Section~5.
By Thm.~5.5, the map
$\Gamma_Z:\; \GFan_Z(I)\longrightarrow \GFan(I\cap K[Y])$ is bijective. Since $I\cap K[Y] =\{0\}$,
we have  $\GFan(I\cap K[Y]) =\{\emptyset\}$. Therefore $\GFan_Z(I)$, which is not empty, 
has cardinality one.
\end{proof}

The following observation will be useful later.

\begin{remark}\label{rem:differentF}
The entries of every $s$-tuple $F_Z$ of elements in~$I$ can be considered as elements in a set of 
generators of~$I$. Consequently, in the setting of the proposition, it is possible that there are 
several $s$-tuples $F_Z$ which are $Z$-separating. However, claim~(c) shows that all of them generate~$I$
and yield the same reduced Gr\"obner basis with respect to any elimination ordering for~$Z$.
Further observations on this topic are contained in Remark~\ref{rem:goodweights} and Example~\ref{ex-123}.
\end{remark}

Let us apply the above proposition in a concrete case.

\begin{example}\label{ex-mingensaffine}
Let $P = \QQ[x,y,z,w]$, let $f_1=  w^2 +x -y +3z$, $f_2 =  zw^2 +w^3 +y$, 
$f_3=w^3 -xz +yz -3z^2 +y$, and let $I =\langle f_1, f_2, f_3\rangle$. 
We have $\dim(P_{\M}/I_{\M}) = 2$ and $\dim(P_{\M}/\Lin_{\M}(I))=2$ which implies 
that $P_{\M} / I_{\M}$ is a regular local ring.
Using Corollary 3.5 in~\cite{KLR3}, we check that $\{f_1, f_2\}$ 
is $(x,y)$-separating.
From the proposition we deduce that~$I$ is minimally generated by $\{f_1, f_2\}$.

Moreover, if $\sigma$ is any elimination ordering for $(x,y)$, 
the reduced $\sigma$-Gr\"obner basis of~$I$ 
is $G = \{ x +z w^2  +w^3  +w^2 +3z,\;  y +zw^2 +w^3  \}$.  Using~$G$,
we construct an isomorphism $\phi:\; P/I \cong \QQ[z, w]$ given by 
$\phi(\bar{x}) = -(z w^2 +w^3 +w^2 +3z)$, $ \phi(\bar{y}) = -(zw^2 +w^3)$, 
$\phi(\bar{z}) = z$, and $\phi(\bar{w})=w$.
\end{example}

\medskip
\subsection*{1.B.\ Checking \textit{Z}-Separating Tuples}
$\mathstrut$

An important algorithmic problem is to determine whether a given tuple
of distinct indeterminates~$Z$ is separating for an ideal~$I$ in~$P$
which is contained in~$\M$ and given by a set of generators. The following
proposition provides a general answer.

\begin{proposition}\label{prop-CheckZ}
Let~$I$ be an ideal in~$P$ which is contained in~$\M$, and
let $Z=(z_1,\dots,z_s)$ be a tuple of distinct indeterminates in~$X$.
Then the following conditions are equivalent.
\begin{enumerate}
\item[(a)] The ideal~$I$ is $Z$-separating.

\item[(b)] For every elimination ordering~$\sigma$ for~$Z$,
we have $\langle Z\rangle \subseteq \LT_\sigma(I)$.

\item[(c)] There exists an elimination ordering~$\sigma$ for~$Z$ such that
$\langle Z\rangle \subseteq \LT_\sigma(I)$.
\end{enumerate}
\end{proposition}

\begin{proof}
To show that~(a) implies~(b), we note that, by definition, there exists
a $Z$-separating term ordering~$\sigma$ for~$I$. By~\cite{KLR3}, Remark 4.3, it
follows that any elimination ordering for~$Z$ is then also a 
$Z$-separating term ordering for~$I$. Now the claim follows 
from~\cite{KLR3}, Proposition~4.2.

Since~(b) obviously implies~(c), it remains to prove (c)$\Rightarrow$(a).
The fact that~$\sigma$ is an elimination ordering for~$Z$ implies that
the reduced $\sigma$-Gr\"obner basis of~$I$ has the shape given in~\cite{KLR3},
Proposition 4.2.d. Thus the claim follows from this proposition.
\end{proof}

Clearly, condition~(b) of this proposition allows us to check effectively whether~$I$
is $Z$-separating for a given tuple~$Z$, but the required Gr\"obner basis computation
may be too costly, in particular if the given ideal is not $Z$-separating. 
However, the following example shows that apparently 
there is no generally applicable better method.

\begin{example}\label{ex-GBnecessary}
Consider the ring $P=\QQ[x,y,z]$, the tuple $Z=(x)$, and the ideal 
$I = \langle f_1,\dots,f_{10}\rangle$, where 
$$
\begin{array}{lcl}
f_1 &=&  xy^2 +\frac{1}{2}y^3 -\frac{1}{2}y^2z -x^2 -\frac{1}{2}xy -y^2 +\frac{1}{2}xz +x, \cr
f_2 &=&  y^2z^2 +3y^3 -4y^2z -xz^2 -3xy +4xz  ,\cr
f_3 &=&  y^3z -xyz -y^2z +xz   ,\cr
f_4 &=&  y^4 -xy^2 -y^3 +xy   ,\cr
f_5 &=&  x^2y^2 -x^3  ,\cr
f_6 &=&  x^3 +\frac{1}{2}x^2y +xy^2 +\frac{1}{2}y^3 -\frac{1}{2}x^2z -\frac{1}{2}y^2z -x^2 -y^2 ,\cr
f_7 &=&  x^2z^2 +y^2z^2 +3x^2y +3y^3 -4x^2z -4y^2z   ,\cr
f_8 &=&  x^2yz +y^3z -x^2z -y^2z  ,\cr
f_9 &=&  x^2y^2 +y^4 -x^2y -y^3   ,\cr
f_{10} &=&  x^4 +x^2y^2   
\end{array}
$$
At first glance, the ideal~$I$ does not appear to be $Z$-separating, even if we use
linear combinations of the generators. However, for an elimination ordering~$\sigma$ for~$Z$, the
reduced $\sigma$-Gr\"obner basis of the ideal~$I$ is $\{ x-y^2,\, y^4+y^2 \}$, 
and this proves that~$I$ is indeed $Z$-separating.
\end{example}

To have more efficient methods for checking if a given ideal is $Z$-sep\-a\-ra\-ting
for a given tuple~$Z$, which may or may not work in special cases, we can try to use one of
the following techniques.

\begin{remark}\label{rem-randomFforZ}
Let $I= \langle f_1,\dots,f_r\rangle$ be an ideal in~$P$ which is generated by polynomials 
$f_i \in \M$, and let~$Z=(z_1,\dots,z_s)$ be a tuple of indeterminates in~$X$.
\begin{enumerate}
\item[(a)] To check condition~(b) of the proposition, it may suffice to compute a partial
Gr\"obner basis of~$I$, or to use a timeout for the calculation.

\item[(b)] We can also hope to find suitable polynomials $f_{i_1},\dots,f_{i_s}\in I$
for which the tuple $(f_{i_1},\dots,f_{i_s})$ is $Z$-separating among the given 
generators of~$I$. More precisely, we may proceed as follows:
\begin{itemize}
\item[(1)] For all tuples $F_Z = (f_{i_1}, \dots, f_{i_s})$
such that $i_1,\dots,i_s \in \{1,\dots,r\}$ are $s$ distinct indices
and $z_j \in \Supp(f_{i_j})$ for $j=1,\dots,s$, perform the following step.

\item[(2)] Using an LPF-solver as in~\cite{KLR3}, Cor. 3.5, check whether the tuple~$F_Z$
is $Z$-separating. If it is, return the tuple of weights calculated by the
LPF-solver and stop.
\end{itemize}

\item[(c)] The loop suggested in~(b) may have too many iterations to be practically
feasible. In our implementation and experiments, we therefore used
the following heuristics: create suitable tuples $(f_{i_1},\dots,
f_{i_s})$ randomly and perform step~(2) for every tuple
after it has been created. If the LPF-solver is not successful,
try another random example until one tuple $F_Z$ works or
until a certain number of unsuccessful tries has been performed.

Even for large examples, these random choices
 lead to success in many cases, and setting
the number of random tries to several hundred (or several thousand in
the case of really large examples) resulted in a high hit rate of
correct answers.
\end{enumerate}
\end{remark}

In view of the difficulties involved in computing Gr\"obner
bases for the ideals considered in this paper, the following algorithm offers
an efficient method for finding $Z$-separating re-embeddings
in the case of ideals generated by polynomials of degree two. This case will be
exactly the one we encounter when we study re-embeddings of border basis schemes.

\begin{algorithm}{\bf (Checking for $Z$-Separating Re-embeddings)}\label{alg-checkZsepI}\\
Let $P=K[x_1,\dots,x_n]$, let $I = \langle f_1,\dots,f_r\rangle$ be an ideal in~$P$ 
generated by polynomials of the form $f_i = \ell_i + p_i \in\M$, where $\ell_i \in P_1$
and $p_i \in P_2$ for $i=1,\dots,r$, and let $Z=(z_1,\dots,z_s)$ be a tuple of distinct
indeterminates in~$X$. Consider the following sequence of instructions.
\begin{enumerate}
\item[(1)] Let $Y = X \setminus Z$, set $w_i=0$ for every $i\in\{1,\dots,n\}$
such that $x_i \in Y$, let $\widehat{Z} = Z$, let $\hat{f}_i=f_i$
for $i=1,\dots,r$, and let $d=1$.

\item[(2)] For $i=1,\dots,r$, delete in~$\hat{f}_i$ every term which is not divisible by an
indeterminate in~$\widehat{Z}$ and write the resulting polynomial in the form
$\hat{f}_i = \hat{\ell}_i + \hat{p}_i$ where $\hat{\ell}_i \in P_1$ and $\hat{p}_i \in P_2$
for $i=1,\dots,r$.

\item[(3)] If $\dim_K( \langle \hat{\ell}_1, 
\dots, \hat{\ell_r} \rangle_K) < s$ then return {\tt "No Weights Found"} and stop.
Otherwise, interreduce the polynomials $\hat{f}_i$ $K$-linearly such that their
linear parts are distinct indeterminates of~$\widehat{Z}$.  
Rename the resulting polynomials and write them in the form
$\hat{f}_i = z_i + \hat{p}_i$ for $i=1,\dots,s$ and $\hat{f}_j = q_{j-s}$ 
for $j=s+1,\dots,r$. Here we have $\hat{p}_i, q_j \in P_2$, 
and we let $u = r - s$ as well as $\hat{s}=s$.

\item[(4)] Interreduce the polynomials $q_1,\dots,q_u$ $K$-linearly 
such that they have distinct leading terms w.r.t.\ {\tt Lex}. Remove
all zero polynomials, redefining~$u$ if necessary, and call the resulting polynomials
$q_1,\dots,q_u$ again.

\item[(5)] Reduce the polynomials $\hat{f}_1,\dots,\hat{f}_{\hat{s}}$ $K$-linearly
against $q_1,\dots, q_u$ and get the polynomials $\tilde{f}_1,\dots, \tilde{f}_{\hat{s}}$.

\item[(6)] If none of the resulting polynomials $\tilde{f}_j$ 
is an indeterminate from~$\widehat{Z}$, return {\tt "No Weights Found"} and stop.
Otherwise, let $\widetilde{Z}$ be the set of all indeterminates $z_i$ which equal one of the
polynomials $\tilde{f}_j$, set $w_k=d$ for each of these indeterminates~$z_i=x_k$,
and remove from~$\widehat{Z}$ all indeterminates from~$\widetilde{Z}$.
Let $\hat{s} = \# \widehat{Z}$, let $\hat{f}_1, \dots,\hat{f}_{\hat{s}}$ 
be the remaining polynomials $\tilde{f}_j$, and denote $2d+1$ by~$d$.

\item[(7)] Repeat the following steps until $\widehat{Z} = \emptyset$. Then return
 $W=(w_1,\dots,w_n)$ and stop.

\item[(8)] For $i=1,\dots,\hat{s}$ and $j=1,\dots,u$, remove from $\hat{f}_i$ and $q_j$
every term which is not divisible by one of the indeterminates in~$\widehat{Z}$.

\item[(9)] Perform steps (4), (5), and~(6).
\end{enumerate}
This is an algorithm which, if it is successful, returns a tuple of non-negative 
weights $W=(w_1,\dots,w_n)\in \mathbb{N}^n$ such that every term ordering~$\sigma$ 
which is compatible with the grading
given by~$W$ is an elimination ordering for~$Z$ with the property that there exists a
coherently $Z$-separating tuple of polynomials $(g_1,\dots,g_s)$ in~$I$ satisfying 
$\LT_\sigma(g_i)=z_i$ for $i=1,\dots,s$. 
\end{algorithm}

\begin{proof}
In view of the two cases examined in step~(6), the algorithm either stops or
the number of elements in~$\widehat{Z}$ decreases in each iteration of the loop
in steps {(7)--(9)}. Hence the algorithm is finite, and it remains to prove correctness
if it returns a tuple $W=(w_1,\dots,w_n)$. 

Let~$\sigma$ be a term ordering which is compatible with the grading given by~$W$.
We want to show that for every indeterminate $z_i\in Z$ there exists
a polynomial $g_i \in I$ such that $\LT_\sigma(g_i)=z_i$ and such that $\{g_1,\dots,g_s\}$
is the reduced $\sigma$-Gr\"obner basis of~$\langle g_1,\dots,g_s\rangle$.
In step~(1), the indeterminates in~$Y = X\setminus Z$ are given weight zero,
and in step~(2) all terms which are products of indeterminates from~$Y$ are removed.
Since we are only performing $K$-linear combinations of the given generators
$f_1,\dots,f_r\in I$, and since every indeterminate in~$Z$ will get a positive weight,
the deleted terms will surely be smaller than any term which still needs to get a weight.

Clearly, if the linear part of $\hat{I} = \langle \hat{f}_1,\dots, \hat{f}_r \rangle$
has a dimension smaller than~$s$, then the desired
polynomials $g_1,\dots,g_s$ do not exist in~$I$. Therefore step~(3) correctly stops
the algorithm in this case. If this linear part has dimension~$s$, then there exist
$K$-linear combinations $\hat{g}_i$ of $\hat{f}_1,\dots, \hat{f}_r$ such that 
$\Lin_\M(\hat{g}_i) = z_i$ for $i=1,\dots,s$ and $\Lin_{\M}(\hat{g}_i)=0$ for $i>s$.
Such polynomials are correctly found in step~(3), and the homogeneous quadratic polynomials
$\hat{g}_{s+1},\dots, \hat{g}_r$ among them are interreduced in step~(4).
Every polynomial~$h_i$ in~$\hat{I}$ which satisfies $\Lin_{\M}(h_i) = z_i$ is of the
form $h_i = \hat{g}_i + k_i$ where $k_i \in \langle \hat{g}_{s+1},\dots, \hat{g}_r\rangle_K$.
In order to have a polynomial of this kind for which~$z_i$ can be assigned the smallest
available weight, it is necessary that its homogeneous component of degree two is zero.
The check whether such polynomials exist is done in steps~(5) and~(6) where the
quadratic parts of the polynomials $\hat{f}_j$ are $K$-linearly reduced against
the homogeneous polynomials of degree two. 

If the quadratic part of no polynomial $\hat{f}_j$ reduces to zero, our attempt to
assign weights fails and step~(6) gives the correct answer. Otherwise, we assign those
indeterminates, namely the ones in~$\widetilde{Z}$, the next available weight~$d$.
In the loop of steps (7)--(9) we repeat this procedure. Setting the terms divisible only
by indeterminates in~$\widetilde{Z}$ to zero and increasing the lowest available weight
to $2d+1$ makes sure that any remaining indeterminate will get a bigger weight than any
of the deleted terms. Once again, interreducing the quadratic polynomials and then
reducing the polynomials $z_i + \hat{p}_i$ against these quadratic polynomials 
decides whether there exists an indeterminate~$z_i$ that can get the next lowest 
available weight. 

If this loop reaches $\widehat{Z} = \emptyset$, we have assigned
weights to all indeterminates in~$Z$ such that for $i=1,\dots,s$ there exist polynomials
$\tilde{g}_i = z_i + p_i$ in~$I$ for which~$z_i$ is the unique term of highest weight 
in~$\tilde{g}_i$ and $p_i$ is homogeneous of degree two. 
In particular, we can renumber these polynomials and the indeterminates in~$Z$
such that $z_1 <_\sigma \cdots <_\sigma z_s$ and $p_i \in K[Y][z_1,\dots,z_{i-1}]_2$
for $i=1,\dots,s$. Starting with $\tilde{g}_1 = z_1 + p_1$, we can then substitute~$z_i$
in all polynomials $\tilde{g}_j$, where $j>i$, by $-p_i$ and get the reduced
$\sigma$-Gr\"obner basis of $\langle \tilde{g}_1,\dots, \tilde{g}_s \rangle$
we were looking for. Altogether, Proposition~2.6 in~\cite{KLR2} implies that
$(g_1,\dots,g_s)$ is coherently $Z$-separating, as we wanted to show.
\end{proof}

Several variants of Algorithm~\ref{alg-checkZsepI} are possible, 
as the following remark shows.

\medskip
\begin{remark}\label{rem-CheckZvariants}
Suppose we are in the setting of Algorithm~\ref{alg-checkZsepI}.
\begin{enumerate}
\item[(a)] The hypothesis that the polynomials $f_1,\dots,f_s \in I$ have degree
two is not strictly necessary. Using some care, one can write down a version
of the algorithm which replaces the quadratic parts $\hat{p}_i,\, q_j$ by the
non-linear parts of the respective polynomials. Since this more general situation
is not relevant to this paper, we leave the details to the interested readers.

\item[(b)] For certain types of ideals, some of the polynomials $\hat{f}_1,\dots,
\hat{f}_r$ in step~(2) of Algorithm~\ref{alg-checkZsepI} are already of the form
$c\, z_i$ with $c\in K\setminus \{0\}$ and $z_i\in \widehat{Z}$. In this case we can
include an execution of step~(6) and replace step~(3) by the following steps.
\begin{itemize}
\item[(3a)] If $\dim_K( \langle \hat{\ell}_1, \dots, \hat{\ell_r} \rangle_K) < s$, 
return {\tt "No Weights Found"} and stop. 

\item[(3b)] Let $\widetilde{Z}$ be the set of all $z_i\in \widehat{Z}$
such that $\Supp(\hat{f}_i) = \{z_i\}$. If we have $\widetilde{Z} \ne \emptyset$ then
perform the following step (3c) and afterwards return to step~(2). Otherwise, continue
with step~(4).

\item[(3c)] Let $w_k=d$ for each indeterminate $z_i=x_k$ in~$\widetilde{Z}$
and remove from~$\widehat{Z}$ all indeterminates from~$\widetilde{Z}$.
Let $\hat{s} = \# \widehat{Z}$, renumber the polynomials $\hat{f}_j$ such that
the remaining ones are $\hat{f}_1, \dots,\hat{f}_{\hat{s}}$, and denote $2d+1$ by~$d$.
\end{itemize}
The advantage of this version of the algorithm is that some indeterminates in~$Z$
which can be eliminated may be detected before the linear reductions in steps (3), (4), and~(5). 
On the other hand, the weights assigned to the remaining indeterminates in~$\widehat{Z}$ 
will tend to be higher.

\item[(c)] If we replace step~(4) by the following sequence of steps, the algorithm
will frequently find more homogeneous quadratic polynomials in the ideal. In this case
it will be able to reduce the polynomials $\hat{f}_i$ better in step~(5) and discover more
indeterminates in~$\widetilde{Z}$ in step~(6).
\begin{itemize}
\item[(4a)] Form all products $y_i \widehat{f}_j$ where $y_i\in Y$ and $j = 1,\dots,\hat{s}$.
Let $S$ be the set of all terms of degree three in the support of one of these
products.

\item[(4b)] Compute a basis $q'_1,\dots,q'_v$ of the set of all $K$-linear 
combinations of the products $y_i \widehat{f}_j$ which do not have any term of~$S$ 
in their support.

\item[(4c)] Interreduce the polynomials $q_1,\dots,q_u,\, q'_1,\dots,q'_v$ $K$-lin\-ear\-ly
such that they have distinct leading terms with respect to~{\tt Lex}. 
Remove all zero polynomials and denote the resulting polynomials by $q_1,\dots,q_u$
again, redefining~$u$ if necessary.
\end{itemize}
The disadvantage of this version is that we may have to compute the kernel of a huge
matrix in step~(4b).

\item[(d)] If we mimic all linear reductions performed in steps~(3), (4), and~(5)
on the initial generators $f_1,\dots,f_r$, and if we renumber everything such that
the indeterminates $z_i$ are given weights in the order $z_1, z_2,\dots,z_s$, then
we obtain a set of polynomials $f'_1,\dots,f'_s$ which are {\bf $Z$-triangular} in the
following sense: $f'_1 = z_1 - h'_1$, $\dots$, $f'_s = z_s - h'_s$ where
$h'_1\in K[Y]_{\le 2}$, $h'_2 \in K[Y][z_1]_{\le 2}$, $\dots$, 
$h'_s \in K[Y][z_1,\dots,z_{s-1}]_{\le 2}$.

In particular, by substituting $z_1 \mapsto h'_1$ in $f'_2,\dots,f'_s$,
we get polynomials of the form $f''_2 = z_2 - h''_2$, etc., 
where $h''_2 \in K[Y]_{\le 3}$, and so on.
Next we continue by substituting $z_2 \mapsto h''_2$ in $f''_3,\dots,f''_s$, and so on.
In this way we finally get polynomials $g_1 = z_1 - \tilde{h}_1$, $\dots$,
$g_s = z_s - \tilde{h}_s$ with $\tilde{h}_1, \dots, \tilde{h}_s \in K[Y]$.
Thus the tuple $(g_1,\dots,g_s)$ is coherently $Z$-separating and defines a
$Z$-separating re-embedding of~$I$.

Notice that the calculation of $f'_1,\dots,f'_s$ involves only linear algebra operations
in the vector space $\langle f_1,\dots,f_r \rangle_K$ and that the substitutions
yielding $g_1,\dots,g_s$ from $f'_1,\dots,f'_s$ involve no reduction steps at all.
Therefore, if Algorithm~\ref{alg-checkZsepI} succeeds, a coherently $Z$-separating
tuple $(g_1,\dots,g_s)$ can be computed efficiently.
\end{enumerate}
\end{remark}

To finish this section, let us apply Algorithm~\ref{alg-checkZsepI} 
in a concrete example.

\begin{example}\label{ex-checkZsep}
Let $P=\QQ[x_1,...,x_{10}]$, and let~$I$ be the ideal of~$P$ generated by
\begin{align*} 
f_1 &=  x_1^2 \, {+} x_1 x_5 \,{+} x_4 x_6 \,{-} x_7 +x_{10},\; 
f_2 = x_3 x_6 \,{+} 2 x_6 x_7 \,{-} x_4 x_9 \,{-} x_5 x_9 - x_1,\\
f_3 &= x_{4}^2 -2x_6x_7 +x_5x_9 +x_{6} x_{9} -x_{5},\; 
f_4 = x_{6} x_{7} +x_{9}^2  -x_{10}^2  +x_{6} -x_{8},\\
f_5 &=  -x_{6} x_{7} +x_{4} x_{9} +x_{5} x_{9} -x_{2} +x_{6},\; 
f_6 =  x_{6} x_{9} +x_{4} x_{10} +x_{2} -x_{8}
\end{align*}
The linear parts of these polynomials are 
$$
-x_{7} +x_{10},\quad -x_1,\quad -x_5,\quad x_6-x_8,\quad -x_2+x_6,\quad x_2-x_8 
$$
and so the $K$-vector space spanned by these linear parts has dimension~5.

\begin{enumerate}
\item[(a)] Let us choose $Z=(x_1, x_2, x_5, x_7, x_8)$, so that the complementary tuple is
$Y =X\setminus Z =(x_3, x_4, x_6, x_9, x_{10})$, and apply 
Algorithm~\ref{alg-checkZsepI} to check whether~$I$ is $Z$-separating.
\begin{itemize}
\item[(1)] First we assign weight zero to the indeterminates in the tuple $Y$.

\item[(2)] Deleting all terms in $\Supp(f_i) \cap K[Y]$ yields 
\begin{align*} 
\hat{f}_1 &\;=\; x_{1}^2  +x_{1} x_{5} -x_{7}, &
\hat{f}_2 &\;=\; 2x_{6} x_{7} -x_{5} x_{9}-x_{1},\\
\hat{f}_3 &\;=\; -2x_{6} x_{7} +x_{5} x_{9}-x_{5}, &
\hat{f}_4 &\;=\; x_{6} x_{7} -x_{8},\\
\hat{f}_5 &\;=\;  -x_{6} x_{7} +x_{5} x_{9} -x_{2}, &
\hat{f}_6 &\;=\;  x_{2} -x_{8}
\end{align*}

\item[(3)] The result of these interreductions is
\begin{align*} 
\hat{f}_1 &\;=\; x_1 -2x_{6} x_{7} +x_{5} x_{9}, &
\hat{f}_2 &\;=\; x_2 +x_{6} x_{7} -x_{5} x_{9},\\
\hat{f}_3 &\;=\; x_5 +2x_{6} x_{7} -x_{5} x_{9},&
\hat{f}_4 &\;=\; x_7 -x_{1}^2  -x_{1} x_{5},\\
\hat{f}_5 &\;=\; x_8 -x_6x_7,&
q_1 &\;=\; -2x_{6} x_{7} +x_{5} x_{9}
\end{align*}

\item[(5)] We obtain  $\tilde{f}_1 =x_1$, $\tilde{f}_2 = x_2 -\tfrac{1}{2}x_{5} x_{9}$,
$\tilde{f}_3 = x_5$, $\tilde{f}_4 = x_7 -x_{1}^2  -x_{1} x_{5}$, and
$\tilde{f}_5 = x_8 -\tfrac{1}{2}x_{5} x_{9}$.

\item[(6)] Hence we have $\widetilde{Z} = \{x_1, x_5\}$ and assign weight~1 to~$x_1$ and~$x_5$.
We set $\widehat{Z} = (x_2, x_7, x_8)$ and $d=3$.

\item[(8)] Here we get $\tilde{f}_1 = x_2$, $\tilde{f}_2 = x_7$, $\tilde{f}_3 = x_8$,
and $q_1 = -2x_6x_7$.

\item[(5)-(7)] Now  $\tilde{f}_1 = x_2$, $\tilde{f}_2 = x_7$, and $\tilde{f}_3 = x_8$
allows us to set the weights of~$x_2$, $x_7$, and $x_8$ to~3 and return the tuple
$$
W=(1,\  3,\  0,\  0,\  1,\  0,\  3,\  3,\  0,\  0)
$$
\end{itemize}
Hence the ideal~$I$ is $Z$-separating and every term ordering compatible with the
grading given by~$W$ is $Z$-separating for~$I$.

\item[(b)] If we apply Algorithm~\ref{alg-checkZsepI} to 
$Z=(x_1,x_2,x_5,x_6,x_7)$, we get $\widetilde{Z} = \{ x_5\}$ in the first
iteration of step~(5) and 
\begin{align*} 
\tilde{f}_1 &\;=\; x_{1} -x_{3}x_{6} - x_{6}x_{9}, &
\tilde{f}_2 &\;=\; x_{2} + x_{6}x_{9},\\
\tilde{f}_3 &\;=\; x_{6} +\tfrac{1}{2}x_{6}x_{9}, &
\tilde{f}_4 &\;=\; x_{7} -x_{1}^2  -x_{1}x_{5} -x_{4}x_{6}
\end{align*}
in the second iteration of step~(5). Hence the algorithm returns
``\texttt{No Weights Found}'' and stops.
\end{enumerate}
In fact, using Algorithm~\ref{alg-checkZsepI} we may check that
there are no other tuples~$Z$ of length~5 such that~$I$ is $Z$-separating.
\end{example}

\bigskip\bigbreak
%
%

\section{Separating Re-embeddings via Gr\"obner Fans of Linear Ideals}
\label{Separating Re-embeddings via Groebner Fans of Linear Ideals}

In this section we introduce a novel method for finding optimal
$Z$-separating re-embeddings which improves on the approach in~\cite{KLR3}.
Let $P=K[x_1,\dots,x_n]$ be a polynomial ring over a field~$K$,
let $\M = \langle x_1,\dots, x_n\rangle$, and let $I\subseteq \M$ 
be an ideal of~$P$.
In the first subsection we show that the knowledge of $\GFan(\langle\Lin_\M(I)\rangle)$, 
i.e., the Gr\"obner fan of the linear part of~$I$, can be used to find $Z$-separating
re-embedding of~$I$, and in the second subsection we study how one can
compute Gr\"obner fans of homogeneous linear ideals.

\medskip
\subsection*{2.A.\ Finding \textit{Z}-Separating Re-embeddings}$\mathstrut$

Given a tuple of distinct indeterminates~$Z$ in~$X$, the methods explained in
Subsection~1.B allow us to check whether the ideal~$I$ is $Z$-separating.
This raises the question how to find such tuples~$Z$. 
In~\cite{KLR2} and~\cite{KLR3} we gave some answers
to this question which use the computation of the Gr\"obner fan of~$I$. 
Unfortunately, this computation may be infeasible for large examples. 
However, a way out is suggested by Proposition~\ref{prop-cohSepPoly}.c. 
Indeed, our next algorithm says that we can get away with knowing
$\GFan(\langle\Lin_\M(I)\rangle)$. 
As mentioned previously, if we have the equality
$I = \langle f_1,\dots,f_r\rangle$ with $f_i\in \M$ then 
we can calculate the linear part of~$I$ via
$\Lin_{\M}(I) = \langle \Lin_{\M}(f_1),\dots,\Lin_{\M}(f_r) \rangle_K$
(see~\cite{KLR2}, Proposition 1.9).

\begin{algorithm}{\bf ($Z$-Separating Re-embeddings 
via ${\boldsymbol \GFan(\langle \Lin_\M(I) \rangle)}$}
\label{alg-Emb}\\
Let $I\subseteq \M$ be an ideal of~$P$, and let $s \le \dim_K( \Lin_\M(I) )$.
Consider the following sequence of instructions.
\begin{enumerate}
\item[(1)] Compute $\GFan(\langle \Lin_\M(I) \rangle)$.

\item[(2)] Form the set~$S$ of all tuples $Z=(z_1,\dots,z_s)$ such that
there is marked reduced Gr\"obner basis~$\overline{G}$ in 
$\GFan(\langle \Lin_\M(I) \rangle)$ for which $z_1,\dots,z_s$ are among the marked terms. 

\item[(3)] If $S = \emptyset$, return 
{\tt "No re-embedding found"}. 
While $S\ne \emptyset$, perform the following steps.

\item[(4)] Choose a tuple $Z = (z_1,\dots,z_s) \in S$ and remove it from~$S$.

\item[(5)] Using one of the methods suggested in Subsection~1.B, check whether the ideal~$I$ is 
$Z$-separating. If it is, return~$Z$ and stop. Otherwise, continue with the next iteration.
\end{enumerate}
This is an algorithm which, if successful, finds a tuple of distinct indeterminates
$Z=(z_1,\dots,z_s)$ in~$X$ such that~$I$ is $Z$-separating. 

Moreover, if $s = \dim_K( \Lin_\M(I) )$ and the algorithm is successful then
the output tuple~$Z$ defines an optimal re-embedding of~$I$.
\end{algorithm}

\begin{proof}
Every tuple~$Z$ such that there exists an $Z$-separating re-em\-bed\-ding
of~$I$ is contained in the tuple of leading terms of a marked reduced Gr\"obner basis
of~$\Lin_{\M}(I)$ by Proposition~\ref{prop-cohSepPoly}.c.
The set of all possible such tuples~$Z$ is computed in steps~(1) and~(2).
If the loop in steps (4)-(5) finds a tuple~$Z$ such that~$I$ is $Z$-separating,
we are done.

In addition, if $s = \dim_K( \Lin_\M(I) )$ and the algorithm is successful,
then~\cite{KLR2}, Corollary 4.2, shows that the $Z$-separating re-embedding 
of~$I$ is optimal.
\end{proof}

If we are looking for optimal re-embeddings and
use the method of Proposition~\ref{prop-CheckZ} to perform step~(5), 
Algorithm~\ref{alg-Emb} is able to certify that no optimal $Z$-separating 
re-embedding of~$I$ exists. However, it may have to compute some huge Gr\"obner bases. 
Moreover, the next remark points out some further limitations.

\begin{remark}\label{rem-onlysufficient}
Notice that Algorithm~\ref{alg-Emb} provides only a 
sufficient condition for detecting optimal re-embeddings of~$I$. 
On one side, it can happen that an optimal re-embedding is obtained using a subset 
of generators of a leading term ideal of $\Lin_\M(I)$ (see~\cite{KLR3}, Example~6.6).
On the other side, it can happen that an optimal re-embedding cannot 
be achieved by a separating re-embedding, as it was already shown in~\cite{KLR2}, Example 3.7.
The following Example~\ref{ex-isotoK[x]} provides an even more intriguing case of a somehow 
unexpected isomorphism.
\end{remark}

\begin{example}\label{ex-isotoK[x]} 
Let $P=\QQ[x,y,z]$, and let $\theta:\; P \To P$ be defined by 
$\theta(x) = z -y^2$, $\theta(y) = x -(z-y^2)^2$, and $\theta(z) =  x +y - (z-y^2)$.
It is easy to check that $\theta(y +x^2) = x$, $\theta(z +x -y -x^2) = y$, and
$\theta(x + (z +x -y -x^2)^2) = z$. Therefore the map~$\theta$ is surjective, and
for dimensional reasons it follows that~$\theta$ is an isomorphism. Its inverse
$\theta^{-1}$ satisfies $\theta^{-1}(x) = y +x^2$, $\theta^{-1}(y) = x-y+z-x^2$, 
and $\theta^{-1}(z) =  x +(x -y +z -x^2)^2$.

Letting $J = \langle y-x^2,\, z+x^4\rangle$, we  construct an isomorphism 
$\lambda: P/J \To \mathbb{Q}[x]$ such that $\lambda(\bar{x}) = x$, 
$\lambda(\bar{y}) = x^2$, and $\lambda(\bar{z}) = -x^4$. 
Next, to get an isomorphism $\bar{\theta}: P/J'  \to P/J$
induced by~$\theta$, we transform the generators of~$J$ 
via $\theta^{-1}$ and get 
\begin{align*}
\theta^{-1}(x^2 -y) \;=\; &x^4 +2x^2y +x^2 +y^2 -x +y -z \\
\theta^{-1}(x^4+z)  \;=\; &x^8 +4x^6y +6x^4y^2 +4x^2y^3 +x^4 +y^4 -2x^3 +2x^2y \\
& -2x^2z +x^2 -2xy +y^2 +2xz -2yz +z^2 +x
\end{align*}
Consequently, we have $J' = \langle \theta^{-1}(x^2-y),\; \theta^{-1}(x^4+z) \rangle$.

Finally, we notice that using the first generator of~$J'$ we can eliminate~$z$. 
We get the isomorphism $\rho: \QQ[x,y]/\langle F\rangle\To \QQ[x,y,z]/J'$,
where 
$$
F = 2x^8 +8x^6y +12x^4y^2 +8x^2y^3 +2y^4 +x 
$$
Letting $\alpha = \lambda \circ \bar{\theta}  \circ \rho$,
we obtain $\alpha : \QQ[x,y]/\langle F\rangle \To \QQ[x]$ which is an isomorphism
of $\mathbb{Q}$-algebras. 
It turns out that $\alpha(\bar{x}) = -2x^4$, $\alpha(\bar{y}) = \bar{x} -4\bar{x}^8$, 
whence $\alpha(\bar{x}^2 + \bar{y}) = x$, and we have $\Ker(\alpha) =\langle F \rangle$.
It is clear that the existence of such an isomorphism cannot be detected using
the methods developed in this paper.
\end{example}

The following example shows Algorithm~\ref{alg-Emb} at work.

\begin{example}\label{ex-searchingforoptembed}
Let $P=\QQ[x,y,z,w]$, let $F =(f_1, f_2, f_3)$, where $f_1 = x -y -w^2$, 
$f_2=x + y -z^2$, and $f_3 = z + w +z^3$, and let $I = \langle f_1,f_2,f_3\rangle$. 
\begin{enumerate}
\item[(1)] We obtain $\Lin_\M(I) = \langle z +w,\, x, \, y\rangle_K$ and the 
methods explained below return the two marked reduced 
Gr\"obner bases $\{ (x,x),\, (y,y),\, (z,z+w)\}$ and $\{ (x,x),\, (y,y),\, (w,w+z)\}$.

\item[(2)] We get $S =\{ (x, y, z),\, (x, y, w) \}$.

\item[(4)] We pick $Z = (x, y, z)$  and delete it from~$S$. 

\item[(5)] We construct an elimination ordering~$\sigma$ for~$Z$ and find that 
the minimal set of generators of~$\Lin_{\M}(\LT_\sigma(I))$ is $L = \{x\}$.
Therefore  $L\ne \{x,y,w\}$ and continue with the next iteration.

\item[(4)] Next we let $Z=(x,y,w)$ and let $S=\emptyset$.

\item[(5)] We construct an elimination ordering~$\sigma$ for~$Z$ 
and compute the minimal set of generators~$L$ of $\Lin_{\M}(\LT_\sigma(I))$. 
Since $L=\{x,y,w\}$, we return $Z = (x, y, w)$ and stop.
\end{enumerate}

Notice that we could also have used the heuristics in part~(b) of 
Remark~\ref{rem-randomFforZ} to execute step~(5) and check whether~$F$ is $Z$-separating
using an LPF-solver. In both cases we would have received the answer that this
is not the case, and we would have missed the optimal $Z$-separating re-embedding
of~$I$ for $Z=(x,y,w)$. To get polynomials which are $Z$-separating, it suffices to 
replace $f_2$ with $f_2' = f_2-f_1$. 

To find the actual polynomials defining the optimal re-embedding, 
we need to compute the reduced $\sigma$-Gr\"obner basis of~$I$. It is
$$
( x -\tfrac{1}{2} z^6 -z^4 -z^2, \;  y + \tfrac{1}{2} z^6 +z^4,\;  w +z^3 +z)
$$
This tuple gives rise to a $\QQ$-algebra isomorphism $P/I\cong \QQ[z]$
via $x \mapsto \tfrac{1}{2} z^6 +z^4 +z^2$, $y \mapsto -\tfrac{1}{2} z^6 -z^4$, and
$w \mapsto -z^3 -z$.
\end{example}

\medskip
\subsection*{2.B.\ The Gr\"obner Fan of a Homogeneous Linear Ideal}$\mathstrut$

In the remainder of this subsection we study ways how to perform step~(1) of 
Algorithm~\ref{alg-Emb}. Therefore we now let~$I_L$ be an ideal in~$P$
generated by linear forms, and our goal is to compute 
the Gr\"obner fan $\GFan(I_L)$, i.e., the set of all marked reduced Gr\"obner 
bases of~$I_L$.

More precisely, let $L=(\ell_1, \dots, \ell_r)$ be a tuple of
linear forms in~$P$, and let $I_L = \langle L\rangle$. 
For $i=1,\dots,r$, we write $\ell_i = a_{i1} x_1 + \cdots + a_{in} x_n$
with $a_{ij} \in K$. Then the matrix $A=(a_{ij}) \in \Mat_{r,n}(K)$ is 
called the {\bf coefficient matrix} of~$L$.
To compute $\GFan(I_L)$, the following proposition will be useful.

\begin{proposition}\label{prop-charLinZ}
Let $L=(\ell_1, \dots, \ell_r)$ be a tuple of $K$-linearly independent linear forms in~$P$,
let $I_L = \langle  L\rangle$, and let $A=(a_{ij})$ be the  coefficient matrix of~$L$.
Moreover, let $s\le r$, let $Z=(z_1,\dots,z_s)$ be a tuple of distinct indeterminates in $X=(x_1,\dots,x_n)$, 
and let $Y = X\setminus Z$. Then the following conditions are equivalent.
\begin{enumerate}
\item[(a)] There exists a term ordering~$\sigma$ such that 
$\langle Z\rangle \subseteq \LT_\sigma(I_L)$.

\item[(b)] There exists a term ordering~$\sigma$ such that
the residue classes of the elements of~$Y$ generate the $K$-algebra $P/\LT_\sigma(I_L)$.

\item[(c)] The residue classes of the elements of~$Y$ generate the 
$K$-algebra $P/I_L$.

\item[(d)] The residue classes of the elements of~$Y$ generate the
$K$-vector space $P_1 / \langle L\rangle_K$.

\item[(e)] Let $i_1,\dots,i_s\in \{1,\dots,n\}$ be the indices such that
$z_j =  x_{i_j}$ for $j=1,\dots,s$. Then the columns $i_1,\dots,i_s$ of~$A$
are linearly independent.
\end{enumerate}
\end{proposition}

\begin{proof}
To show (a)$\Rightarrow$(b), it suffices to note that $\langle Z \rangle 
\subseteq\LT_\sigma(I_L)$ implies that the canonical map $K[Y] \cong P/\langle Z\rangle
\longrightarrow P/\LT_\sigma(I_L)$ is surjective. 
The implication (b)$\Rightarrow$(c) follows 
from Macaulay's Basis Theorem (see~\cite{KR1}, Theorem 1.5.7).

To prove (c)$\Rightarrow$(d), we note that~$I_L$
is generated by linear forms, and therefore $P/I_L$ is isomorphic to a polynomial ring.
Thus the residue classes of a tuple of indeterminates~$Y$ are a $K$-algebra
system of generators of the ring $P/I_L$ if and only if they are a system of generators
of the $K$-vector space given by its homogeneous component $P_1/(I_L)_1$ of degree one, 
where $(I_L)_1 = \langle L\rangle_K$.
The assumption in~(d) implies that the indeterminates in~$Z$ can be expressed as linear
combinations of the indeterminates in~$Y$. Consequently, any elimination ordering~$\sigma$
for~$Z$ satisfies $\langle Z\rangle \subseteq \LT_\sigma(I_L)$. This proves (d)$\Rightarrow$(a).

Finally, we show that~(d) and~(e) are equivalent. Notice that both conditions
imply $s\le r$. The tuple $Y = X \setminus Z$
is a system of generators of the vector space $P_1 / \langle L \rangle_K$ 
if and only if $\{\ell_1,\dots,\ell_r\}$
together with~$Y$ is a system of generators of~$P_1$. This means that, if we extend~$A$ with
$n-s$ rows that are unit vectors having their non-zero entries at the positions
of the indeterminates in~$Y$, the resulting matrix~$\overline{A}$ of size $(r+n-s)\times n$ 
has the maximal rank~$n$. Now we consider the matrix $A_{(i_1,\dots,i_s)}$ 
consisting of columns $i_1,\dots,i_s$ of~$A$. By renumbering the 
indeterminates, we may assume that $i_1=1$, $\dots$, $i_s=s$, and hence
that the extended matrix is upper block triangular of the form
$$
\overline{A} \;=\; \begin{pmatrix}
A_{(i_1,\dots,i_s)} & \ast \\
0 & I_{n-s}
\end{pmatrix}
$$
where $I_{n-s}$ is the identity matrix of size $n-s$.
Now it is clear that the rows of~$\overline{A}$ generate $K^n$ if and
only if the rows of $A_{(i_1,\dots,i_s)}$ generate $K^s$, and this is equivalent
to $A_{(i_1,\dots,i_s)}$ having maximal rank~$s$. This concludes the proof
of the proposition.
\end{proof}

As a special case, we get the following characterization of tuples~$Z$
which are leading term tuples of a marked reduced Gr\"obner basis in $\GFan(I_L)$.

\begin{corollary}\label{cor-CharLinZ}
In the setting of the proposition, assume that $s=r$. Then the following
conditions are equivalent.
\begin{enumerate}
\item[(a)] There exists a term ordering~$\sigma$ such that 
$\langle Z\rangle = \LT_\sigma(I_L)$.

\item[(b)] There exists a term ordering~$\sigma$ such that
the residue classes of the elements of~$Y$ generate the $K$-algebra 
$P/\LT_\sigma(I_L)$ minimally.

\item[(c)] The residue classes of the elements of~$Y$ generate the 
$K$-algebra $P/I_L$ minimally.

\item[(d)] The residue classes of the elements of~$Y$ are a $K$-basis 
of $P_1 / \langle L\rangle_K$.

\item[(e)] Let $i_1,\dots,i_s\in \{1,\dots,n\}$ be the indices such that
$z_j =  x_{i_j}$ for $j=1,\dots,s$. Then the columns $i_1,\dots,i_s$ of~$A$
form an invertible matrix of size $s\times s$.
\end{enumerate}
\end{corollary}

Our next goal is to construct a bijection between the Gr\"obner fan
of~$I_L$ and the non-zero maximal minors of~$A$. The following
terminology will prove useful.

\begin{definition}\label{def-genGFan}
Let $J$ be an ideal in~$P$. 
\begin{enumerate}
\item[(a)] For a marked reduced Gr\"obner basis
$$
G = \{ (\LT_\sigma(g_1), g_1), \dots, (\LT_\sigma(g_r), g_r)\}
$$ 
of~$J$,
we call $\LT_\sigma(G) = \{ \LT_\sigma(g_1),\dots, \LT_\sigma(g_r)\}$
the {\bf leading term set} of~$G$.

\item[(b)] The set $\LTGFan(J)$ of all leading term sets of marked
reduced Gr\"obner bases in $\GFan(J)$ is called the {\bf leading term Gr\"obner fan}
of~$J$. 
\end{enumerate}
\end{definition}

The following lemma collects some preparatory results.

\begin{lemma}\label{lem-easyassociated}
Let $L=(\ell_1, \dots, \ell_s)$ be a tuple of  $K$-linearly independent linear forms in~$P$,
let $I_L = \langle  L\rangle$, and let $A=(a_{ij})$ be the  coefficient matrix of~$L$. 
Moreover, let $L' = (\ell'_1, \dots, \ell'_r)$  be a further tuple 
of linear forms in~$I_L$, and let $A' \in\Mat_{r,n}(K)$ be its coefficient matrix.
\begin{enumerate}
\item[(a)]  The tuple~$L'$ is a minimal system of generators of~$I_L$ if and only if
$r=s$ and there exists a matrix $U \in \GL_s(K)$ such that $L' = U\cdot L$, 
or equivalently  such that $A' = U\cdot A$.

\item[(b)] A set of pairs $\{ (x_{i_1}, \ell'_1) \dots, (x_{i_r}, \ell'_r) \}$,
where $1\le i_1 < \cdots < i_r \le n$, 
is a marked reduced Gr\"obner basis of~$I_L$ iff $r=s$, the matrix $A_{(i_1,\dots,i_s)}$
is invertible, and $L' = A_{(i_1, \dots, i_s)}^{-1}\cdot L$.
\end{enumerate}
\end{lemma}

\begin{proof} 
Claim~(a) follows from the fact that every tuple of 
minimal generators of~$I_L$ is also a basis of the $K$-vector space~$(I_L)_1$.

To prove~(b) we observe that a minimal Gr\"obner basis of a linear homogeneous ideal 
is also a minimal set of generators of~$I_L$. This yields $r=s$. Moreover, it is reduced 
if and only if the submatrix $A'_{(i_1, \dots, i_s)}$ of~$A'$
is the identity matrix, and hence the conclusion follows from~(a).
\end{proof}

Now we are ready to present the key result for computing $\GFan(I_L)$.

\begin{theorem}{\bf (The Gr\"obner Fan of a Homogeneous Linear 
Ideal)}\label{thm-bijection}\\
Let $L=(\ell_1, \dots, \ell_s)$ be a tuple of $K$-linearly independent linear forms in~$P$,
let $I_L = \langle L\rangle$, and let $A=(a_{ij})$ be the  coefficient matrix of~$L$.  
Furthermore, let $M$ be the set of tuples $(i_1,\dots,i_s)$
such that $1\le i_1 < \cdots < i_s \le n$ and such that the corresponding
maximal minor of the matrix $A=(a_{ij})$ is non-zero.
\begin{enumerate}
\item[(a)] The map $\phi:\; \LTGFan(I_L) \longrightarrow M$ given by
$\phi(Z) = (i_1,\dots,i_s)$ for $Z = \{ x_{i_1},\dots,x_{i_s} \} \in \LTGFan(I_L)$
and $i_1 < \cdots < i_s$ is well-defined and bijective.

\item[(b)] The  map $\psi:\; \GFan(I_L) \longrightarrow M$ given by
$\psi(G) = \phi(\LT_\sigma(G))$ for every $G\in \GFan(I)$ is well-defined and bijective.

\end{enumerate}
\end{theorem}

\begin{proof}
First we prove~(a).
To begin with, let us check that~$\phi$ is well-defined. For an element
$G\in \GFan(I_L)$, the tuple $Z=\LT_\sigma(G)$ satisfies $\langle Z\rangle
=\LT_\sigma(I_L)$. Let $Z=\{x_{i_1},\dots,x_{i_s}\}$ with $1\le i_1< \cdots < i_s \le n$. 
Then Corollary~\ref{cor-CharLinZ}.e shows $\det(A_{(i_1,\dots,i_s)}) \ne 0$,
and therefore~$\phi$ is well-defined.

Since the map~$\phi$ is clearly injective, we still need to show that it is surjective.
Given $(i_1,\dots,i_s)\in M$, part~(b) of the lemma implies that
$L' = A_{(i_1,\dots,i_s)}^{-1} \cdot L$ is a reduced Gr\"obner basis of~$I_L$
with leading term set $Z=\{ x_{i_1}, \dots, x_{i_s} \}$. 
Therefore $\phi(Z) = (i_1,\dots,i_s)$ proves the desired surjectivity.

To show~(b), we note that the map~$\psi$ is clearly well-defined and injective.
By the definition of $\LTGFan(I_L)$, it is also surjective. 
\end{proof}

Based on the preceding results, there exist several ways of computing
the Gr\"obner fan of an ideal generated by linear forms.

\begin{remark}{\bf (Computation of $\GFan(I_L)$)}\label{rem-compGFan}\\
Let $I_L=\langle \ell_1,\dots,\ell_s \rangle$ be an ideal in~$P$ generated
by linearly independent linear forms $\ell_1,\dots,\ell_s \in P_1$ as above.
\begin{enumerate}
\item[(a)] Using the preceding proposition, we can calculate $\GFan(I_L)$
by computing all maximal minors of the matrix~$A$. If we have $\det(A_{(i_1,\dots,i_s)}) \ne 0$
for some $1\le i_1 < \cdots < i_s \le n$, the tuple $L' = A_{(i_1,\dots,i_s)}^{-1}
\cdot L$ is the reduced Gr\"obner basis of~$I_L$, 
and $Z = \{ x_{i_1}, \dots, x_{i_s} \}$
is its leading term set with respect to any elimination ordering for~$Z$.

\item[(b)] Another way of viewing the task to compute $\GFan(I_L)$ is
obtained by applying Corollary~\ref{cor-CharLinZ}. The complements
of the leading term sets of reduced Gr\"obner bases of~$I_L$ correspond uniquely
to sets of terms whose residue classes form a $K$-basis of $P_1 / (I_L)_1$.
All sets of indeterminates whose residue classes are linearly independent 
in~$P_1/(I_L)_1$ are the independent sets of a linear matroid, and maximal
such sets are the bases of the matroid. The task of computing the bases of a linear
matroid has been studied intensively, and many algorithms are known,
see for instance the reverse search technique of D.\ Avis and K. Fukuda (cf.~\cite{AF}).
\end{enumerate}
\end{remark}

Let us recall the definition of a restricted Gr\"obner fan from~\cite{KLR3}
and introduce the corresponding notion of a restricted leading term Gr\"obner fan.

\begin{definition}\label{def-restrictedgenGFan}
Let~$J$ be an ideal in~$P$, and let $Z=(z_1,\dots,z_q)$ be a tuple of distinct indeterminates
in~$X$.
\begin{enumerate}
\item[(a)] The set of all marked reduced Gr\"obner bases
$$
\overline{G} \;=\; \{ (\LT_\sigma(g_1),\, g_1),\; \dots,\; (\LT_\sigma(g_r),\, g_r) \}
$$
such that for every $j\in \{1,\dots,q \}$ there exists an index $i_j \in \{1,\dots,r\}$
with $z_j = \LT_\sigma(g_{i_j})$ is called the {\bf $Z$-restricted Gr\"obner fan} of~$J$
and denoted by $\GFan_Z(J)$.

\item[(b)] The set $\LTGFan_Z(J)$ of all leading term sets of marked
reduced Gr\"obner bases in $\GFan_Z(J)$ is called the {\bf $Z$-restricted 
leading term Gr\"obner fan} of~$J$.  
\end{enumerate}
\end{definition}

The following corollary follows easily from Proposition~\ref{thm-bijection}.

\begin{corollary}{\bf (The Restricted GFan of a Homogeneous Linear 
Ideal)}\label{cor-restrictedGBLinear}\\
Let $L=(\ell_1, \dots, \ell_s)$ be a tuple of $K$-linearly independent linear forms in~$P$,
let $I_L = \langle  L\rangle$, and let $A=(a_{ij})$ be the  coefficient matrix of~$L$.   
Furthermore, let $Z=(z_1,\dots,z_q)$ be a tuple of distinct indeterminates
in~$X$, and let $M_Z$ be the set of all tuples $(i_1,\dots,i_s)$ such that
$1\le i_1 < \cdots < i_s \le n$, such that $\det(A_{(i_1,\dots,i_s)}) \ne 0$,
and such that $\{z_1,\dots,z_q\} \subseteq \{ x_{i_1}, \dots, x_{i_s} \}$.  
\begin{enumerate}
\item[(a)]  The map $\phi$ induces a bijection
$\phi_Z:\; \LTGFan_Z(I_L) \longrightarrow M_Z$.

\item[(b)] The map $\psi$ induces a bijection $\psi_Z:\; \GFan_Z(I_L)  
\longrightarrow  M_Z$.
\end{enumerate}
\end{corollary}

To conclude this section, let us calculate the Gr\"obner fan of an
explicit homogeneous linear ideal.

\begin{example}\label{ex-restrictedGFan}
Let $P = \QQ[x, y, z, w]$, let $\ell_1=x+y -z +4w$, $\ell_2 =x-y-z$,
and let $I_L = \langle \ell_1, \ell_2\rangle$.
The set $\{\ell_1, \ell_2\}$ is a set of minimal generators of~$I_L$ 
and its coefficient matrix is
$$
A \;=\; \begin{pmatrix} 1 & \ \ 1 & -1 & 4 \cr  1 & -1 & -1& 0     
\end{pmatrix}
$$
One $2\times2$-submatrix is singular. The others are 
$$
M_1 = \left(\begin{smallmatrix} 1 & \;\; 1  \cr  1 & -1  \end{smallmatrix}\right),\,
M_2 = \left(\begin{smallmatrix} 1 & 4  \cr  1 & 0  \end{smallmatrix}\right),\,
M_3 = \left(\begin{smallmatrix} \;\; 1 & -1 \cr -1 & -1  \end{smallmatrix}\right),\,
M_4 = \left(\begin{smallmatrix} \;\; 1 &  4 \cr -1 & 0  \end{smallmatrix}\right),\,
M_5 = \left(\begin{smallmatrix} -1 & 4 \cr -1 & 0 \end{smallmatrix}\right)
$$
Multiplying their inverses by $A$ we get the matrices
$$
\Bigl( \begin{smallmatrix} 1 & 0 & -1 & 2 \cr 0 & \mathstrut 1 &\;\; 0 & 2 
   \end{smallmatrix}\Bigr) ,\,
\Bigl( \begin{smallmatrix} 1 & -1 & -1 & 0 \cr 0 & \;\; 1/2 & \;\; 0 & 1 
   \end{smallmatrix}\Bigr) ,\,
\Bigl( \begin{smallmatrix} \;\; 0 & 1 & 0 & \;\; 2 \cr -1 & \mathstrut 0 & 1 & -2  
   \end{smallmatrix}\Bigr) ,\,
\Bigl( \begin{smallmatrix} -1 & \; 1 & \;\; 1 & 0 \cr \;\; 1/2 & \; 0 & -1/2 & 1
   \end{smallmatrix}\Bigr) ,\,
\Bigl( \begin{smallmatrix} -1 & 1 & \; 1 & 0 \cr \;\; 0 & 1/2 & \; 0 & 1 
   \end{smallmatrix}\Bigr) 
$$
They correspond to  the following marked reduced Gr\"obner bases of~$I_L$ 
which form the Gr\"obner fan of~$I_L$:
\begin{align*}
\{(x, x-z +2w), (y, y +2w) \}, & \quad \{(x, x -y - z), (w, w + \tfrac{1}{2}y)\}, \cr
\{(y, y+2w), (z, z-x-2w)\}, & \quad \{ (y, y-x+z), (w, w +\tfrac{1}{2}x -\tfrac{1}{2}z\}, \cr
\{ (z, z-x+y), (w,w +\tfrac{1}{2}y)\} \;\; &
\end{align*}
If we let $Z = (z)$, the $Z$-restricted Gr\"obner fan of~$I_L$ is 
$$ 
\bigl\{ \; \{(y, y+2w), (z, z-x-2w)\},\quad \{ (z, z-x+y), (w,w +\tfrac{1}{2}y)\} 
\; \bigr\}
$$
\end{example}

%
%

\bigskip\bigbreak
\section{The Defining Ideal of a Border Basis Scheme}
\label{The Defining Ideal of a Border Basis Scheme}

As before, we work over an arbitrary field~$K$, and we let
$P=K[x_1,\dots,x_n]$ be a polynomial ring in~$n$ indeterminates over~$K$.
The following kind of subsets of~$\mathbb{T}^n$ will play a central role.
The name ``order ideal''  comes from the fact that~$\mathbb{T}^n$ is a {\bf poset}
(i.e., a partially ordered set) with respect to divisibility, and that an order
ideal is a finite poset ideal in~$\mathbb{T}^n$ according to the usual terminology in
the theory of posets. The notions of the border, the rim and the interior of an order
ideal have the following intuitive geometric interpretations.

\begin{definition}
Let $\OO = \{t_1,\dots,t_\mu\}$ be a finite subset of~$\mathbb{T}^n$.
\begin{enumerate}
\item[(a)] The set~$\OO$ is called an {\bf order ideal} in~$\mathbb{T}^n$
if for every term~$t\in\OO$ and every term $t'\in\mathbb{T}^n$ dividing~$t$,
we have $t'\in\OO$.

\item[(b)] For an order ideal~$\OO$, we call
$$
\partial\OO \;=\; (x_1\OO \cup \cdots \cup x_n\OO) \;\setminus\; \OO
$$
the {\bf border} of~$\OO$. The border of an order ideal $\OO = \{ t_1,\dots,t_\mu\}$
will always be denoted by $\partial\OO = \{ b_1,\dots, b_\nu\}$.

\item[(c)] For an order ideal~$\OO$, the set of terms
$$
\nabla\OO = \{ t \in\OO \mid x_i\, t \in \partial\OO \hbox{\ \rm for some\ }i
\in\{1,\dots,n\}\}
$$
is called the {\bf rim} of~$\OO$.

\item[(d)] For an order ideal~$\OO$, the set of terms $\OO^\circ = \OO \setminus \nabla\OO$
is called the {\bf interior} of~$\OO$.
\end{enumerate}
\end{definition}

\begin{example}{\bf (The (2,2)-Box)}\label{ex-box22}\\
In the polynomial ring $P=K[x,y]$ over a field~$K$, 
we consider the order ideal $\OO = \{1,\, y,\, x,\, xy\}$ which we call the {\bf (2,2)-box}.
Let us illustrate the above notions as follows.
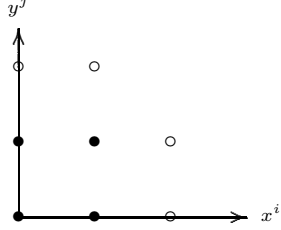
\begin{center}
\makebox{\beginpicture
\setcoordinatesystem units <1cm,1cm>
\setplotarea x from 0 to 3, y from 0 to 2.5
\axis left /
\axis bottom /
\arrow <2mm> [.2,.67] from  2.5 0  to 3 0
\arrow <2mm> [.2,.67] from  0 2  to 0 2.5
\put {$\scriptstyle x^i$} [lt] <0.5mm,0.8mm> at 3.1 0.1
\put {$\scriptstyle y^j$} [rb] <1.7mm,0.7mm> at 0 2.6
\put {$\bullet$} at 0 0
\put {$\bullet$} at 1 0
\put {$\bullet$} at 0 1
\put {$\bullet$} at 1 1

\put {$\circ$} at 0 2
\put {$\circ$} at 2 0
\put {$\circ$} at 2 1
\put {$\circ$} at 1 2
\endpicture}
\end{center}
Here the bullet points correspond to the terms in the order ideal~$\OO$
and the circles to its border $\partial\OO = \{ y^2,\, x^2,\, xy^2,\, x^2y\}$.
The rim of~$\OO$ is given by the set of terms $\nabla\OO = \{ y,\, x,\, xy\}$
just inside the border, and the interior of~$\OO$ is simply $\OO^\circ = \{ 1\}$.
\end{example}

Order ideals provide convenient $K$-bases of 0-dimensional affine $K$-algebra
in the following way.

\begin{definition}
Let $\OO=\{t_1,\dots,t_\mu\}$ be an order ideal in~$\mathbb{T}^n$ with border
$\partial\OO= \{b_1,\dots, b_\nu\}$.
\begin{enumerate}
\item[(a)] A set $G=\{g_1,\dots,g_\nu\}$ of polynomials in~$P$ of the form
$g_j = b_j - \sum_{i=1}^\mu \gamma_{ij} t_i$ with $\gamma_{ij}\in K$
is called an {\bf $\OO$-border prebasis}.

\item[(b)] An $\OO$-border prebasis~$G$ is called an {\bf $\OO$-border basis}
if the residue classes of the elements of the set~$\OO$ form a $K$-basis of the
0-dimensional affine $K$-algebra $P/\langle G\rangle$.
\end{enumerate}
\end{definition}

It is easy to see that the residue classes of the elements of~$\OO$ generate
the $K$-vector space $P/\langle G\rangle$ when~$G$ is an $\OO$-border prebasis.
For the basic theory of border bases, we refer the reader to~\cite{KR2}, Section~6.4.
Border bases allow us to parametrize large sets of 0-dimensional polynomial ideals in~$P$.
The corresponding moduli space is defined as follows.

\begin{definition}\label{def-BO}
Let $\OO=\{t_1,\dots,t_\mu\}$ be an order ideal in~$\mathbb{T}^n$ with border
$\partial\OO= \{b_1,\dots, b_\nu\}$.
\begin{enumerate} 
\item[(a)] Let $C=\{c_{ij} \mid i\in\{1,\dots,\mu\}, j\in\{1,\dots,\nu\}\}$
be a set of new indeterminates, and let $K[C]=K[c_{11},\dots,c_{\mu\nu}]$.
Then the set of polynomials $G=\{g_1,\dots,g_\nu\}$ in $K[C][x_1,\dots,x_n]$,
where  $g_j \;=\; b_j - c_{1j} t_1 - \cdots - c_{\mu j} t_\mu$
for $j=1,\dots,\nu$, is called the {\bf generic $\OO$-border prebasis}.

\item[(b)] For $r=1,\dots,n$, the matrix
$\mathcal{A}_r = (a_{ij}^{(r)}) \in \Mat_\mu(K[C])$, where
$$
a_{ij}^{(r)} \;=\; \begin{cases}
\delta_{im} & \hbox{\ \rm if\ } x_r t_j = t_m\\
c_{im} & \hbox{\ \rm if\ } x_r t_j = b_m
\end{cases}
$$
is called the $r$-th {\bf generic multiplication matrix} for~$\OO$.

\item[(c)] Consider the ideal~$I(\BO)$ in $K[C]$ which is
generated by all entries of the commutators
$\mathcal{A}_r \mathcal{A}_s - \mathcal{A}_s \mathcal{A}_r$ 
such that $1\le r<s\le n$.
Then the closed subscheme of the affine space $\mathbb{A}^{\mu\nu}_K
= \Spec(K[C])$ defined by this ideal is called the
{\bf $\OO$-border basis scheme} and denoted by~$\BO$.
The ideal~$I$ is denoted by $I(\BO)$ and called the
{\bf defining ideal} of~$\BO$.
The affine coordinate ring of~$\BO$ will be denoted by $B_\OO = K[C]/I(\BO)$.

\item[(d)] The ring homomorphism
$$
B_\OO \;\longrightarrow\; U_\mathcal{O} :=
B_\OO[x_1,\dots,x_n] / \langle g_1,\dots,g_\nu\rangle
$$
is called the {\bf universal $\OO$-border basis family}.
\end{enumerate}
\end{definition}

By~\cite{KR2}, Theorem~6.4.10, choosing $(\gamma_{ij}) \in K^{\mu\nu}$ and substituting
$c_{ij} \mapsto \gamma_{ij}$ in the generic $\OO$-border prebasis yields an $\OO$-border
basis in~$P$ if and only if the point $(\gamma_{ij})$ is a $K$-rational point of~$\BO$.
In this way the border basis scheme parametrizes all $\OO$-border bases. 
It is a very good moduli space, since~$\OO$ is a $B_\OO$-basis 
of~$U_\OO$ (cf.~\cite{KR3}, Thm.~3.4)
whence the universal $\OO$-border basis family is flat.

For the various order ideals in~$\mathbb{T}^n$ of a given length~$\mu$, the
corresponding border basis schemes form an open covering of the Hilbert scheme
$\Hilb^\mu(\mathbb{A}^n_K)$ which parametrizes all 0-dimensional subschemes 
of~$\mathbb{A}^n_K$ of length~$\mu$. Therefore border basis schemes allow us
to study large parts of this Hilbert scheme via the easily describable, 
explicitly given ideals $I(\BO)$. Let us construct the natural
system of generators of~$I(\BO)$ in an alternative way.

\begin{definition}\label{def-neighbours}
Let $\OO=\{t_1,\dots,t_\mu\}$ be an order ideal in~$\mathbb{T}^n$ with border
$\partial\OO= \{b_1,\dots, b_\nu\}$.
Let $\mathcal{A}_1,\dots,\mathcal{A}_n$ be the set of
generic multiplication matrices, and for $j=1,\dots,\mu$
let $c_j=(c_{1j},\dots,c_{\mu j})\tr$ be the $j$-th column
of~$(c_{ij})$.
\begin{enumerate}
\item[(a)] Let $j, j'\in\{1,\dots,\nu\}$ be such that
$b_j=x_\ell b_{j'}$ for some $\ell\in\{1,\dots,n\}$. Then $b_j,b_{j'}$
are called {\bf next-door neighbours} and the tuple of polynomials
$(c_j - \mathcal{A}_\ell c_{j'})\tr$
is denoted by $\ND(j, j')$.

\item[(b)] The union of all entries of the tuples $\ND(j, j')$ is called the set of
{\bf next-door generators} of~$I(\BO)$ and is denoted by
$\ND_\OO$.

\item[(c)] Let $j, j'\in\{1,\dots,\nu\}$ be such that
$b_j=x_\ell b_m$ and $b_{j'}=x_k b_m$ for some $m\in \{1,\dots,\nu\}$.
Then $b_j,b_{j'}$ are called {\bf across-the-corner neighbours}.

\item[(d)] Let $j, j'\in\{1,\dots,\nu\}$ be such that $b_j = x_\ell t_m$ and
$b_{j'} = x_k t_m$ for some $m\in\{1,\dots,\mu\}$. Then $b_j,b_{j'}$
are called {\bf across-the-rim neighbours} and the tuple of polynomials
$(\mathcal{A}_k c_j - \mathcal{A}_\ell c_{j'})\tr$
is denoted by $\AR(j, j')$.

\item[(e)] The union of all entries of the tuples $\AR(j, j')$ is called the set
of {\bf across-the-rim generators} of~$I(\BO)$ and is denoted by $\AR_\OO$.
\end{enumerate}
\end{definition}

The following example illustrates this definition. It will be reconsidered 
in the last section (see Example~\ref{ex-theLshapeCont}).

\begin{example}{\bf (The L-Shape)}\label{ex-theLshape}\\
The order ideal $\OO = \{1,\, y,\,  x,\, y^2,\, x^2\}$ in $P=\mathbb{Q}[x,y]$ 
will be  called the {\bf L-shape}.
Its border is $\partial\OO = \{ xy,\, y^3,\, xy^2,\, x^2y, x^3\}$.
\begin{figure}[ht]
	\makebox{\beginpicture
		\setcoordinatesystem units <1cm,1cm>
		\setplotarea x from 0 to 4, y from 0 to 3.5
		\axis left /
		\axis bottom /
		\arrow <2mm> [.2,.67] from  3.5 0  to 4 0
		\arrow <2mm> [.2,.67] from  0 3  to 0 3.5
		\put {$\scriptstyle x^i$} [lt] <0.5mm,0.8mm> at 4.1 0.1
		\put {$\scriptstyle y^j$} [rb] <1.7mm,0.7mm> at 0 3.6
		\put {$\bullet$} at 0 0
		\put {$\bullet$} at 1 0
		\put {$\bullet$} at 0 1
		\put {$\bullet$} at 0 2
		\put {$\bullet$} at 2 0
		\put {$\scriptstyle 1$} [lt] <-1mm,-1mm> at 0 0
		\put {$\scriptstyle t_1$} [rb] <-1.3mm,0.4mm> at 0 0
		\put {$\scriptstyle t_3$} [rb] <-1.3mm,0.4mm> at 1 0
		\put {$\scriptstyle t_2$} [rb] <-1.3mm,0mm> at 0 1
		\put {$\scriptstyle t_4$} [rb] <-1.3mm,0mm> at 0 2
		\put {$\scriptstyle t_5$} [lb] <0.8mm,0.4mm> at 2 0
		\put {$\scriptstyle b_1$} [rb] <4.3mm,0mm> at 1 1
		\put {$\scriptstyle b_2$} [lb] <0.8mm,0mm> at 0 3
		\put {$\scriptstyle b_3$} [lb] <0.8mm,0mm> at 1 2
		\put {$\scriptstyle b_4$} [lb] <0.8mm,0mm> at 2 1
		\put {$\scriptstyle b_5$} [lb] <0.8mm,0.7mm> at 3 0
		\put {$\times$} at 0 0
		\put {$\circ$} at 0 3
		\put {$\circ$} at 2 1
		\put {$\circ$} at 1 2
		\put {$\circ$} at 1 1
		\put {$\circ$} at 3 0
		\endpicture}
\end{figure}
Then the set of triples $(j, j',\ell)$ associated to next-door neighbours 
$b_j = x_\ell b_{j'}$ is $\{(3,1,2)\}, (4,1,1) \}$, 
since $b_3 = y b_1$ and $b_4 = x b_1$.
Moreover, the set of quintuples $(j, j', k, \ell, m)$ associated to 
across-the-corner neighbours $b_j=x_k b_m$ and $b_{j'}=x_\ell b_m$
is $\{(4, 3, 1, 2, 1)\}$, since $b_4 = x b_1$ and $b_3 = y b_1$.
Finally, the set of quintuples $(j, j', k, \ell, m)$ associated to
across-the-rim neighbours $b_j = x_k t_m$ and $b_{j'} = x_\ell t_m$ 
is  $\{(3, 2, 1, 2, 4),\, (5,4,1,2,5)\}$, since we have
$b_3 = x t_4$ and $b_2 = y t_4$, as well as $b_5 = x t_5$ and $b_4 = y t_5$.
\end{example}

\goodbreak
The main property of the sets of polynomials associated to the neighbour pairs
is the following result shown in~\cite{KR3}, Prop.~4.1 and~\cite{KSL}, Prop.~2.5.

\begin{proposition}
The set $\ND_\OO \cup \AR_\OO$
contains (up to sign) the non-zero entries of all commutators
$\mathcal{A}_k\mathcal{A}_\ell -\mathcal{A}_\ell \mathcal{A}_k$ with
$1\le k<\ell \le n$. In particular, this set is a set of generators
of~$I(\BO)$. We call it the {\bf set of natural generators} of~$I(\BO)$.
\end{proposition}

Next, we recall some results of~\cite{KSL} about the
homogeneous components of the natural generators of~$I(\BO)$.

\begin{proposition}\label{prop-lin_cij}
Let $\OO=\{t_1,\dots,t_\mu\}$ be an order ideal in~$\mathbb{T}^n$ with
$\partial\OO= \{b_1,\dots, b_\nu\}$, and let $C=\{c_{ij} \mid i=1,\dots,\mu;\; j=1,\dots,\nu\}$.
\begin{enumerate}
\item[(a)] Let $\tau_{\alpha\beta}^{(k\ell)}$ be the $(\alpha,\beta)$-entry of the commutator
$\mathcal{A}_k \mathcal{A}_\ell - \mathcal{A}_\ell \mathcal{A}_k$.
Then $\tau_{\alpha\beta}^{(k\ell)}$ is of the form $\tau_{\alpha\beta}^{(k\ell)} =
\lambda_{\alpha\beta}^{(k\ell)} + q_{\alpha\beta}^{(k\ell)}$ with homogeneous polynomials
$\lambda_{\alpha\beta}^{(k\ell)} \in K[C]_1$ and $q_{\alpha\beta}^{(k\ell)} \in K[C]_2$.

\item[(b)] Let $i\in \{1,\dots,\mu\}$. For the natural generator $\ND(j, j')_i$,
where $b_j = x_\ell b_{j'}$ for some $j,j'\in\{1,\dots,\nu\}$ and $\ell\in\{1,\dots,n\}$, 
there are two possibilities:
\begin{itemize}
\item[(1)] If $t_i = x_\ell t_{i'}$ for some $i, i'\in\{1,\dots,\mu\}$ then the linear
part is (up to sign) given by $c_{ij}-c_{i'j'}$.

\item[(2)] If $t_i$ is not divisible by~$x_\ell$ then the linear part is (up to sign)
given by~$c_{ij}$.
\end{itemize}

\item[(c)] Let $m\in \{1,\dots,\mu\}$. For the natural generator $\AR(j,j')_m$ of~$I(\BO)$,
where $b_j = x_\ell t_{m'}$ and $b_{j'} = x_k t_{m'}$ for some $j,j'\in \{1,\dots,\nu\}$ and
$m'\in\{1,\dots,\mu\}$, there are several possibilities:
\begin{itemize}
\item[(1)] If $t_m = x_k t_i = x_\ell t_{i'}$ for some $i,i'\in \{1,\dots,\mu\}$ then
the linear part of $\AR(j,j')_m$ is (up to sign) given by $c_{ij}-c_{i'j'}$.

\item[(2)] If $t_m = x_k t_i$ and $t_m$ is not divisible by~$x_\ell$,
then the linear part of $\AR(j,j')_m$ is (up to sign) given by $c_{ij}$.

\item[(3)] If $t_m = x_\ell t_{i'}$ and $t_m$ is not divisible by~$x_k$
then the linear part of $\AR(j,j')_m$ is (up to sign) given by $c_{i' j'}$.

\item[(4)] If neither $x_k$ nor $x_\ell$ divide~$t_m$ then the linear part 
of $\AR(j,j')_m$ is zero.
\end{itemize}
\end{enumerate}
\end{proposition}

\begin{proof} Claim~(a) was shown in~\cite{KSL}, Prop.~2.7.c, and claims~(b), (c)
follow from~\cite{KSL}, Cor.~2.8.
\end{proof}

The following example illustrates this situation.

\begin{example}\label{ex-122}
In the polynomial ring $P=\mathbb{Q}[x,y]$, consider the order ideal 
${\OO = \{1,\, y,\, x,\, y^2,\, xy\}}$ and its border 
$\partial\OO = \{ x^2,\, y^3,\, xy^2,\, x^2y\}$.
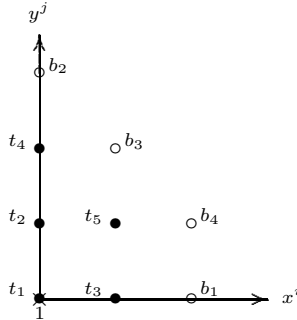
\begin{center}
	\makebox{\beginpicture
		\setcoordinatesystem units <1cm,1cm>
		\setplotarea x from 0 to 3, y from 0 to 3.5
		\axis left /
		\axis bottom /
		\arrow <2mm> [.2,.67] from  2.5 0  to 3 0
		\arrow <2mm> [.2,.67] from  0 3  to 0 3.5
		\put {$\scriptstyle x^i$} [lt] <0.5mm,0.8mm> at 3.1 0.1
		\put {$\scriptstyle y^j$} [rb] <1.7mm,0.7mm> at 0 3.6
		\put {$\bullet$} at 0 0
		\put {$\bullet$} at 1 0
		\put {$\bullet$} at 0 1
		\put {$\bullet$} at 1 1
		\put {$\bullet$} at 0 2
		\put {$\scriptstyle 1$} [lt] <-1mm,-1mm> at 0 0
		\put {$\scriptstyle t_1$} [rb] <-1.3mm,0.4mm> at 0 0
		\put {$\scriptstyle t_3$} [rb] <-1.3mm,0.4mm> at 1 0
		\put {$\scriptstyle b_1$} [lb] <0.8mm,0.4mm> at 2 0
		\put {$\scriptstyle t_2$} [rb] <-1.3mm,0mm> at 0 1
		\put {$\scriptstyle t_4$} [rb] <-1.3mm,0mm> at 0 2
		\put {$\scriptstyle t_5$} [rb] <-1.3mm,0mm> at 1 1
		\put {$\scriptstyle b_2$} [lb] <0.8mm,0mm> at 0 3
		\put {$\scriptstyle b_3$} [lb] <0.8mm,0mm> at 1 2
		\put {$\scriptstyle b_4$} [lb] <0.8mm,0mm> at 2 1
		\put {$\times$} at 0 0
		\put {$\circ$} at 0 3
		\put {$\circ$} at 2 0
		\put {$\circ$} at 2 1
		\put {$\circ$} at 1 2
		\endpicture}
\end{center}
Here the bullet points correspond to the terms in the order ideal~$\OO$
and the circles to its border $\partial\OO$.
The rim of~$\OO$ is given by $\nabla\OO = \{ x,\, y^2,\, xy\}$
and the interior of~$\OO$ is $\OO^\circ = \{ 1,\ y\}$.
Then the list of triples $(j,j',\ell)$ of next-door neighbours 
$b_j=x_\ell b_{j'}$ is $\{(4,1,2)\}$,
and the list of quintuples $(j,j',k,\ell,m')$ of across-the rim 
neighbours $x_k b_j = x_\ell b_{j'}$ is $\{(2,3,1,2,4),\, (3,4,1,2,5)\}$.

Now we want to look at the linear parts of the across-the-rim generators
for the quintuple $(2,3,1,2,4)$. In this case we have $(j,j')=(2,3)$, $(k,\ell)=(1,2)$,
$b_2=y^3$, $b_3=xy^2$, and $t_{m'}=t_4=y^2 \in \nabla\OO$.
The list $\AR(2,3)$ is given by 
\begin{align*}
\AR(2,3)_1 &=  c_{11}c_{32} +c_{13}c_{42} -c_{12}c_{43} 
               +c_{14}c_{52} -c_{13} c_{53},\\
\AR(2,3)_2 &=  c_{21}c_{32} +c_{23}c_{42} -c_{22}c_{43} 
               +c_{24}c_{52} -c_{23}c_{53} -c_{13},\\
\AR(2,3)_3 &= c_{31}c_{32} +c_{33}c_{42} -c_{32}c_{43} 
               +c_{34} c_{52} -c_{33} c_{53} +c_{12},\\
\AR(2,3)_4 &= c_{32}c_{41} +c_{44}c_{52} -c_{43}c_{53} -c_{23},\\
\AR(2,3)_5 &= c_{32}c_{51} -c_{43}c_{52} +c_{42}c_{53} -c_{53}^2  
              +c_{52} c_{54} +c_{22} -c_{33}
\end{align*}
In accordance with Proposition~\ref{prop-lin_cij}.c, we make the following observations. 
\begin{enumerate}
\item[(1)] $\AR(2,3)_5$: We have $m=5$ and $t_5=xy$. Then 
$t_5 = x_k\cdot t_2 = x_\ell\cdot t_3$, where 
$t_2=y$ and $t_3=x$. Letting $i=2$ and $i'=3$,
the linear part of $\AR(2,3)_5$ is given by
$$
\Lin(\AR(2,3)_5) = c_{22} -c_{33} = c_{i j} -c_{i' j'}
$$

\item[(2)] $\AR(2,3)_3$: We have $m=3$ and $t_3=x$. Then 
$t_3 = x_k\cdot t_1 = x\cdot t_1$ and $y=x_\ell \nmid t_3$. 
Hence we have $i=1$ and the linear part of $\AR(2,3)_3$ is given by
$$
\Lin(\AR(2,3)_3) = c_{12} = c_{i j}
$$

\item[(3)] $\AR(2,3)_4$: We have $m=4$ and $t_4=y^2$. Then we see that
$t_4 = x_\ell\cdot t_2 = y\cdot t_2$ and $x\nmid t_4$, 
and hence $i'=2$.
The linear part of $\AR(2,3)_4$ is given by
$$
\Lin(\AR(2,3)_4) = -c_{23} = -c_{i' j'}
$$

\item[(4)] $\AR(2,3)_1$: We have $m=1$ and $t_1=1$ is neither
divisible by $x_k=x$ nor by $x_\ell =y$.
Thus the linear part of $\AR(2,3)_1$ is zero.
\end{enumerate}
\end{example}

\bigskip\bigbreak
%
%

\section{The Arrow Degree}
\label{The Arrow Degree}

Let us continue to use the notation and definitions of the
preceding section. In particular, let $\OO = \{t_1,\dots,t_\mu\}$
be an order ideal in $\mathbb{T}^n$, let $\partial\OO =
\{b_1,\dots, b_\nu\}$, and let $C=\{ c_{ij} \mid i=1,\dots,\mu;
\; j=1,\dots,\nu\}$ be the set of coefficients of the generic 
$\OO$-border prebasis $G=\{g_1,\dots,g_\nu\}$, where
$g_j = b_j - \sum_{i=1}^\mu c_{ij}\, t_i$ for $j=1,\dots,\nu$.

For the very detailed analysis of the natural set of generators of~$I(\BO)$
performed in the subsequent parts of this paper, the following grading
will be very useful (see~\cite{Hai} and~\cite{Hui2}).

\begin{definition}\label{def-arrdeg}
Let $\OO \!=\! \{t_1,\dots,t_\mu\}$ be an order ideal in~$\mathbb{T}^n$, let
$\partial\OO\!=\! \{b_1,\dots, b_\nu\}$, and let
$K[C] = K[c_{ij} \mid i=1,\dots,\mu;\; j=1,\dots,\nu]$, where~$K$ is a field.
\begin{enumerate}
\item[(a)] The $\mathbb{Z}^n$-grading on~$K[C]$ defined by
$W\in\Mat_{n,\mu\nu}(\mathbb{Z})$, where we let
$\deg_W(c_{ij}) = \log(b_j) - \log(t_i)$ for $i=1,\dots,\mu$
and for $j=1,\dots,\nu$, is called the {\bf arrow grading}.

\item[(b)] The $\mathbb{Z}$-grading on~$K[C]$ defined by
letting $\deg(c_{ij}) = \deg(b_j) - \deg(t_i)$ for $i=1,\dots,\mu$
and $j=1,\dots,\nu$ is called the {\bf total arrow grading}.
\end{enumerate}
\end{definition}

The following result is shown in~\cite{KSL}, Lemma~3.4.

\begin{proposition}\label{prop-arrowHomog}
The ideal~$I(\BO)$ is homogeneous with respect to the arrow grading
on~$K[C]$.
In particular, since the total arrow grading is coarser than the arrow grading,
the ideal $I(\BO)$ is also homogeneous with respect to the total arrow grading.
\end{proposition}

As for the arrow degrees of the indeterminates~$c_{ij}$, we have the
following easy observation (see also~\cite{Hui2}, Lemma 4.1.1).

\begin{proposition}\label{prop-poscomp}
For all $i\in\{1,\dots,\mu\}$ and $j\in\{1,\dots,\nu\}$, the arrow degree
$\deg_W(c_{ij})$ has at least one positive component.
\end{proposition}

\begin{proof}
Suppose that $i\in\{1,\dots,\mu\}$ and $j\in\{1,\dots,\nu\}$
such that all components of $\log(b_j) -\log(t_i)$ are negative.
This means that $b_j \mid t_i$, and therefore the fact that~$\OO$
is an order ideal implies $b_j\in\OO$, a contradiction.
\end{proof}

This proposition inspires the following definition which is in accordance
with~\cite{Hui2}, Section~4.1.

\begin{definition}
For $i\in \{1,\dots,\mu\}$ and $j\in\{1,\dots,\nu\}$, the
indeterminate $c_{ij}$ is said to be a {\bf standard indeterminate},
and $\deg_W(c_{ij})$ is called a {\bf standard arrow}, if
$\deg_W(c_{ij})$ contains exactly one positive component.

If $\deg_W(c_{ij})$ contains at least two positive entries,
we say that $c_{ij}$ is a {\bf positive indeterminate} in~$C$
and that $\deg_W(c_{ij})$ is a {\bf positive arrow.}
\end{definition}

Let us examine these notions in a concrete case.

\begin{example}\label{ex-StdPosIndets}
In $P=\mathbb{Q}[x,y]$, look at the order ideal 
$\OO = \{t_1,\dots,t_5\}$, where $t_1=1$, $t_2=y$, $t_3=x$, $t_4=y^2$, 
and $t_5=xy$. Then the border of~$\OO$ is given by $\partial\OO = \{b_1,\dots,b_4\}$,
where $b_1=x^2$, $b_2=y^3$, $b_3=xy^2$, and $b_4=x^2y$ (see Example~\ref{ex-122}).
We have $\deg_W(c_{11})=(2,0)$, $\deg_W(c_{12})=(0,3)$,
$\deg_W(c_{13})=(1,2)$, $\deg_W(2,1)=(2,-1)$, etc.
The set of standard indeterminates is 
$$
\{c_{11}, c_{12}, c_{22}, c_{24}, c_{31}, c_{33}, c_{42}, c_{43},
c_{53}, c_{54} \}
$$ 
and the set of positive indeterminates is 
$\{ c_{13}, c_{14}, c_{23}, c_{34} \}$.
\end{example}

Of course, knowing that the ideal $I(\BO)$ is homogeneous with respect to the
arrow grading is useful. It is even better to know the
arrow degrees of the various generators. This information is provided
by the next proposition.

\goodbreak
\begin{proposition}
Let $\OO=\{t_1,\dots,t_\mu\}$ be an order ideal in~$\mathbb{T}^n$
with border $\partial\OO = \{b_1,\dots, b_\nu\}$.
\begin{enumerate}
\item[(a)] Let $j, j'\in\{1,\dots,\nu\}$ be such that $b_j,b_{j'}$
are next-door neighbours with $b_j = x_\ell b_{j'}$, and let $i\in\{1,\dots,\mu\}$.
Then the arrow degree of $\ND(j, j')_i$ is $\deg_W(c_{ij})$.

\item[(b)] Let $j, j'\in \{1,\dots,\nu\}$ be such that $b_j, b_{j'}$
are across-the-rim neighbours with $b_j = x_\ell t_{m'}$ and $b_{j'} = x_k t_{m'}$
for some $m'\in\{1,\dots,\mu\}$.
Then the arrow degree of $\AR(j, j')_m$ is given by $\deg_W(c_{mj}) + e_k =
\deg_W(c_{mj'})+e_\ell$.
\end{enumerate}
\end{proposition}

\begin{proof}
This follows from Proposition~\ref{prop-lin_cij}.
\end{proof}

Using the arrow degrees, we can visualize parts~(b) and~(c)
of Proposition~\ref{prop-lin_cij} as follows.

\begin{remark}\label{rem-LinParts}
Let $\OO=\{t_1,\dots,t_\mu\}$ be an order ideal in~$\mathbb{T}^n$
with border $\partial\OO = \{b_1,\dots, b_\nu\}$.
\begin{enumerate}
\item[(a)] Let $j, j'\in\{1,\dots,\nu\}$ be such that $b_j,b_{j'}$
are next-door neighbours with $b_j = x_\ell b_{j'}$, and let $i\in\{1,\dots,\mu\}$.
Then the polynomial $\ND(j,j')_i$ has (up to sign) the linear part $c_{ij} - c_{i'j'}$
if $x_\ell \mid t_i$. This can be visualized as the equivalence
of the dashed arrows in the picture below.

\item[(b)] Let $j, j'\in \{1,\dots,\nu\}$ be such that $b_j, b_{j'}$
are across-the-rim neighbours with $b_j = x_\ell t_{m'}$ and $b_{j'} = x_k t_{m'}$
for some $m'\in\{1,\dots,\mu\}$, and let $m\in \{1,\dots,\mu\}$.
If $t_m = x_k t_i  = x_\ell t_{i'}$ for some $i,i'\in \{1,\dots,\mu\}$
then the polynomial $\AR(j,j')_m$ has (up to sign) the linear part $c_{ij} - c_{i' j'}$.
This can be visualized as the equivalence of the dashed arrows in the picture below.
\end{enumerate}
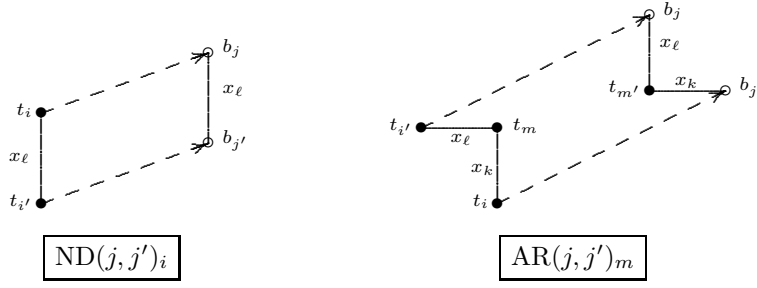
\begin{center}
	\makebox{\beginpicture
		\setcoordinatesystem units <1cm,1cm>
		\put {\fbox{$\mathrm{ND}(j,j')_i$}} [lt] <-1.7mm,0.8mm> at 0.2 -0.5
		\put {$\scriptstyle t_{i'}$} [rt] <-1.7mm,0.8mm> at 0.1 0.05
		\put {$\scriptstyle t_{i}$} [rt] <-1.7mm,0.8mm> at 0.1 1.25
		\put {$\bullet$} at 0 0
		\put {$\bullet$} at 0 1.2
		\put {$\scriptstyle x_\ell$} [rt] <-1.7mm,0.8mm> at 0.06 0.6
		\put {\fbox{$\mathrm{AR}(j,j')_m$}} [lt] <-1.7mm,0.8mm> at 6.2 -0.5
		\put {$\scriptstyle t_i$} [rt] <-1.7mm,0.8mm> at 6.1 0.05
		\put {$\scriptstyle t_{i'}$} [rt] <-1.7mm,0.8mm> at 5.1 1.05
		\put {$\scriptstyle t_{m}$} [lt] <0.7mm,0.8mm> at 6.1 1.05
		\put {$\scriptstyle x_{\ell}$} [rt] <0.5mm,0.8mm> at 5.65 0.85
		\put {$\scriptstyle x_{k}$} [rt] <0.5mm,0.8mm> at 5.95 0.45
		\put {$\bullet$} at 6 0
		\put {$\bullet$} at 5 1		
		\put {$\bullet$} at 6 1	
		\put {$\scriptstyle b_{j'}$} [lt] <0.5mm,0.8mm> at 2.3 0.9
		\put {$\scriptstyle b_{j}$} [lt] <0.5mm,0.8mm> at 2.3 2.1
		\put {$\circ$} at 2.2 2
		\put {$\circ$} at 2.2 0.8
		\put {$\scriptstyle x_\ell$} [lt] <0.5mm,0.8mm> at 2.3 1.5
		\put {$\scriptstyle b_{j'}$} [lt] <0.5mm,0.8mm> at 9.1 1.6
		\put {$\scriptstyle b_{j}$} [lt] <0.5mm,0.8mm> at 8.1 2.6
		\put {$\scriptstyle t_{m'}$} [rt] <-1.7mm,0.8mm> at 8.1 1.55
		\put {$\scriptstyle x_{\ell}$} [lt] <0.5mm,0.8mm> at 8.05 2.1
		\put {$\scriptstyle x_{k}$} [lt] <0.5mm,0.8mm> at 8.25 1.6
		\put {$\circ$} at 8 2.5
		\put {$\circ$} at 9 1.5
		\put {$\bullet$} at 8 1.5	
		\setlinear
		\plot 0 0  0 1.2 /
		\plot 2.2 2  2.2 0.8 /
		\plot 6 0  6 1 /
		\plot 5 1  6 1 /
		\plot 8 2.5  8 1.5 /
		\plot 9 1.5  8 1.5 /
		\arrow <2mm> [.2,.67] from 2 0.727 to 2.2 0.8
		\arrow <2mm> [.2,.67] from 2 1.927 to 2.2 2
		\arrow <2mm> [.2,.67] from 8.9 1.45 to 9 1.5
		\arrow <2mm> [.2,.67] from 7.9 2.45 to 8 2.5
		\setdashes 
		\plot 0 0  2.2 0.8 /
		\plot 0 1.2  2.2 2 /
		\plot 6 0  9 1.5 /
		\plot 5 1  8 2.5 /
		\dimen0=1cm \dimen0=1cm
		\endpicture}
\end{center}
For the cases that $x_\ell \nmid t_i$ in~(a) or $x_\ell \nmid t_m$ in~(b)
or $x_k \nmid t_m$ in~(b), we leave the analogous visualizations to the
interested reader.
\end{remark}

Our next result yields important restrictions for the indeterminates
appearing in the quadratic parts of the natural generators of~$I(\BO)$.
It uses the following notion (see also~\cite{Hui1}, Section 4.1).

\begin{definition}\label{def-exposedterm}
Let $j\in\{1,\dots,\nu\}$, and let $\ell\in\{1,\dots,n\}$.
Then the border term~$b_j$ is called {\bf $x_\ell$-exposed}
if it is of the form $b_j = x_\ell \, t_i$ with
$i \in \{1,\dots,\mu\}$. In this case we also say that~$t_i$
{\bf $x_\ell$-exposes} the border term~$b_j$.
\end{definition}

Now we are ready to formulate the proposition.

\begin{proposition}\label{prop-QuadrParts}
Let $\OO=\{t_1,\dots,t_\mu\}$ be an order ideal in~$\mathbb{T}^n$
with border $\partial\OO = \{b_1,\dots, b_\nu\}$.
\begin{enumerate}
\item[(a)] Let $j, j'\in\{1,\dots,\nu\}$ be such that $b_j,b_{j'}$
are next-door neighbours with $b_j = x_\ell b_{j'}$, and let $i\in \{1,\dots,\mu\}$.
Assume that $b_{\lambda_1}, \dots, b_{\lambda_s} \in\partial\OO$ are the
$x_\ell$-exposed border terms, and let $b_{\lambda_p} = x_\ell \, t_{\rho_p}$
for $p=1,\dots,s$. Then the quadratic terms in the support of~$\ND(j,j')_i$
are the products $c_{i \lambda_p}\, c_{\rho_p j'}$ with $p=1,\dots,s$.

\item[(b)] Let $j,j'\in \{1,\dots,\nu\}$ be such that $b_j, b_{j'}$
are across-the-rim neighbours with $b_j = x_\ell t_{m'}$ and $b_{j'} = x_k t_{m'}$
for some $m'\in\{1,\dots,\mu\}$. For every $m\in\{1,\dots,\mu\}$,
the quadratic terms in the support of $\AR(j,j')_m$ are of the following types.
\begin{itemize}
\item[(1)] Let $b_{\kappa_1},\dots, b_{\kappa_s} \in \partial\OO$ be the
$x_k$-exposed border terms, and let us write $b_{\kappa_p} = x_\ell\, t_{\rho_p}$
for $p=1,\dots,s$. 
Then the terms  $c_{m \kappa_p} \, c_{\rho_p j}$ may appear in the support of $\AR(j,j')_m$.

\item[(2)] Let $b_{\lambda_1}, \dots,b_{\lambda_u} \in \partial\OO$ be the
$x_\ell$-exposed border terms, and let us write $b_{\lambda_q} = x_\ell t_{\sigma_q}$
for $q=1,\dots,u$. Then the terms  $c_{m \lambda_q}\, c_{\sigma_q j'}$
may appear in the support of $\AR(j,j')_m$.

\end{itemize}
\end{enumerate}
\end{proposition}

\begin{proof}
Let $G=\{g_1,\dots,g_\nu\}$ be the generic $\OO$-border prebasis.

First we prove~(a). The polynomials in $\ND(j,j')$ are the coefficients of the
terms $t_1,\dots,t_\mu$ in the reduction $x_\ell g_{j'} - g_j$, viewed as a polynomial in
$K[C][x_1,\dots,x_n]$. The leading terms $x_\ell b_{j'}$ and $b_j$ cancel by
definition. Hence we only have to reduce the $x_\ell$-exposed border terms
in $x_\ell g_{j'}$. The coefficient of~$b_{\lambda_p}$ in $x_\ell g_{j'}$
is $c_{\rho_p j'}$. The coefficient of~$t_i$ in~$g_{\lambda_p}$ is $c_{i \lambda_p}$.
Therefore the coefficient of~$t_i$ in $c_{\rho_p j'}\, g_{\lambda_p}$
is $c_{\rho_p j'} c_{i \lambda_p}$ and it will be a part of the coefficient
of~$t_i$ in the final result of the reduction, i.e., in $\ND(j,j')_i$.

Next we prove~(b). Again the entries of $\AR(j,j')$ are obtained
as the coefficients of $t_1,\dots,t_\mu$ in the reduction of $x_k g_j - x_\ell g_{j'}$,
and again the highest terms $x_k b_j$ and $x_\ell b_{j'}$ cancel.
To reduce $x_k g_j$, we have to reduce the $x_k$-exposed border terms
$b_{\kappa_1},\dots, b_{\kappa_s}$. The coefficient of $b_{\kappa_p}$ in $x_k g_j$ is
$c_{\rho_p j}$. The polynomial $\AR(j,j')_m$ has arrow degree $\deg_W(c_{ij})+e_k$.
It is therefore the coefficient of~$t_m$ in the reduction of $x_k g_j - x_\ell
g_{j'}$. The coefficient of $t_m$ in $g_{\kappa_p}$ is~$c_{m \kappa_p}$.
Hence the coefficient of~$t_m$ in $c_{\rho_p j}\, g_{\kappa_p}$ is
$c_{\rho_p j} c_{m \kappa_p}$, and it may appear in $\AR(j,j')_m$.
These quadratic terms are the ones listed in case~(1).

The analysis of the terms in the reduction of $x_\ell g_{j'}$
is completely analogous and leads to the quadratic terms in case~(2).
\end{proof}

To illustrate this proposition, we apply it to the above example.

\begin{example}\label{ex-QuadraticTerms}
Assume the setting of Example~\ref{ex-122}. Let $P=\mathbb{Q}[x,y]$, and 
let $\OO = \{t_1,\dots,t_5\}$ be the order ideal given by $t_1=1$, $t_2=y$, 
$t_3=x$, $t_4=y^2$, and $t_5=xy$. Then the $x$-exposed border terms are 
$b_1=x^2$, $b_3=xy^2$, and $b_4=x^2y$, while the $y$-exposed border terms are 
$b_2=y^3$ and $b_3=xy^2$. For the next-door neighbour relation $b_4 = yb_1$, 
the list $\ND(4,1)$ is given by
\begin{align*}
\ND(4,1)_1 &= -c_{12} c_{41} -c_{13} c_{51} +c_{14}\\ 
\ND(4,1)_2 &= -c_{22} c_{41} -c_{23} c_{51} -c_{11} +c_{24}\\
\ND(4,1)_3 &=-c_{32} c_{41} -c_{33} c_{51} +c_{34}\\
\ND(4,1)_4 &=-c_{41} c_{42} -c_{43} c_{51} -c_{21} +c_{44}\\
\ND(4,1)_5 &=-c_{41} c_{52} -c_{51} c_{53} -c_{31} +c_{54}
\end{align*}
By Proposition~\ref{prop-QuadrParts}.a, the quadratic terms $c_{1\lambda_p}c_{\varrho_p1}$ 
in the support of $\ND(4,1)_1$ satisfy $\lambda_p\in\{2,3\}$ and $\varrho_p \in\{4,5\}$.
They are given by $c_{12} c_{41}$ and $c_{13} c_{51}$. The arrow degree 
of $\ND(4,1)_1$ is given by $\deg_W(\ND(4,1)_1) = \deg_W(c_{14}) =(2,1)$.

Next we consider the across-the rim neighbour relation $xb_2 = yb_3$. 
Using $m=2$ in Example~\ref{ex-122}, we have 
$$
\AR(2,3)_2 =  c_{21}c_{32} +c_{23}c_{42} -c_{22}c_{43} 
+c_{24}c_{52} -c_{23}c_{53} -c_{13}
$$
According to Proposition~\ref{prop-QuadrParts}.b, the terms 
$c_{2\kappa_p}c_{\varrho_p2}$ of $\AR(2,3)_2$ 
with $x$-exposed border indices $\kappa_p\in\{1,3,4\}$ are 
$c_{21}c_{32}$, $c_{23}c_{42}$ and $c_{24}c_{52}$,
while the quadratic terms $c_{2\lambda_q}c_{\sigma_q3}$ 
of $\AR(2,3)_2$ with $y$-exposed border indices $\lambda_q\in\{2,3\}$ 
are $c_{22}c_{43}$ and $c_{23}c_{53}$.
In this case we have $\deg_W(\AR(2,3)_2) = \deg_W(c_{13}) = (1,2)$.
\end{example}

For later applications, it is useful to reformulate the results of this
proposition using the following terminology.

\begin{definition}\label{def-interiorrim}
Let $\OO = \{t_1,\dots,t_\mu\}$ be an order ideal in~$\mathbb{T}^n$ with border  
$\partial\OO \!=\! \{b_1,\dots,b_\nu\}$, and let
$C = \{ c_{ij} \mid i=1,\dots,\mu;\; j=1,\dots,\nu\}$. 
\begin{enumerate}
\item[(a)] The set~$C^\circ$ of all indeterminates $c_{ij}$ such that
$t_i \in \OO^\circ$ is called the set of {\bf interior indeterminates}.

\item[(b)] The set $C^{\nabla}$ of all indeterminates $c_{ij}$ such that
$t_i \in \nabla \OO$ is called the set of {\bf rim indeterminates}.
\end{enumerate}
\end{definition}

Once again, let us check the meaning of this terminology using
Example~\ref{ex-122}.

\begin{example}\label{ex-InteriorRimIndets}
In the setting of Example~\ref{ex-122}, we have the equalities
$\OO^\circ = \{ t_1,\ t_2\}$ and $\nabla\OO = \{ t_3,\, t_4,\, t_5\}$.
Consequently, we deduce that the set of interior indeterminates is 
$C^\circ =\{c_{ij} \mid i=1,2; \, j=1,2,3,4\}$
and  the set of rim indeterminates is 
$C^{\nabla} =\{c_{ij} \mid i=3,4,5;\, j=1,2,3,4\}$.
\end{example}

Using this definition, we can combine the statements of
Propositions~\ref{prop-lin_cij} and~\ref{prop-QuadrParts} as follows.

\begin{corollary}\label{cor-interpretQuadrPart}
Let $\OO=\{t_1,\dots,t_\mu\}$ be an order ideal in~$\mathbb{T}^n$
with border $\partial\OO = \{b_1,\dots, b_\nu\}$.
\begin{enumerate}
\item[(a)] Let $j, j'\in\{1,\dots,\nu\}$ be such that $b_j,b_{j'}$
are next-door neighbours with $b_j = x_\ell b_{j'}$, and let $i\in \{1,\dots,\mu\}$.
Then the linear part of $\ND(j,j')_i$ is (up to sign) given by
$$
\begin{cases}
c_{ij}-c_{i' j'}  &  \text{ if $x_\ell$ divides $t_i$,} \\
c_{ij} & \text{ otherwise,}
\end{cases}
$$
and all terms in the quadratic part are of the form $c_{i\lambda} c_{\rho j'}$
with $\lambda\in \{1,\dots,\nu\}$ and a rim indeterminate $c_{\rho j'}$.

\item[(b)]  Let $j,j'\in \{1,\dots,\nu\}$ be such that $b_j, b_{j'}$
are across-the-rim neighbours with $b_j = x_\ell t_{m'}$ and $b_{j'} = x_k t_{m'}$
for some $m'\in\{1,\dots,\mu\}$. For $m\in\{1,\dots,\mu\}$, the linear part of
$\AR(j,j')_m$ is (up to sign) given by
$$
\begin{cases} c_{ij}-c_{i'j'} & \text{ if $x_k$ and $x_\ell$ divide $t_m$,}\\
c_{ij} & \text{ if $x_k$ divides $t_m$, but $x_\ell$ doesn't,}\\
c_{i'j'} & \text{ if $x_\ell$ divides $t_m$, but $x_k$ doesn't,}
\end{cases}
$$ 
and the terms in the quadratic part
are of the form $c_{m\lambda} c_{\rho j}$ or $c_{m\lambda} c_{\rho j'}$
with $\lambda \in \{1,\dots,\nu\}$ and rim indeterminates $c_{\rho j}$ and 
$c_{\rho j'}$, respectively.
\end{enumerate}
\end{corollary}

\begin{proof}
This follows immediately from Proposition~\ref{prop-lin_cij} 
and Proposition~\ref{prop-QuadrParts}.
\end{proof}

\bigskip\bigbreak
%
%

\section{Cotangent Equivalence Classes and Optimal Embeddings}
\label{Cotangent Equivalence Classes and Optimal Embeddings}

As before, we let $P=K[x_1,\dots,x_n]$ be a polynomial ring over a field~$K$,
let $\OO = \{ t_1,\dots,t_\mu \}$ be an order ideal in~$\mathbb{T}^n$,
let $\partial\OO = \{b_1,\dots,b_\nu \}$ be its border, and let
$C = \{ c_{ij} \mid i=1,\dots,\mu;\; j=1,\dots,\nu \}$ be the
set of coefficients of the generic $\OO$-border prebasis.
Recall that the affine coordinate ring of~$\BO$ is
$B_\OO = K[C] / I(\BO)$, and let $\m = \langle c_{ij} +I(\BO) \mid
i=1,\dots,\mu;\; j=1,\dots,\nu \rangle$ be the maximal ideal of~$B_\OO$
generated by the residue classes of the indeterminates~$c_{ij}$.

In Section~\ref{Separating Re-embeddings via Groebner Fans of Linear Ideals}
we discussed how to use $\GFan(\langle\Lin_\M(I(\BO))\rangle)$ to 
detect $Z$-sep\-arating, possibly optimal,  re-embeddings of~$I(\BO)$.
In the framework of border basis schemes further tools are available.
Example~\ref{ex-Section5} at the end of this section illustrates 
their applicability.

\begin{definition}\label{def-cotequiv}
For every indeterminate $c_{ij}\in C$, let $\bar{c}_{ij}$ denote
its residue class in the cotangent space $\m / \m^2$ 
of $B_\OO$ at the origin.
\begin{enumerate}
\item[(a)] The relation~$\sim$ on~$C$ defined
by $c_{ij} \sim c_{i'j'} \Leftrightarrow \bar{c}_{ij} = \bar{c}_{i'j'}$
is an equivalence relation called {\bf cotangent equivalence}.

\item[(b)] An indeterminate $c_{ij}\in C$ is called {\bf trivial}
if $\bar{c}_{ij}=0$. The trivial indeterminates form the {\bf trivial
cotangent equivalence class} in~$C$.

\item[(c)] A non-trivial indeterminate $c_{ij} \in C$ is called {\bf basic} if
its cotangent equivalence class consists only of~$c_{ij}$. In this case, 
the cotangent equivalence class $\{c_{ij}\}$ is also called {\bf basic}.

\item[(d)] A non-trivial indeterminate $c_{ij} \in C$ is called {\bf proper}
if its cotangent equivalence class contains at least two elements. In this case,
the cotangent equivalence class of~$c_{ij}$ is also called {\bf proper}.

\end{enumerate}
\end{definition}

The meaning of these notions will become clear in the next theorems.
First we need a lemma which provides further information about the above definition.

\begin{lemma}\label{lem-equivclasses}
Let us assume to be in the above setting.
\begin{enumerate}
\item[(a)] The union~$U$ of the  supports of the elements in  
a minimal set of generators of the ideal $\langle \Lin_\M(I(\BO))\rangle$ 
does not depend on the choice of a minimal set of generators. 

\item[(b)] The set of basic indeterminates is $C\setminus U$.

\item[(c)] The union of the sets of trivial and proper indeterminates is~$U$.

\item[(d)] The set of basic indeterminates is contained in any set of indeterminates
whose residue classes modulo $I(\BO)$ generate the $K$-algebra $B_\OO$.
\end{enumerate}
\end{lemma}

\begin{proof}
To prove (a) we note that any minimal set of generators of
$\langle \Lin_\M(I(\BO))\rangle$ is also a minimal set of generators 
of the $K$-vector space $\Lin_\M(I(\BO))$.
Let~$A$ and~$B$ be two such sets. Since every linear form~$\ell$ in~$A$ 
is a linear combination of linear forms in~$B$, each indeterminate in $\Supp(\ell)$ 
is in the support of some linear form in~$B$. By interchanging the roles of~$A$ and~$B$, 
the conclusion follows.

Since $\sim$ is an equivalence relation on~$C$, to prove claims~(b) and~(c) it suffices 
to show that a basic indeterminate is not in~$U$, while trivial and proper indeterminates
are in~$U$.
For a contradiction, assume that $c_{ij} \in U$ is a basic indeterminate.
From $c_{ij}\in U$ and the fact that $c_{ij}$ is the only element in its equivalence class, 
we deduce that there is a polynomial in~$I(\BO)$ of the form $c_{ij} +q$ 
with $q \in \M^2$. Hence we get $\bar{c}_{ij} = 0$, a contradiction.
Next let $c_{ij}$ be trivial.  Then there is a polynomial in $I(\BO)$ of the form $c_{ij} +q$
with $q \in \M^2$, and hence we get $c_{ij} \in U$. Next, let $c_{ij}$ be proper. 
Then there exist another indeterminate $c_{i'j'}$ and a polynomial in $I(\BO)$ of the form 
$c_{ij} -c_{i'j'}+q$ with $q \in \M^2$. Hence we get $c_{ij} \in U$.

Finally, we prove~(d). From~(b) we know that the basic indeterminates do not belong to the 
support of any generator of $\Lin_\M(I(\BO))$. Hence the claim follows from 
Proposition~\ref{prop-cohSepPoly}.c.
\end{proof}

\begin{definition}\label{def-leadsetofPr}
Let $E= \{c_{i_1 j_1}, \dots, c_{i_p j_p} \}$ be a proper equivalence class in~$C$,
and let $\sigma$ be a term ordering on ${\mathbb T}^{\mu\nu}$
with $c_{i_1 j_1} >_\sigma \cdots  >_\sigma c_{i_p j_p} $. Then the set
$E\setminus\{c_{i_p j_p}\}$ is called the {\bf $\sigma$-leading set} of~$E$ and
denoted by~$E^\sigma$.
\end{definition}

\begin{notation}
Let $\sigma$ be a term ordering on~${\mathbb T}^{\mu\nu}$.
In accordance with the notation introduced in Proposition~\ref{prop-cohSepPoly},
the unique minimal set of terms generating the ideal 
$\LT_\sigma( \langle\Lin_ \M(I(\BO))\rangle )$ will be denoted by $S_\sigma$.
\end{notation}

In the following theorem we give an explicit representation of~$S_\sigma$
for every term ordering~$\sigma$ on $\mathbb{T}^{\mu\nu}$
and show the importance of $\sigma$-leading sets.
This implies  a description of $\LTGFan(\langle \Lin_\M(I(\BO))\rangle)$
which has several advantages when compared with the general 
description given in Theorem~\ref{thm-bijection}. 
For instance, it does not have to deal with huge matrices.

\begin{theorem}\label{thm-shapeofSsigma}
Let $E_0$ be the trivial equivalence class and $E_1, \dots, E_q$ 
the proper equivalence classes in~$C$. Let $I(\BO)_\lin = \langle \Lin_\M(I(\BO))\rangle$
be the ideal generated by the linear parts of the polynomials in~$I(\BO)$.
\begin{enumerate}
\item [(a)] Let~$\sigma$ be a term ordering on $\mathbb{T}^{\mu\nu}$.
Then we have $S_\sigma = E_0 \cup E_1^\sigma \cup \cdots \cup E_q^\sigma$,
and hence $\# S_\sigma  = \# E_0 + \sum_{i=1}^q \# E_i - q$.

\item[(b)] For $i=1,\dots,q$, let $E_i^\ast$ be a set obtained from~$E_i$ 
by deleting one of its elements. Then there exists a term ordering~$\sigma$
on $\mathbb{T}^{\mu\nu}$ such that we have
$S_\sigma = E_0 \cup E_1^\ast \cup \cdots \cup E_q^\ast$.

\item[(c)]  Let $\Sigma$ be the set of all sets of the form 
$E_0 \cup E_1^\ast \cup \cdots \cup E_q^\ast$, where~$E_i^\ast$ is obtained
from the set~$E_i$ by deleting one of its elements.
Then the map $\phi: \LTGFan(I(\BO)_\lin )  \longrightarrow \Sigma$ given by 
$\phi(S_\sigma) = E_0 \cup E_1^\sigma \cup \cdots \cup E_q^\sigma$ 
is well-defined and bijective.

\item[(d)]  We have
$\#\LTGFan(I(\BO)_\lin) = \#\GFan(I(\BO)_\lin)= \prod\limits_{i=1}^q \# E_i$.
\end{enumerate}
\end{theorem}

\begin{proof} 
To prove claim~(a) we observe that the inclusion $E_0 \subseteq S_\sigma$ follows from 
$E_0 \subseteq \Lin_\M(I(\BO))$, and that  the inclusion 
$S_\sigma \subseteq E_0 \cup E_1\cup \cdots \cup E_q$
follows from Lemma~\ref{lem-equivclasses}.c.

For $k\in \{1,\dots,q\}$, we write the proper equivalence class 
$E_k= \{c_{i_1j_1}, \dots, c_{i_p j_p} \}$ 
such that $c_{i_1 j_1} >_\sigma \cdots  >_\sigma c_{i_p j_p}$. 
Using the definition of a proper equivalence class, it follows that 
$c_{i_\ell j_\ell} - c_{i_p j_p} \in \Lin_\M(I(\BO))$, so that
$\LT_\sigma(c_{i_\ell j_\ell} - c_{i_p j_p}) = c_{i_\ell j_\ell}$ 
for every $\ell \in \{1,\dots,p-1\}$.
Hence we have proved 
$$
E_0 \cup E_1^\sigma \cup \cdots \cup E_q^\sigma \;\subseteq\;  
S_\sigma \;\subseteq\;   E_0 \cup E_1\cup \cdots \cup E_q
$$ 

Therefore the $\sigma$-smallest 
element in each proper equivalence class does not belong to~$S_\sigma$, 
which finishes the proof of~(a).
 
Claim~(b) follows from~(a) if we show that there exists a term ordering~$\sigma$ 
such that $E_i^\sigma = E_i^\ast$ for $i=1,\dots,q$.
By definition, the sets~$E_i$ are pairwise disjoint. Consequently, a term ordering 
which solves the problem can be chosen as a block term ordering, and hence it 
suffices to consider the case $q=1$. 
So, let $E_1 = \{c_{i_1j_1}, \dots, c_{i_p j_p} \}$  and w.l.o.g. assume 
that $E_1^\ast = \{c_{i_2 j_2}, \dots, c_{i_p j_p} \}$. As we observed before, we have 
$c_{i_\ell j_\ell} - c_{i_1 j_1} \in \Lin_\M(I(\BO))$ for $\ell=2,\dots,p$.
To finish the proof it suffices to take a term ordering~$\sigma$ such that 
$c_{i_\ell j_\ell} >_\sigma c_{i_1 j_1}$ for $\ell=2,\dots,p$.
  
Claim~(c) follows immediately from~(b).
To prove claim~(d) we note that the first equality is obvious.  
To show $\# \LTGFan(I(\BO)_\lin ) =  \prod_{i=1}^q \# E_i$, it suffices to 
deduce from~(b) that the number of the leading term ideals of $I(\BO)_\lin$ 
equals the number of $q$-tuples of indeterminates $c_{ij}$, exactly one chosen in
each proper equivalence class.
\end{proof}

Now we are ready to classify the indeterminates in~$C$ which can be used for a 
$Z$-separating re-embedding of~$I(\BO)$ as follows.

\begin{theorem}\label{thm-classifycij}
Let~$Z$ be a tuple of indeterminates from~$C$ such that there exists a
$Z$-separating re-embedding of~$I(\BO)$ and let $Y = C \setminus Z$.
\begin{enumerate}
\item[(a)] The basic indeterminates of~$C$ are contained in~$Y$.

\item[(b)] Each proper equivalence class in~$C$ contains at least one element of~$Y$.

\item[(c)] All basic indeterminates of~$C$ are rim indeterminates.

\item[(d)] Each proper cotangent equivalence class in~$C$ contains a rim indeterminate.

\item[(e)] If $\# Z = \dim _K(\Lin_\M(I(\BO))$, then the $Z$-separating 
re-em\-bed\-ding of~$I(\BO)$ is optimal, the trivial indeterminates of~$C$ are contained 
in~$Z$, and each proper equivalence class in~$C$ contains exactly one element of~$Y$.

\end{enumerate}
\end{theorem}

\begin{proof}
To prove~(a), let~$\sigma$ be a $Z$-separating term ordering for $I(\BO)$,
let~$S_\sigma$ be the minimal set of terms generating 
$\LT_\sigma(\langle \Lin_{\M}(I(\BO)\rangle )$,
and let~$U$ be the union of the supports of the elements in a minimal set of generators 
of $\langle \Lin_{\M}(I(\BO)\rangle$. Note that Proposition~\ref{prop-cohSepPoly}.c 
implies $Y \supseteq C \setminus S_\sigma$. From $S_\sigma \subseteq U$ we deduce the inclusion 
$Y \supseteq C \setminus U$, and thus the claim follows from Lemma~\ref{lem-equivclasses}.b.
 
Claim~(b) follows from the mentioned inclusion $Y \supseteq C \setminus S_\sigma$ and 
Theorem~\ref{thm-shapeofSsigma}.a.

For the proof of~(c), we assume that $c_{ij} \in C$ is a basic indeterminate and that
$t_i\in \OO^\circ$. If there exists an indeterminate $x_\ell$ such that
$x_\ell b_j \in \partial\OO$ then we let $j'\in\{1,\dots,\nu\}$ such that 
$b_{j'} = x_\ell b_j$. As we have $t_i\in\OO^\circ$, there exists an index
$i'\in\{1,\dots,\mu\}$ such that $t_{i'} = x_\ell t_i \in \OO$.
Hence we have a next-door neighbour pair $\ND(j,j')$ and 
Proposition~\ref{prop-lin_cij} implies $c_{ij} \sim c_{i' j'}$. This contradicts
the hypothesis that $c_{ij}$ is a basic indeterminate. Thus no next-door
neighbour pair $(b_j,b_{j'})$ exists.
Since the border~$\partial\OO$ is connected with respect to the neighbour
relations, it follows that there exists an across-the-rim neighbour pair $(b_j,b_{j'})$.
Thus we may assume that there are $m'\in\{1,\dots,\mu\}$ and $k,\ell \in\{1,\dots,n\}$
such that $b_j = x_\ell t_{m'}$ and $b_{j'} = x_k t_{m'} \in\partial\OO$
for some $j'\in \{1,\dots,\nu\}$. Now, if $t_i$ is not divisible by~$x_\ell$,
then Proposition~\ref{prop-QuadrParts} says that $c_{ij}$ is a trivial indeterminate,
in contradiction to the hypothesis. Hence $t_i$ has to be divisible by~$x_\ell$.
The term $t_i x_k/x_\ell$ cannot be in the border of~$\OO$, because then also
$t_i x_k$ would be outside~$\OO$, in contradiction to the hypothesis $t_i\in\OO^\circ$.
Therefore there exists an index $i'\in\{1,\dots,\mu\}$ such that $x_\ell t_{i'} = x_k t_i$.
By Proposition~\ref{prop-lin_cij}, we get $c_{ij} \sim c_{i'j'}$, in contradiction
to the hypothesis that $c_{ij}$ is basic.

Now we show~(d). For a contradiction, assume that $c_{ij}\in C$ is proper
and that every element in the cotangent equivalence class of~$c_{ij}$
is an interior indeterminate. By Proposition~\ref{prop-poscomp}, there exists
an index $\ell\in\{1,\dots,n\}$ such that the arrow degree $\deg_W(c_{ij})$
has a positive $\ell$-th component. Let $j'\in\{1,\dots,\nu\}$ be such that
the $\ell$-th component of $\log(b_{j'})$ is zero. (For instance, let $k\ne \ell$
and consider the unique term of the form $b_{j'} = x_k^N \in\partial\OO$ with 
$N\ge 1$.) Since the border is connected, we can find a sequence of border terms
$b_j = b_{j_1} \sim \cdots \sim b_{j_q} = b_{j'}$ such that $b_{j_p}$ and $b_{j_{p+1}}$
are next-door or across-the-rim neighbours for $p=1,\dots,q-1$.

Notice that it is not possible to find indeterminates $c_{i_1 j_1} \sim \cdots\sim c_{i_q j_q}$
since the equality of arrow degrees $\deg_W(c_{i_q j_q}) = \deg_W(c_{i_1 j_1})$ 
would imply that the $x_\ell$-degree of $t_{i_q}$ is negative. Let us try to construct
such a sequence of cotangent equivalences inductively. When we have found terms 
$t_{i_1},\dots,t_{i_p}\in\OO$ with the property that $c_{i_1 j_1} \sim \cdots\sim c_{i_p j_p}$
and try to find a term $t_{i_{p+1}}\in \OO$ such that 
$c_{i_p j_p} \sim c_{i_{p+1} j_{p+1}}$, three things can happen:
\begin{enumerate}
\item[(1)] A term $t_{i_{p+1}}$ of the desired kind exists in~$\OO^\circ$.
\item[(2)] A term $t_{i_{p+1}}$ of the desired kind exists in~$\OO^\nabla$.
\item[(3)] No term $t_{i_{p+1}}$ of the desired kind exists in~$\OO$ because
one of the components of $\log(t_{i_{p+1}})$ would be negative.
\end{enumerate}
In case~(1), we can continue our inductive construction for one further step.
By the hypothesis that $c_{ij}$ is not equivalent to a rim indeterminate,
case~(2) never occurs. Hence case~(3) has to happen for some $p\in\{ 1,\dots,q-1\}$.
In this case we have $\bar{c}_{i_p j_p} = 0$ by Proposition~\ref{prop-lin_cij},
and thus $\bar{c}_{ij} = 0$. Hence we have arrived at a contradiction to the hypothesis 
that~$c_{ij}$ is proper, and the proof is complete.

Finally, we prove (e). 
If $\#(Z) = \dim _K(\Lin_\M(I(\BO))$,  then Corollary~\ref{cor-cohSepPoly}.a implies that 
the $Z$-separating re-embedding of~$I(\BO)$ is optimal. Moreover,
Corollary~\ref{cor-cohSepPoly}.b implies that $Z = S_\sigma$ and hence the conclusion 
follows from Theorem~\ref{thm-shapeofSsigma}.a.
\end{proof}

Based on the preceding results and on the algorithm in~\cite{KSL}, Alg. 3.8, for
computing the cotangent equivalence classes, we can now check effectively
whether a given border basis scheme admits a $Z$-separating
embedding. Notice that we are excluding some trivial cases (namely $n=1$ and
$\mu \le 2$) in order to be able to apply~\cite{KSL}, Alg.~3.8, but these
cases can be dealt with easily by a direct computation.

\begin{algorithm}{\bf ($Z$-separating Re-embeddings Using Cotangent Equivalence Classes)}
\label{alg-compEmb}\\
Let $\OO = \{ t_1,\dots,t_\mu \}$ be an order ideal in~$\mathbb{T}^n$,
where $n\ge 2$ and $\mu\ge 3$,
let $\partial\OO = \{b_1,\dots,b_\nu \}$, and let
$C = \{ c_{ij} \mid i=1,\dots,\mu;\; j=1,\dots,\nu \}$ be the
set of coefficients of the generic $\OO$-border prebasis.
Consider the following sequence of instructions.
\begin{enumerate}
\item[(1)] Compute the set of natural generators of~$I(\BO)$.

\item[(2)] Compute the trivial cotangent equivalence class~$E_0$
and also the proper cotangent equivalence classes $E_1,\dots,E_q$.

\item[(3)] Let $S=\emptyset$.

\item[(4)] Turn each set $Z_0\cup Z_1\cup \cdots \cup Z_q$ such that 
$Z_0 \subseteq E_0$ and $Z_i \subset E_i$ for $i=1,\dots,q$ into a tuple~$Z$ 
and perform the following steps. 

\item[(5)] Using one of the methods given in Subsection~1.B, check whether the
ideal $I(\BO)$ is $Z$-separating. If it is, append~$Z$ to~$S$.

\item[(6)]  Continue with step (4) using the next tuple $Z$ until all tuples have been dealt with. 
Then return $S$ and stop.
\end{enumerate}
\begin{enumerate}
\item[(a)] This is an algorithm which computes the set~$S$ of all tuples~$Z$ 
of distinct indeterminates in~$C$ such that there exists a $Z$-sep\-a\-ra\-ting 
re-embedding of~$I(\BO)$.

\item[(b)] Assume that step~(4) is replaced by the following step. 

\begin{enumerate}
\item [(4')] Turn each set $E_0\cup E_1^\ast\cup \cdots \cup E_q^\ast$, 
where $E_i^\ast$ is obtained from~$E_i$ by deleting one element, 
into a tuple~$Z$ and perform the following steps. 
\end{enumerate}
Then the result is an algorithm which computes the set~$S$ of all tuples~$Z$ 
of distinct indeterminates in~$C$ such that there exists an optimal
$Z$-separating re-embedding of~$I(\BO)$.
\end{enumerate}
\end{algorithm}

\begin{proof}
All claims follow from Theorem~\ref{thm-shapeofSsigma} 
and Theorem~\ref{thm-classifycij}.a,b,e.  
\end{proof}

Notice that, if we use Proposition~\ref{prop-CheckZ} in step~(5), 
its execution may be very slow. The speed-ups given in Remark~\ref{rem-randomFforZ},
Algorithm~\ref{alg-checkZsepI}, and Remark~\ref{rem-CheckZvariants}
allow us to treat much larger examples, but do not guarantee to find a
tuple~$Z$ which yields a re-embedding of $I(\BO)$, even if it exists.

\medskip

The following example illustrates the results of this section.
The readers may also compare it to the example in~\cite{Hui1}, Remark 7.5.3.

\begin{example}\label{ex-Section5}
In $P =\QQ[x,y]$, consider the order ideal
$\OO = \{t_1,\dots,t_8\}$ given by $t_1=1$, $t_2=y$, $t_3=x$, $t_4=y^2$, $t_5=xy$, $t_6=x^2$, 
$t_7=y^3$, and $t_8 =xy^2$. Then we have $\partial\OO = \{b_1,\dots,b_5\}$ with $b_1=x^2y$, 
$b_2=x^3$, $b_3=y^4$, $b_4 = xy^3$, and $b_5=x^2y^2$. 
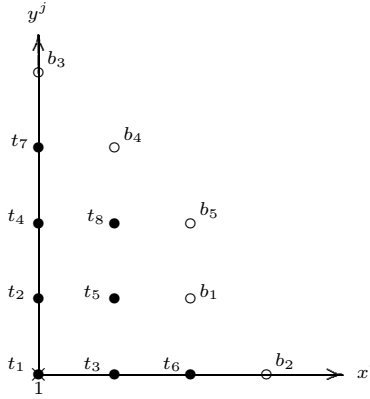
\begin{center}
	\makebox{\beginpicture
		\setcoordinatesystem units <1cm,1cm>
		\setplotarea x from 0 to 4, y from 0 to 4.5
		\axis left /
		\axis bottom /
		\arrow <2mm> [.2,.67] from  3.5 0  to 4 0
		\arrow <2mm> [.2,.67] from  0 4  to 0 4.5
		\put {$\scriptstyle x^i$} [lt] <0.5mm,0.8mm> at 4.1 0.1
		\put {$\scriptstyle y^j$} [rb] <1.7mm,0.7mm> at 0 4.6
		\put {$\bullet$} at 0 0
		\put {$\bullet$} at 1 0
		\put {$\bullet$} at 0 1
		\put {$\bullet$} at 1 1
	         \put {$\bullet$} at 1 2
		\put {$\bullet$} at 0 2
		\put {$\bullet$} at 0 3
		\put {$\bullet$} at 2 0
		\put {$\scriptstyle 1$} [lt] <-1mm,-1mm> at 0 0
		\put {$\scriptstyle t_1$} [rb] <-1.3mm,0.4mm> at 0 0
		\put {$\scriptstyle t_3$} [rb] <-1.3mm,0.4mm> at 1 0
		\put {$\scriptstyle t_2$} [rb] <-1.3mm,0mm> at 0 1
		\put {$\scriptstyle t_4$} [rb] <-1.3mm,0mm> at 0 2
		\put {$\scriptstyle t_5$} [rb] <-1.3mm,0mm> at 1 1
		\put {$\scriptstyle t_6$} [lb] <-4mm,0.4mm> at 2 0
		\put {$\scriptstyle t_7$} [lb] <-4mm,0mm> at 0 3
		\put {$\scriptstyle t_8$} [lb] <-4mm,0mm> at 1 2
		\put {$\scriptstyle b_1$} [lb] <0.8mm,0mm> at 2 1
		\put {$\scriptstyle b_2$} [lb] <0.8mm,0.8mm> at 3 0
	        \put {$\scriptstyle b_3$} [lb] <0.8mm,0.8mm> at 0 4
	        \put {$\scriptstyle b_4$} [lb] <0.8mm,0.8mm> at 1 3
	        \put {$\scriptstyle b_5$} [lb] <0.8mm,0.8mm> at  2 2
		\put {$\times$} at 0 0
	
		\put {$\circ$} at 3 0
		\put {$\circ$} at 2 2
		\put {$\circ$} at 2 1
		\put {$\circ$} at 1 3
		\put {$\circ$} at 0 4
		\endpicture}
\end{center}
Thus $\QQ[C] =\QQ[c_{11}, \dots, c_{85}]$ is a polynomial ring in~$40$ indeterminates. 
Notice that the dimension of~$\BO$ is $\dim(\BO) = \mu n = 16$ and that there are 
36 natural generators of the ideal $I(\BO)$. The linear parts of these generators are 
$$
\begin{array}{l}
c_{65},\;  c_{51} - c_{85},\;  c_{45},\;   c_{44},\;  c_{55},\;  c_{43} - c_{54},\;   
c_{42},\;  c_{41} - c_{75},\;  c_{52} - c_{75},\;  c_{35},  \cr 
c_{34},\;  c_{33},\;  c_{31},\;  c_{25},\;  c_{24},\;  c_{23},\;  c_{22},\;  c_{21},\;  
c_{32},\;  c_{15},\;  c_{14},\;  c_{13},\;  c_{12},\;   c_{11} 
\end{array}
$$
If we let $U$ be the union of the supports of these elements, we get 
$$
C\setminus U = \{c_{53},\,  c_{61},\,  c_{62},\,  c_{63},\,  c_{64},\,  c_{71},\,  
c_{72},\,  c_{73},\,  c_{74},\,  c_{81},\,  c_{82},\,  c_{83},\,  c_{84} \}
$$
which is exactly the set of basic indeterminates by Lemma~\ref{lem-equivclasses}.b.
Moreover, note that $C\setminus U$ is contained in the set 
of rim indeterminates (see Theorem~\ref{thm-classifycij}.c).

For the trivial cotangent equivalence class~$E_0$ and the 
proper cotangent equivalence classes $E_1,\dots,E_q$, we get
\begin{align*}
E_0 = & \{ c_{11},  c_{12},  c_{13},  c_{14},  c_{15},  c_{21},  c_{22},  c_{23},  
c_{24},  c_{25},  c_{31},  c_{32},  c_{33},  c_{34},  c_{35}, \\
 &\; c_{42},  c_{44},  c_{45},  c_{55},  c_{65} \}\\
E_1 = & \{c_{51},  c_{85}\},\qquad  E_2 = \{ c_{43}, c_{54}\},\qquad
E_3 = \{c_{41}, c_{52}, c_{75}\}
\end{align*}
Using $\sigma= \tt DegRevLex$, we obtain $E_1^\sigma = \{c_{51}\}$, $E_2^\sigma = \{c_{43}\}$,
and $E_3^\sigma = \{c_{41},  c_{52} \}$.
Next we compute the set $S_\sigma$ of the minimal generators of the leading term ideal 
of $I(\BO)_\lin$ and get
$$
\begin{array}{l}
S_\sigma = \{  c_{11},  c_{12},  c_{13},  c_{14},  c_{15},  c_{21},  c_{22},  c_{23},  
c_{24},  c_{25},  c_{31},  c_{32},  c_{33},  c_{34},  c_{35}, \cr  
\qquad \quad c_{41}, c_{42},  c_{43},  c_{44},  c_{45},  c_{51},  c_{52},  c_{55},  c_{65}
  \}
\end{array}
$$
which coincides with $E_0\cup E_1^\sigma \cup E_2^\sigma \cup E_3^\sigma$ 
(see the first claim of Theorem~\ref{thm-shapeofSsigma}.a).
Here we have $\#S_\sigma = 24$, $\#E_0 = 20$, $\#E_1 = 2$, $\#E_2 = 2$,  $\#E_3 = 3$,
and therefore 
$$
\#S_\sigma = \#E_0 + \#E_1 +\#E_2 +\#E_3 -3
$$
in accordance with the second claim of Theorem~\ref{thm-shapeofSsigma}.a.

The minimal sets of terms generating the sets in $\LTGFan(I(\BO)_\lin)$ are
\begin{align*}
	Z_1 = & E_0\cup \{ c_{51} \} \cup \{ c_{43} \} \cup \{ c_{41},\ c_{52}\} \\
	Z_2 = & E_0\cup \{ c_{51} \} \cup \{ c_{43} \} \cup \{ c_{41},\ c_{75}\} \\
	Z_3 = & E_0\cup \{ c_{85} \} \cup \{ c_{43} \} \cup \{ c_{41},\ c_{52}\} \\
	Z_4 = & E_0\cup \{ c_{85} \} \cup \{ c_{43} \} \cup \{ c_{41},\ c_{75}\} \\
	Z_5 = & E_0\cup \{ c_{51} \} \cup \{ c_{54} \} \cup \{ c_{41},\ c_{52}\} \\
	Z_6 = & E_0\cup \{ c_{51} \} \cup \{ c_{54} \} \cup \{ c_{41},\ c_{75}\} \\
	Z_7 = & E_0\cup \{ c_{85} \} \cup \{ c_{54} \} \cup \{ c_{41},\ c_{52}\} \\
	Z_8 = & E_0\cup \{ c_{85} \} \cup \{ c_{54} \} \cup \{ c_{41},\ c_{75}\} \\
	Z_9 = & E_0\cup \{ c_{51} \} \cup \{ c_{43} \} \cup \{ c_{52},\ c_{75}\} \\
	Z_{10} = & E_0\cup \{ c_{85} \} \cup \{ c_{43} \} \cup \{ c_{52},\ c_{75}\} \\
	Z_{11} = & E_0\cup \{ c_{51} \} \cup \{ c_{54} \} \cup \{ c_{52},\ c_{75}\} \\
	Z_{12} = & E_0\cup  \{ c_{85} \} \cup \{ c_{54} \} \cup \{ c_{52},\ c_{75}\}
\end{align*}
Thus we see that $\prod_{i=1}^3 \#E_i =12= \# (\LTGFan(I(\BO)_\lin)$
(see Theorem~\ref{thm-shapeofSsigma}.d).
As remarked before, $C\setminus U$ is the set of basic indeterminates, and we notice  
that $Z_i\cap (C\setminus U) =\emptyset$ for $i=1,\dots,12$ (see Theorem~\ref{thm-classifycij}.c).

The complements of the sets~$Z_i$ in~$C$ are 
$$
\begin{array}{ll}
Y_1&\!\!\!= \{ c_{53}, c_{54}, c_{61}, c_{62}, c_{63}, c_{64}, c_{71}, c_{72}, 
c_{73}, c_{74}, c_{75}, c_{81}, c_{82}, c_{83}, c_{84}, c_{85} \} \cr
Y_2&\!\!\!= \{  c_{52}, c_{53}, c_{54}, c_{61}, c_{62},  c_{63},  c_{64},  
c_{71},  c_{72},  c_{73},  c_{74},  c_{81},  c_{82},  c_{83},  c_{84},  
c_{85}  \} \cr
Y_3&\!\!\!= \{ c_{51},  c_{53},  c_{54},  c_{61},  c_{62},  c_{63},  c_{64},  
c_{71},  c_{72},  c_{73},  c_{74},  c_{75},  c_{81},  c_{82},  c_{83},  
c_{84} \} \cr
Y_4&\!\!\!= \{ c_{51},  c_{52},  c_{53},  c_{54},  c_{61},  c_{62},  c_{63},  
c_{64},  c_{71},  c_{72},  c_{73},  c_{74},  c_{81},  c_{82},  c_{83},  
c_{84} \} \cr
Y_5&\!\!\!= \{ c_{43},  c_{53},  c_{61},  c_{62},  c_{63},  c_{64},  c_{71},  
c_{72},  c_{73},  c_{74},  c_{75},  c_{81},  c_{82},  c_{83},  c_{84},  
c_{85} \} \cr
Y_6&\!\!\!=\{ c_{43},  c_{52},  c_{53},  c_{61},  c_{62},  c_{63},  c_{64},  
c_{71},  c_{72},  c_{73},  c_{74},  c_{81},  c_{82},  c_{83},  c_{84},  
c_{85}  \} \cr
Y_7&\!\!\!=\{ c_{43},  c_{51},  c_{53},  c_{61},  c_{62},  c_{63},  c_{64},  
c_{71},  c_{72},  c_{73},  c_{74},  c_{75},  c_{81},  c_{82},  c_{83},  
c_{84} \} \cr
Y_8&\!\!\!= \{ c_{43},  c_{51},  c_{52},  c_{53},  c_{61},  c_{62},  c_{63},  
c_{64},  c_{71},  c_{72},  c_{73},  c_{74},  c_{81},  c_{82},  c_{83},  
c_{84}  \} \cr
Y_9&\!\!\!=\{ c_{41},  c_{53},  c_{54},  c_{61},  c_{62},  c_{63},  c_{64},  
c_{71},  c_{72},  c_{73},  c_{74},  c_{81},  c_{82},  c_{83},  c_{84},  
c_{85}  \} \cr
Y_{10}&\!\!\!=\{ c_{41},  c_{51},  c_{53},  c_{54},  c_{61},  c_{62},  c_{63},  
c_{64},  c_{71},  c_{72},  c_{73},  c_{74},  c_{81},  c_{82},  c_{83},  
c_{84}  \} \cr
Y_{11}&\!\!\!=\{ c_{41},  c_{43},  c_{53},  c_{61},  c_{62},  c_{63},  c_{64},  
c_{71},  c_{72},  c_{73},  c_{74},  c_{81},  c_{82},  c_{83},  c_{84},  
c_{85} \} \cr
Y_{12}&\!\!\!= \{ c_{41},  c_{43},  c_{51},  c_{53},  c_{61},  c_{62},  
c_{63},  c_{64},  c_{71},  c_{72},  c_{73},  c_{74},  c_{81},  c_{82},  
c_{83},  c_{84} \}
\end{array}
$$
It is straightforward to verify that $C\setminus U$, the set of basic indeterminates, 
is contained in each set~$Y_i$ (see Theorem~\ref{thm-classifycij}.a).
Since we have $\#Z_i = \dim_\QQ(\Lin_\M(I(\BO)))$ for every~$i=1,\dots 12$, we are in the 
situation considered in Theorem~\ref{thm-classifycij}.e.
Using Algorithm~\ref{alg-compEmb}, we check that for each~$Z_i$ there exists 
an optimal $Z_i$-separating re-embedding of~$I(\BO)$.
Finally, for $i\in \{1,\dots,12\}$, we use $\#Z_i = 24$ and conclude 
that the set~$Z_i$ defines an isomorphism $\BO\cong \QQ[Y_i]$, 
where $\QQ[Y_i]$ is a polynomial ring having 16 indeterminates.
\end{example}

\bigskip\bigbreak
%
%

\section{Re-Embeddings of Planar Border Basis Schemes}
\label{Re-Embeddings of Planar Border Basis Schemes}

When we are given an order ideal~$\OO$ in two indeterminates, we say
that the corresponding scheme $\BO$ is a {\bf planar} border basis scheme.
This is the setting we examine more closely in the current section.

Before starting out with this task, we introduce some additional
terminology in the general setting. So, let $P=K[x_1,\dots,x_n]$ be once again 
a polynomial ring over a field~$K$, let $\OO = \{t_1,\dots,t_\mu\}$ 
be an order ideal in~$\mathbb{T}^n$ with border $\partial\OO = \{b_1,\dots,b_\nu\}$,
and let $C = \{ c_{ij} \mid i=1,\dots,\mu;\; j=1,\dots,\nu \}$ 
be a set of indeterminates which may serve as the coefficients 
of the generic $\OO$-border prebasis.

Recall that a term $b_j \in\partial\OO$ is called {\bf $x_\ell$-exposed} if it is of 
the form $b_j = x_\ell t_i$ with $t_i\in\OO$ (see Definition~\ref{def-exposedterm}).
The following definition extends this notion to the indeterminates in~$C$.

\begin{definition}\label{def-exposedIndets}
In the above setting, let $\ell \in \{1,\dots,n\}$.
\begin{enumerate}
\item[(a)] For $i\in \{1,\dots,\mu\}$ and $j\in \{1,\dots,\nu\}$, we say that the 
indeterminate $c_{ij} \in C$ is {\bf $x_\ell$-exposed} if
$x_\ell t_i \in \partial \OO$ and if there exists next-door neighbours $b_j, b_{j'}$
such that $b_{j'} = x_\ell b_j$ or across-the street neighbours $b_j, b_{j'}$
such that $x_k b_{j'} = x_\ell b_j$. 

\item[(b)] The union of the sets of all $x_\ell$-exposed indeterminates $c_{ij}$,
where $\ell\in \{1,\dots,n\}$, is called the set of all {\bf exposed indeterminates} in~$C$ 
and denoted by~$C^{\rm ex}$. The indeterminates in $C \setminus C^{\rm ex}$ 
are called {\bf non-exposed}.

\end{enumerate}
\end{definition}

\begin{remark}\label{rem:intnonexp}
Notice that all interior indeterminates of~$C$ are non-exposed, i.e., that
all exposed indeterminates are rim indeterminates.
\end{remark}

The following simple example illustrates this definition.

\begin{example}{\bf (The (2,1)-Box})\label{ex-12exposed}\\
In $P=\mathbb{Q}[x,y]$, consider the order ideal $\OO = \{t_1,t_2\}$, where 
$t_1=1$ and $t_2=x$. Then the border of~$\OO$ is $\partial\OO = \{b_1, b_2, b_3\}$, 
where $b_1=y$, $b_2=xy$, and $b_3 = x^2$.
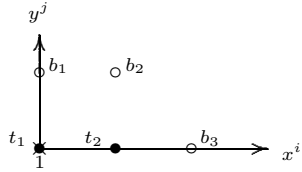
\begin{center}
	\makebox{\beginpicture
		\setcoordinatesystem units <1cm,1cm>
		\setplotarea x from 0 to 3, y from 0 to 1.5
		\axis left /
		\axis bottom /
		\arrow <2mm> [.2,.67] from  2.5 0  to 3 0
		\arrow <2mm> [.2,.67] from  0 1  to 0 1.5
		\put {$\scriptstyle x^i$} [lt] <0.5mm,0.8mm> at 3.1 0
		\put {$\scriptstyle y^j$} [rb] <1.7mm,0.7mm> at 0 1.6
		\put {$\bullet$} at 0 0
		\put {$\bullet$} at 1 0
		\put {$\scriptstyle 1$} [lt] <-1mm,-1mm> at 0 0
		\put {$\scriptstyle t_1$} [rb] <-1.3mm,0.4mm> at 0 0
		\put {$\scriptstyle t_2$} [rb] <-1.3mm,0.4mm> at 1 0
		\put {$\scriptstyle b_3$} [lb] <0.8mm,0.4mm> at 2 0
		\put {$\scriptstyle b_2$} [rb] <4.3mm,0mm> at 1 1
		\put {$\scriptstyle b_1$} [lb] <0.8mm,0mm> at 0 1
  		\put {$\circ$} at 0 1
		\put {$\circ$} at 1 1
		\put {$\times$} at 0 0
		\put {$\circ$} at 2 0
		\endpicture}
\end{center}
Here the bullets correspond to the terms in~$\OO$
and the circles to terms in~$\partial\OO$. 
The only index~$i$ such that $x t_i \in \partial\OO$ is $i=2$, where $x t_2 = b_3$.
There are two neighbour pairs, namely the next-door pair $b_2 = x b_1$ 
and the across-the-street pair $y b_3 = x b_2$. As they involve 
$x b_1$ and $x b_2$, we get $j\in \{1,2\}$.
Consequently, the $x$-exposed indeterminates $c_{ij}$ in~$C$
satisfy $i\in \{ 2\}$  and $j\in \{1,2\}$, i.e., they are $c_{21}, c_{22}$.
Similarly, we see that the $y$-exposed indeterminates satisfy $i\in \{1,2 \}$ 
and $j \in \{3\}$, i.e., they are $c_{13}, c_{23}$.
\end{example}

The following proposition provides two characterizations of exposed
indeterminates which motivate the above definition.

\begin{proposition}\label{prop:CharExposed}
Assuming the above setting, let $G = \{g_1,\dots,g_\nu\}$ be the generic $\OO$-border 
prebasis, where $g_j \;= \; b_j - \sum_{i=1}^\mu c_{ij} t_i$ for $j=1\dots,\nu$.
Then, for all $i\in \{1,\dots,\mu\}$ and $j\in \{1,\dots,\nu\}$, 
the following conditions are equivalent.
\begin{enumerate}
\item[(a)] There exists an index $\ell \in \{1,\dots,n\}$ such that
the indeterminate $c_{ij} \in C$ is $x_\ell$-exposed.

\item[(b)] There exists an index $j'\in \{1,\dots,\nu\}$ such that $b_j, b_{j'}$
is a next-door neighbour pair with $b_{j'} = x_\ell b_{j}$ or an across-the-street
neighbour pair with $x_k b_{j'} = x_\ell b_j$ and $k\in \{1,\dots,n\}$, and
such that if we write 
$$
x_\ell g_j \;=\; x_\ell b_j  + c_{i_1 j} g_{j_1} + \cdots + c_{i_r j} g_{j_r} 
+ \tsum_{\lambda=1}^\mu f_\lambda t_\lambda
$$
where $i_1,\dots,i_r\in \{1,\dots,\mu\}$, where $j_1,\dots,j_r\in\{1,\dots,\nu\}$, 
and where $f_\lambda \in K[C]$, then we have $i \in \{i_1, \dots, i_r \}$.

\item[(c)] There exists an index $j'\in \{1,\dots,\nu\}$ such that $b_j, b_{j'}$
is a next-door neighbour pair with $b_{j'} = x_\ell b_j$ or an across-the-street
neighbour pair with $x_k b_{j'} = x_\ell b_j$ and $k\in \{1,\dots,n\}$, and
such that the lifting of the neighbour syzygy $e_{j'} - x_\ell e_j$ resp.\ 
$x_k e_{j'} - x_\ell e_j$ is of the form 
$$
x_k^\epsilon e_{j'} - x_\ell e_j - \tsum_{\kappa=1}^\nu \bar{f}_\kappa  e_{\kappa}
\in \Syz_{B_\OO}( \bar{g}_1, \dots, \bar{g}_\nu )
$$
where $\epsilon\in \{0,1\}$ and $\bar{f}_\kappa \in B_\OO$ are the residue classes 
of monomials or binomials in $K[C]$ such that one of their supports contains $c_{ij}$.

\end{enumerate}
\end{proposition}

\begin{proof}
To show (a)$\Rightarrow$(b), we write $x_\ell t_i = b_{j''}$ with 
$j'' \in \{1,\dots,\nu\}$. Clearly, there exists a neighbour pair $b_j, b_{j'}$
such that $x_k^\epsilon b_{j'} = x_\ell b_j$ for some $j, j' \in \{1,\dots,\nu\}$
and $k\in \{1,\dots,n\}$ and $\epsilon \in \{0,1\}$. (For instance, we can take the
pure $x_\ell$-power $b_j = x_\ell^\kappa$ in~$\partial\OO$.) Then we have
$$
x_\ell g_j \;=\; x_\ell b_j - \tsum_{\lambda=1}^\mu  c_{\lambda j} x_\ell t_\lambda
$$
and if we denote those terms $x_\ell t_{\lambda}$ which are in $\partial\OO$ by
$x_\ell t_{i_1} = b_{j_1}$, $\dots$, $x_\ell t_{i_r} = b_{j_r}$ then
the fact that $x_\ell t_i \in \partial\OO$ shows $i \in \{i_1,\dots,i_r\}$.

Conversely, the condition $f_\lambda \in K[C]$ implies that the terms in the
support of $x_\ell (g_j-b_j)$ 
which are contained in $K[C] \cdot \partial \OO$ are of the form $c_{i_\kappa j} b_{j_\kappa}$
with $\kappa \in \{1,\dots,r\}$. Therefore $b_{j_\kappa} = x_\ell t_{i_\kappa}$
for $\kappa = 1,\dots,r$ and $c_{ij} \in \{ c_{i_1 j_1}, \dots, c_{i_r j_r} \}$ 
is $x_\ell$-exposed.

The equivalence of~(b) and~(c) follows from the fact that the representation of
$x_k^\epsilon g_{j'} - x_\ell g_j$ which follows from~(b) is exactly the lifting
of the neighbour syzygy $x_k^\epsilon e_{j'} - x_\ell e_j$ in 
$\Syz_{B_\OO}( \bar{g}_1, \dots, \bar{g}_\nu )$ by~\cite{KR3}, Prop. 4.1.
\end{proof}

Let us reconsider Example~\ref{ex-12exposed} to illustrate this proposition.

\begin{example}{\bf (The (2,1)-Box Continued})\label{ex-12exposed-continued}\\
In the setting of Example~\ref{ex-12exposed}, the generic $\OO$-border prebasis is
$$
\begin{array}{lllll}
g_1 &=& y\ \; -c_{11} -c_{21} x &=& b_1 -c_{11} -c_{21} x \cr
g_2 &=& xy -c_{12} -c_{22} x &=& b_2 -c_{12} -c_{22} x \cr
g_3 &=& x^2\,   -c_{13} -c_{23} x &=& b_3 -c_{13} -c_{23} x 
\end{array}
$$
The next-door neighbour pair $b_2 = x b_1$ and the
across-the-rim-neigh\-bour pair $x b_2 = y b_1$ yield
$$
\begin{array}{lll}
x g_1&= &x b_1 -  c_{11}  x  - c_{21}x^2  =   x b_1 -c_{11} x  - c_{21}  (g_3 + c_{13} +c_{23} x)     \cr
&= & x b_1 -c_{21}g_3  -c_{13}c_{21}-(c_{11} +c_{21}c_{23})x \cr
x g_2 &= &x b_2 - c_{12} x -c_{22} x^2 =  x b_2 -c_{12} x  -c_{22} (g_3 + c_{13} +c_{23} x) \cr
&= & x b_2  -c_{22} g_3  -(c_{12} +c_{13} c_{22})  +c_{22}c_{23} x
\end{array}
$$
Thus the set of $x$-exposed indeterminates is $\{c_{21},\, c_{22}\}$.
A similar computation shows that the set of $y$-exposed indeterminates is $\{c_{13},\, c_{23}\}$.
Hence the set of exposed indeterminates is 
$C^{\rm ex} = \{ c_{13}, c_{21}, c_{22}, c_{23} \}$.

The liftings of the two neighbour syzygies $x e_1 - e_2$ and $x e_2 -y e_3$ are
$$
x e_1 - e_2 + \bar{c}_{21}e_3  \qquad\hbox{\rm and}\qquad 
x e_2 - y e_3 - \bar{c}_{13}e_1  - \bar{c}_{23} e_2  + \bar{c}_{22} e_3
$$
This yields the same sets of exposed indeterminates once again.

Notice that all six indeterminates in~$C$ are rim-indeterminates. 
Thus the inclusion $C^{\rm ex} \subseteq C^\nabla$ of Remark~\ref{rem:intnonexp} 
can be strict. 
\end{example}

\medskip
In the remainder of this section we let $n=2$, i.e., we consider planar border 
basis schemes. Recall that we usually write~$x$ for $x_1$ and~$y$ for~$x_2$ in this case.
The following additional terminology will help us formulate our next algorithm.

\begin{definition} 
Let $\OO = \{ t_1,\dots,t_\mu \}$ be an order ideal in~$\mathbb{T}^2$ with border
$\partial\OO = \{b_1,\dots,b_\nu \}$.
\begin{enumerate}
\item[(a)] A term is called {\bf $x$-free} (resp. {\bf $y$-free}) if~$x$ (resp.~$y$)
does not divide it.

\item[(b)] A term $b_j\in\partial\OO$ is said
to have an {\bf up-neighbour} $b_{j'}\in\partial\OO$ if
$b_{j'} = x_i b_j$ for some indeterminate $x_i \in X$. In this case~$b_j$
is also called a {\bf down-neighbour} of~$b_{j'}$.

\item[(c)] A sequence of border terms $b_{j_1},\dots,b_{j_k}$ 
is said to form a {\bf plateau} if the pairs
$(b_{j_1},b_{j_2})$, $\dots$, $(b_{j_{k-1}}, b_{j_k})$ form a maximal
chain of across-the-rim neighbour pairs with $x b_{j_\ell} = y b_{j_{\ell+1}}$
for $\ell = 1,\dots,k-1$, and if $b_{j_1}, b_{j_k}$ have no
up-neighbours.

\item[(d)] Let $b_{j_1},\dots,b_{j_k}$ be a plateau, 
and let $b_{j_1} = x b_{j'_1}$ with $b_{j'_1}\in\partial\OO$.
Let $b_{j'_2},\dots,b_{j'_m}$ be a maximal sequence of border terms
such that $(b_{j'_{\ell+1}}, b_{j'_\ell})$ are AR or ND neighbours
in $x$-direction, i.e., such that $x b_{j'_{\ell+1}} \in 
\{ b_{j'_\ell}, y b_{j'_\ell} \}$. Then the sequence $b_{j'_1},\dots,
b_{j'_m}$ is called the {\bf $x$-leg} of the  plateau $b_{j_1},\dots,b_{j_k}$.
The {\bf $y$-leg} of the plateau $b_{j_1},\dots,b_{j_k}$ is defined
analogously.
\end{enumerate}
\end{definition}

Let us illustrate these definitions using Example~\ref{ex-122}.

\begin{example}\label{ex-122-2} 
In $P=\mathbb{Q}[x,y]$, look at the order ideal $\OO = \{t_1,\dots,t_5\}$, where 
$t_1=1$, $t_2=y$, $t_3=x$, $t_4=y^2$, and $t_5=xy$. Then the border of~$\OO$
is $\partial\OO = \{b_1,\dots,b_4\}$, where $b_1=x^2$, $b_2=y^3$, 
$b_3=xy^2$, and $b_4=x^2y$ (see Example~\ref{ex-122}).
There is only one border term with an up-neighbour, namely $b_1=x^2$.
The terms $(b_2,b_3,b_4)$ form a plateau which has
the $y$-leg $(b_1)$ and no $x$-leg.
\end{example}

For planar border basis schemes, our next goal is to reprove and make explicit 
some beautiful results from~\cite{Hui1}. 
More precisely, let~$Z$ be the tuple of non-exposed 
indeterminates $C \setminus C^{\rm ex}$, ordered arbitrarily. To show that
there exists a $Z$-separating re-embedding of~$I(\BO)$ which uses the
natural set of generators of~$I(\BO)$, we proceed as follows:
\begin{enumerate}
\item[(1)] Define a set of weights for all $c_{ij}$ such that
$\wt(c_{ij}) = 0$ if $c_{ij}$ is an exposed indeterminate, such that
$\wt(c_{ij}) > 0$ if $c_{ij}$ is a non-exposed indeterminate,
and such that for each non-exposed indeterminate $c_{ij}$ there
exists a natural generator~$f$ of~$I(\BO)$ for which $c_{ij}$
is the unique term of highest weight in~$\Supp(f)$.

\item[(2)] Define an elimination ordering~$\sigma$ for~$Z$
given by a matrix whose first row consists of the weights assigned
in Step~(1).
\end{enumerate}

Then, by Prop.~\ref{prop-cohSepPoly}, we have a $Z$-separating
re-embedding of~$I(\BO)$. In this way, we have eliminated all non-exposed
indeterminates. To execute Step~(1), we exhibit the following explicit algorithm.

\begin{algorithm}{\bf (Weight Assignment Algorithm)}\label{alg-WeightAssign}\\
In the above setting, let $\delta = \max \{ \deg(t_i) \mid t_i\in \OO\}$.
Consider the following instructions.
\begin{enumerate}
\item[(1)] Let $S=\emptyset$ and let $\wt(c_{ij}) = 0$ for all $c_{ij}\in C^{\rm ex}$.

\item[(2)] For $d=\delta, \delta-1, \dots,0$, execute the loop defined by 
the following steps (3)-(7).

\item[(3)] Let $C_d$ be the set of all indeterminates~$c_{ij}$ in~$C$
such that $\deg(t_i)=d$ and such that $c_{ij}$ has not yet been assigned
a weight.

\item[(4)] Let $c_{ij} \in C_d$ be an indeterminate for which~$b_j$
has an up-neighbour $b_{j'} = x_k b_j$.
Let $i'\in\{1,\dots,\mu\}$ such that $t_{i'} = x_k t_i$ and let
$b_{j_1},\dots,b_{j_m}$ be the $x_k$-exposed border terms.
Then we let
$$
\wt(c_{ij}) \;=\; 1 + \wt(c_{i'j'}) + \wt(c_{i' j_1}) + \dots + \wt(c_{i' j_m})
$$

\item[(5)] Let $c_{ij} \in C_d$ be an indeterminate for which~$b_j$
is an element of a plateau $b_{j_1},\dots,b_{j_k}$. To define
the weight of~$c_{ij}$, we distinguish several cases.
In steps (6a) to (6e) we treat the cases where $\deg_W(c_{ij})$
has a positive $x$-component. In steps (7a) to (7e) we treat the cases
where $\deg_W(c_{ij})$ has a non-positive $x$-component and a positive 
$y$-component.

\item[(6a)] Suppose that $j=j_\ell$ for some $\ell\in\{2,\dots,k\}$
and $t_i$ is $x$-free. Moreover, let $m\in\{1,\dots,\mu\} $ be such that
$t_m = y t_i$. Then we set
$$
\wt(c_{ij}) \;=\; 1 + \wt(c_{m1}) + \cdots + \wt(c_{m\nu})
$$

\item[(6b)] Suppose there exists $\ell\in\{2,\dots,k\}$
such that $j=j_\ell$ and such that we have $t_i = x^\alpha y^\beta$
with $\alpha\le \ell -2$. Starting with $\ell=2$, the case treated in
step~(6a), we define $\wt(c_{ij})$ by induction on~$\ell$.
Let $i'\in\{1,\dots,\mu\}$ such that $c_{i' j_{\ell-1}} \sim c_{i j_\ell}= c_{ij}$
and let $m\in\{1,\dots,\mu\}$ such that $t_m = y t_i$. Then we define
$$
\wt( c_{ij} ) \;=\; 1 + \wt(c_{i' j_{\ell-1}}) + \wt(c_{m1}) + \cdots + \wt(c_{m\nu})
$$

\item[(6c)] Suppose that $j=j_1$ and that the plateau has an $x$-leg at~$b_{j_1}$
starting with $b_{j_1} = x b_{j'_1}$. Moreover, suppose that $t_i$ is $x$-free.
Let $b_{j''_1},\dots,b_{j''_m}$ be the $x$-exposed border terms.
Then we introduce a new indeterminate $p_{ij}$, we set $\wt(c_{ij}) = p_{ij}$
and we append the following inequalities to~$S$: 
$p_{ij} > \wt(c_{i j''_1})$, $\dots$, $p_{ij} > \wt(c_{i j''_m})$.

\item[(6d)] [Going down the $x$-leg]\\
Suppose that $j=j_1$ and that $b_{j'_1},\dots,b_{j'_\kappa}$ is the $x$-leg
of the plateau. Then we introduce a new indeterminate $p_{ij}$ and
we set $\wt(c_{ij}) = p_{ij}$. Let $\ell\le \kappa$ and let 
$i'_1,\dots,i'_\ell \in\{1,\dots,\mu\}$ be those indices for which 
$c_{i'_1 j'_1} \sim \cdots \sim c_{i'_\ell j'_\ell}$. Then we set
$$
\wt(c_{i'_1 j'_1}) \;=\; p_{ij}-1,\;\; \wt(c_{i'_2 j'_2}) \;=\; p_{ij}-2,\;\dots,\;
\wt(c_{i'_\ell j'_\ell}) \;=\; p_{ij}-\ell
$$
and we add inequalities of the form $p_{ij} \ge N$ to~$S$ 
to express the fact that $c_{i'_\lambda j'_\lambda}$ is the term of the 
highest weight in the polynomial $\ND(j'_{\lambda+1}, j'_\lambda)_{i'_\lambda}$
resp. $\AR(j'_\lambda, j'_{\lambda+1})_{m_\lambda}$.

\item[(6e)] Suppose that $j=j_\ell$ for some $\ell\in\{ 1,\dots,k\}$
and suppose that $t_i = x^\alpha y^\beta$ with $\alpha \ge \ell-1$.
Starting with the case $\ell=1$, treated in
step~(6c), we define $\wt(c_{ij})$ by induction on~$\ell$.
Choose $i'\in\{1,\dots,\mu\}$ such that we have $c_{i' j_{\ell-1}} \sim c_{i j_\ell}= c_{ij}$,
and let $m\in\{1,\dots,\mu\}$ such that $t_m = y t_i$. Then we define
$$
\wt( c_{ij} ) \;=\; 1 + \wt(c_{i' j_{\ell-1}}) + \wt(c_{m1}) + \cdots + \wt(c_{m\nu})
$$

\item[(7a)] Suppose that $j=j_\ell$ for some $\ell\in\{1,\dots,k-1\}$
and $t_i$ is $y$-free. Let $m\in\{1,\dots,\mu\}$ be such that
$t_m = x t_i$. Then we set
$$
\wt(c_{ij}) \;=\; 1 + \wt(c_{m1}) + \cdots + \wt(c_{m\nu})
$$

\item[(7b)] Suppose there exists $\ell\in\{1,\dots,k-1\}$
such that $j=j_\ell$ and such that $t_i = x^\alpha y^\beta$
with $\beta\le k-1-\ell$. Starting with the case $\ell=k-1$, treated in
step~(7a), we define $\wt(c_{ij})$ by downward induction on~$\ell$.
Moreover, let $i'\in\{1,\dots,\mu\}$ be such that $c_{i' j_{\ell+1}} \sim c_{i j_\ell}= c_{ij}$
and let $m\in\{1,\dots,\mu\}$ such that $t_m = x t_i$. Then we define
$$
\wt( c_{ij} ) \;=\; 1 + \wt(c_{i' j_{\ell+1}}) + \wt(c_{m1}) + \cdots + \wt(c_{m\nu})
$$

\item[(7c)] Suppose that $j=j_k$ and that the plateau has a $y$-leg at~$b_{j_k}$
starting with $b_{j_k} = y b_{j'_1}$. Moreover, suppose that $t_i$ is $y$-free.
Let $b_{j''_1},\dots,b_{j''_m}$ be the $y$-exposed border terms.
Then we introduce a new indeterminate $p_{ij}$, we set $\wt(c_{ij}) = p_{ij}$
and append the following inequalities to~$S$: 
$p_{ij} > \wt(c_{i j''_1})$, $\dots$, $p_{ij} > \wt(c_{i j''_m})$.

\item[(7d)] [Going down the $y$-leg]\\
Suppose that $j=j_k$ and that $b_{j'_1},\dots,b_{j'_\kappa}$ is the $y$-leg
of the plateau. Then we introduce a new indeterminate $p_{ij}$ and
we set $\wt(c_{ij}) = p_{ij}$. Let $\ell\le \kappa$ and let 
$i'_1,\dots,i'_\ell \in\{1,\dots,\mu\}$ be those indices for which 
$c_{i'_1 j'_1} \sim \cdots \sim c_{i'_\ell j'_\ell}$. Then we set
$$
\wt(c_{i'_1 j'_1}) \;=\; p_{ij}-1,\;\; \wt(c_{i'_2 j'_2}) \;=\; p_{ij}-2,\;\dots,\;
\wt(c_{i'_\ell j'_\ell}) \;=\; p_{ij}-\ell
$$
and we add inequalities of the form $p_{ij} \ge N$ to~$S$ 
to express the fact that $c_{i'_\lambda j'_\lambda}$ is the term of the 
highest weight in the polynomial $\ND(j'_\lambda, j'_{\lambda+1})_{i'_\lambda}$
resp. $\AR(j'_\lambda, j'_{\lambda+1})_{m_\lambda}$.

\item[(7e)] Suppose that $j=j_\ell$ for some $\ell\in\{ 1,\dots,k\}$
and suppose that $t_i = x^\alpha y^\beta$ with $\beta \ge k-\ell$.
Starting with the case $\ell=k$ which was treated in
step~(7c), we define $\wt(c_{ij})$ by induction on~$\ell$.
Choose $i'\in\{1,\dots,\mu\}$ such that $c_{i' j_{\ell+1}} \sim c_{i j_\ell}= c_{ij}$,
and let $m\in\{1,\dots,\mu\}$ such that $t_m = x t_i$. Then we define
$$
\wt( c_{ij} ) \;=\; 1 + \wt(c_{i' j_{\ell+1}}) + \wt(c_{m1}) + \cdots + \wt(c_{m\nu})
$$

\item[(8)] Finally, consider the set of inequalities~$S$. Choose the numbers $p_{ij}$
large enough so that all inequalities are satisfied. In this way, all 
indeterminates~$c_{ij}$ will be assigned explicit weights.
\end{enumerate}
This is an algorithm which assigns to each $c_{ij}\in C$
a weight $\wt(c_{ij})\ge 0$ such that the following conditions are satisfied:
\begin{enumerate}
\item[(a)] All exposed indeterminates $c_{ij}\in C^\nabla$ satisfy
$\wt(c_{ij}) = 0$.

\item[(b)] All non-exposed indeterminates $c_{ij} \in C^\circ$ satisfy
$\wt(c_{ij}) > 0$.

\item[(c)] For each non-exposed indeterminate $c_{ij}$ there
exists a natural generator~$f$ of the ideal~$I(\BO)$ for which $c_{ij}$
is the unique term of highest weight in~$\Supp(f)$.
\end{enumerate}
\end{algorithm} 

\begin{proof}
Since the algorithm is clearly finite, it suffices to show that
the assigned weights lead to the desired terms of highest weight
in specifically chosen natural generators of~$I(\BO)$.
\begin{enumerate}
\item[(a)] The definition of the weights in step~(4) implies that $c_{ij}$
is the highest weight term in $\ND(j,j')_{i'}$.

\item[(b)] In case (6a), the term $c_{ij} = c_{i j_\ell}$ is the linear part
and the highest weight term of the polynomial $\AR(j_\ell,j_{\ell-1})_m$.

\item[(c)] In case (6b), the linear part of $\AR(j_\ell,j_{\ell-1})_m$ is
(up to sign) the polynomial $c_{ij} - c_{i' j_{\ell-1}}$ and the highest
weight term is $c_{ij}$.

\item[(d)] In case (6c), the term $c_{ij}$ is the linear part of 
$\ND(j'_1,j_1)_i$ and its weight~$p_{ij}$ is the highest weight of
any term in this polynomial by the given inequalities for~$p_{ij}$.

\item[(e)] In case (6d), when we move from~$b_{j'_\lambda}$ to
$b_{j'_{\lambda+1}}$ while going down the $x$-leg, we go down a next-door
neighbour pair using the polynomial $\ND(j'_{\lambda+1}, j'_\lambda)_{i'_\lambda}$
or an across-the-rim neighbour pair using the polynomial 
$\AR(j'_\lambda,j'_{\lambda+1})_{m_\lambda}$ and weight given 
to~$c_{i_\lambda j_\lambda}$ as well as the given inequalities for $p_{ij}$
make sure that the term of highest weight is $c_{i_\lambda j_\lambda}$
in both cases.

\item[(f)] In case (6e), we are again using the linear part 
of $\AR(j_\ell,j_{\ell-1})_m$ to move along the plateau. It is
(up to sign) the polynomial $c_{ij} - c_{i' j_{\ell-1}}$ and the highest
weight term is~$c_{ij}$ by the given definition.

\item[(g)] In steps (7a)-(7e) analogous arguments apply.
\end{enumerate}
Finally, we note that the required set of inequalities for the
weights $p_{ij}$ can always be solved: when we move into and down an $x$-leg in
steps (6c) and (6d), the weight $p_{ij}$ has to be at least as high as 
some weights $\wt(c_{i'j'})$ with $i'=i$ or $b_{j'}$ in the $x$-tail of the
plateau. On the other hand, the inequalities in steps (7c) and (7d) involve
indices~$i$ which correspond to the other end of the plateau or its $y$-leg.
Therefore we will never get inequalities which interfere with each other
and by choosing $p_{ij}$ large enough all inequalities will be satisfied.
\end{proof}

Based on this algorithm and our strategy, we have now shown the following
general result.

\begin{theorem}{\bf (Elimination of Non-Exposed Indeterminates)}\label{thm-elimNex}\\
Let $K$ be a field, let $P=K[x,y]$,
let $\OO = \{ t_1,\dots,t_\mu \}$ be an order ideal in~$\mathbb{T}^2$, 
let $\partial\OO = \{b_1,\dots,b_\nu \}$ be its border, and let
$C = \{ c_{ij} \mid i=1,\dots,\mu;\; j=1,\dots,\nu \}$  
denote the set of coefficients of the generic $\OO$-border prebasis.
Let $Z=(z_1,\dots,z_s)$ be a tuple consisting of the non-exposed 
indeterminates $C\setminus C^{\rm ex}$.

Using Algorithm~\ref{alg-WeightAssign}, determine
weights $p_{ij}\in\mathbb{N}$ for all $c_{ij}\in C$, and let~$\sigma$
be an elimination ordering for~$Z$ defined by putting the weights~$p_{ij}$
into the first row of a matrix of size $\mu\nu$ and extending it to
a term ordering matrix.

Then the natural generators $f_1,\dots,f_s$ of~$I(\BO)$ given in the proof of 
Algorithm~\ref{alg-WeightAssign} satisfy $\LT_\sigma(f_i) = z_i$
for $i=1,\dots,s$ and define a $Z$-separating re-embedding of~$I(\BO)$.
In other words, the non-exposed indeterminates have been eliminated. 
\end{theorem}

\begin{proof}
This follows from the algorithm and Proposition~\ref{prop-cohSepPoly}.
\end{proof}

This theorem has the following immediate consequences.

\begin{corollary}\label{cor-exposedgenerate}
In the setting of the theorem, the following claims hold.
\begin{enumerate}
\item[(a)] The residue classes modulo $I(\BO)$ of the exposed indeterminates generate 
the $K$-algebra $B_\OO$.

\item[(b)] Every basic indeterminate is exposed and cannot be contained in a tuple~$Z$
for which $I(\BO)$ has a $Z$-separating re-embedding.

\end{enumerate}
\end{corollary}

\begin{proof} 
Claim (a) follows immediately from the theorem, and claim~(b) 
is a consequence of~(a) and Lemma~\ref{lem-equivclasses}.d.
\end{proof}

Let us compare Theorem~\ref{thm-elimNex} and its corollary to Prop.~7.1.2 in~\cite{Hui1}.

\begin{remark}\label{rem-Hui1approach}
In~\cite{Hui1}, Sections~6 and~7, the author uses an ingenious argument
based on the Hilbert-Burch Theorem to show that the residue classes modulo $I(\BO)$
of the exposed indeterminates form a $K$-algebra system of generators of~$B_\OO$.
In particular, using the maximal minors of the matrix~$M$ whose rows
represent the liftings of the ND and AR neighbour syzygies to 
$\Syz_{B_\OO}(\bar{g}_1,\dots, \bar{g}_\nu)$, where $G=\{ g_1,\dots,g_\nu\}$
is the generic $\OO$-border prebasis, as for instance in Example~\ref{ex-box22-weights},
one can find explicit representations of the non-exposed indeterminates~$c_{ij}$ 
as polynomial expressions $c_{ij} = f_{ij}(C^{\rm ex})$ in the exposed indeterminates. 
The advantage of this method is that those expressions form
a Gr\"obner basis of the ideal $\langle c_{ij} - f_{ij}(C^{\rm ex})\rangle$,
and this Gr\"obner basis is obtained by expanding the minors of~$M$.

On the other hand, the advantage of the approach based on the 
Weight Assignment Algorithm is that it finds a subset of the natural set of 
generators of~$I(\BO)$ which yields the re-embedding, and these tend to be rather 
sparse polynomials. However, it is clear that both approaches pay the price that 
the polynomials in the corresponding reduced Gr\"obner basis, which is necessary to 
compute the re-embedding, may have large supports.
\end{remark}

Of particular interest are order ideals $\OO$ such that the non-exposed indeterminates
produce an isomorphism of~$B_\OO$ with a polynomial ring. They can be described
as follows.

\begin{remark}\label{rem-optZsepPlanar}
In the above setting, suppose that the $Z$-separating re-embedding
of~$I(\BO)$ is optimal. By Corollary~\ref{cor-affinescheme}, this implies that 
the border basis scheme~$\BO$ is isomorphic to $\mathbb{A}^{2\mu}_K$.
Moreover, we have $\# Z  = \# C - 2\mu = \mu \nu - 2\mu$.

The special cases in which the $Z$-separating re-embedding is optimal for 
the tuple~$Z$ of non-exposed indeterminates, i.e., the cases where $\# C^{\rm ex} = 2\mu$, 
are classified in~\cite{Hui1}, Prop.~7.3.1.
They occur if and only if the order ideal~$\OO$ has the ``sawtooth'' form described there.
However, as already mentioned in~\cite{Hui1}, Remark 7.5.3.(1), there are cases 
where~$\BO$ is an affine cell, but $\OO$ does not have the ``sawtooth'' form (see
also Example~\ref{ex-Section5}).
\end{remark}

Let us apply the Weight Assignment Algorithm~\ref{alg-WeightAssign} and
Huibregtse's method of Remark~\ref{rem-Hui1approach} in the setting of
Example~\ref{ex-box22} and compare them.

\begin{example}{\bf (The (2,2)-Box Continued)}\label{ex-box22-weights}\\
Let $K$ be a field, let $P=K[x,y]$, and let $\OO = \{t_1,\, t_2,\, t_3,\, t_4\}$
be the order ideal given by $t_1=1$, $t_2=y$, $t_3=x$, and $t_4=xy$.
Then we have $\partial\OO = \{b_1,\, b_2,\, b_3,\, b_4\}$, where
$b_1=y^2$, $b_2=x^2$, $b_3=xy^2$, and $b_4=x^2y$.

In this situation there are three neighbour pairs: next-door neighbours
$(b_1,b_3)$ with $b_3 = x_1 b_1$, next-door neighbours $(b_2,b_4)$ with $b_4 = y b_2$,
and across-the-rim neighbours $(b_3,b_4)$ with $xb_3 = yb_4$. Using 
Definition~\ref{def-exposedIndets}, we conclude that 
$$
C^{\rm ex} \;=\; \{ c_{22},\, c_{24},\, c_{31},\, c_{33},\, c_{41},\, c_{42},\, c_{43},\, c_{44}\}
$$
is the set of exposed indeterminates. Let us follow the steps of 
Algorithm~\ref{alg-WeightAssign} and assign weights to all indeterminates $c_{ij}$.
\begin{enumerate}
\item[(1)] Let $\wt(c_{22}) =\wt(c_{24}) = \wt(c_{31})= \wt(c_{33}) = \wt(c_{41}) =
\wt(c_{42})\allowbreak = \wt(c_{43}) = \wt(c_{44})=0$.

\item[(2)] For $d=2$, we have $C_2=\emptyset$. Hence we let $d=1$
and $C_1 = \{c_{21}, c_{23}, c_{32}, c_{34} \}$ next.

\item[(4)] For the indeterminates~$c_{21}$ and~$c_{32}$, the border terms have
up-neighbours and we get
\begin{align*}
\wt(c_{21}) &\;=\; 1 + \wt(c_{42}) + \wt(c_{43}) + \wt(c_{44}) \;=\; 1\\
\wt(c_{32}) &\;=\; 1 + \wt(c_{41}) + \wt(c_{43}) + \wt(c_{44}) \;=\;  1
\end{align*}

\item[(6c)] For the indeterminate~$c_{23}$, we let $\wt(c_{23})=p_{23}$
and put the inequalities $p_{23} > \wt(c_{22}) = 0$ as well as
$p_{23} > \wt(c_{24}) = 0$ into~$S$.

\item[(6e)] We get $\wt(c_{34}) = 1 + \sum_{j=1}^4 \wt(c_{4j}) + \wt(c_{23}) = p_{23}+1$.

\item[(2)] Finally, we move to the case $d=0$ and $C_0 = \{c_{11},c_{12},c_{13},c_{14}\}$.

\item[(4)] For $c_{11}$ and $c_{12}$ we can use the up-neighbours of~$b_1$ and~$b_2$
and get
\begin{align*}
\wt(c_{11}) &\;=\; 1 + \wt(c_{33}) + \wt(c_{32}) + \wt(c_{34}) \;=\; p_{23} + 3\\
\wt(c_{12}) &\;=\; 1 + \wt(c_{24}) + \wt(c_{21}) + \wt(c_{23}) \;=\; p_{23} + 2
\end{align*}

\item[(6a)] Since $t_1$ is $x$-free, we use $m=2$ and get
$$
\wt(c_{14}) \;=\; 1 + \wt(c_{21}) + \wt(c_{22}) + \wt(c_{23}) + \wt(c_{24}) 
\;=\; p_{23} + 2  
$$ 

\item[(6c)] Let $\wt(c_{13}) = p_{13}$ and put the inequalities
$p_{13} > \wt(c_{12}) = p_{23}+2$ as well as $p_{13} > \wt(c_{14}) = p_{23}+2$
into the set~$S$.

\item[(8)] Finally, we see that $p_{23}=1$ and $p_{13}=4$ satisfy all inequalities in~$S$.
\end{enumerate}
By examining the neighbours used in these assignments,  if we use the tuple
$Z=(c_{11}, c_{12}, c_{13}, c_{14}, c_{21}, c_{23}, c_{32}, c_{34})$ and the
weight vector $w=(4,3,4,3,1,1,1,2)$ then the following natural generators are
coherently $Z$-separating for any elimination ordering~$\sigma$ for~$Z$ 
which is compatible with the grading given by~$w$:
\begin{align*} 
\ND(1,3)_3 &\;=\;  c_{11} - c_{33} + c_{31}c_{32} +c_{34}c_{41} \\
\ND(2,4)_2 &\;=\; c_{12} - c_{24} + c_{21}c_{22} +c_{23}c_{42} \\
\ND(1,3)_1 &\;=\; - c_{13} + c_{12}c_{31} +c_{14}c_{41}  \\
\AR(3,4)_2 &\;=\; c_{14} + c_{21}c_{24} -c_{22}c_{33} -c_{24}c_{43} +c_{23}c_{44}\\
\ND(1,3)_4 &\;=\; c_{21} - c_{43} + c_{31}c_{42} +c_{41}c_{44}\\
\ND(1,3)_2 &\;=\; - c_{23} + c_{22}c_{31} +c_{24}c_{41} \\
\ND(2,4)_4 &\;=\; c_{32} - c_{44} + c_{22}c_{41} +c_{42}c_{43} \\
\ND(3,4)_4 &\;=\; c_{34} - c_{23} + c_{24}c_{41} -c_{33}c_{42}
\end{align*}

Now let us compare these polynomials to the ones computed via the approach 
described in~\cite{Hui1}, Section~7 (see Remark~\ref{rem-Hui1approach}).
First we have to calculate the three neighbour syzygies of the
polynomials in the generic $\OO$-border prebasis. We get
$$
A =\begin{pmatrix}
x       &  c_{31} &  -1        &  c_{41}     \cr
-c_{24} &  c_{33} &  x -c_{44} &  -y +c_{43} \cr
-c_{22}  &  -y     &  -c_{42}   &  1 
\end{pmatrix}
$$
Up to sign, the maximal minors of~$A$ are
\begin{align*}
y^2 &+ (c_{33}c_{41}c_{42} {-} c_{31}c_{42}c_{43} {+} c_{31}c_{44} {-} c_{33})
   +(c_{31}c_{42} {+} c_{41}c_{44} {-} c_{43})\, y \\
   &-c_{31}\, x  -c_{41}\, xy  \\
x^2 &+ (c_{24}c_{41}c_{42} {-} c_{22}c_{41}c_{44} {+} c_{22}c_{43} {-} c_{24}) -c_{22}\,y
   +(c_{22}c_{41} {+} c_{42}c_{43}\\ 
   &-c_{44})\, x -c_{42}\, xy  \\
xy^2  &+ (-c_{22}c_{33}c_{41} {+} c_{22}c_{31}c_{43} {-} c_{24}c_{31})  
   +(-c_{22}c_{31} {-} c_{24}c_{41})\,y\\  
   &-c_{33}\,x -c_{43}\, xy \\
x^2y  &+ (-c_{24}c_{31}c_{42} {+} c_{22}c_{31}c_{44} {-} c_{22}c_{33}) -c_{24}\, y\\  
   &+(-c_{22}c_{31} {-} c_{33}c_{42})\, x -c_{44}\, xy  
\end{align*}
Comparing the coefficients of these polynomials to those of the generic
$\OO$-border prebasis, we find the following equalities:
$$
\begin{array}{l}
c_{11} =    -(c_{33}c_{41}c_{42} -c_{31}c_{42}c_{43} +c_{31}c_{44}-c_{33} )\cr
c_{12} =  -(c_{24}c_{41}c_{42} -c_{22}c_{41}c_{44} +c_{22}c_{43} -c_{24}) \cr
c_{13} =   -(-c_{22}c_{33}c_{41} +c_{22}c_{31}c_{43} -c_{24}c_{31} ) \cr
c_{14}  =   -(-c_{24}c_{31}c_{42} +c_{22}c_{31}c_{44}-c_{22}c_{33}) \cr
c_{21} =  -(c_{31}c_{42} +c_{41}c_{44} -c_{43}) \cr
c_{23} = -(-c_{22}c_{31} -c_{24}c_{41}) \cr
c_{32} =  -(c_{22}c_{41} +c_{42}c_{43}  -c_{44}) \cr
c_{34} = -(-c_{22}c_{31} -c_{33}c_{42})
\end{array}
$$
The set  $\{ c_{11}+ c_{33}c_{41}c_{42} -c_{31}c_{42}c_{43} 
+c_{31}c_{44} -c_{33}, \dots, c_{34}-c_{22}c_{31} -c_{33}c_{42} \}$
which results from the equalities
is exactly the reduced $\sigma$-Gr\"obner basis of the ideal
generated by the polynomials $\ND(1,3)_3$, $\dots$, $\ND(3,4)_4$ given above.

Notice that in this example the scheme~$\BO$ is an affine cell by~\cite{Hui1}, Prop.~7.3.1.
Hence the elimination ideal is $\langle 0\rangle$ and the above polynomials
are a minimal system of generators of~$I(\BO)$.
\end{example}

A particular family of order ideals in which we can describe the results of
the Weight Assignment Algorithm explicitly is given as follows.

\begin{remark}{\bf (Planar Boxes)}\label{rem-twoboxes}\\
For $\alpha,\beta \in \mathbb{N}_+$, the order ideal $\OO = \mathbb{T}^2 \setminus 
\langle x^\alpha, y^\beta \rangle$ is called the {\bf planar $(\alpha,\beta)$-box}.
With the obvious meaning, there are~$\alpha$ terms in the {\bf upper border} 
and in the {\bf upper rim}, of~$\OO$, and there are~$\beta$ terms in the 
{\bf right border} and the {\bf right rim}. Notice that the term 
$x^{\alpha-1} y^{\beta-1}$ belongs to the upper rim and to the right rim. 

Let us look at the $x$-exposed and $y$-exposed indeterminates $c_{ij}$ in this setting.
\begin{enumerate}
\item[(a)] All terms $b_j$ in the upper border, except for $b_j=x^{\alpha-1} y^\beta$,
have an up-neighbour in $x$-direction, i.e., there is a next-door neighbour pair
$x b_j = b_{j'}$. The term $b_j=x^{\alpha-1} y^\beta$ has an across-the-street neighbour
$b_{j'} =  x^\alpha y^{\beta-1}$, since $x b_j = y b_{j'}$. Altogether, each term~$b_j$ 
in the upper border and each term~$t_i$ in the right rim give rise to an $x$-exposed
indeterminate~$c_{ij}$. In this way we get $\alpha \beta$ $x$-exposed indeterminates.

\item[(b)] Analogously, every term~$b_j$ in the right border, combined with a term~$t_i$
in the upper rim, yields a $y$-exposed indeterminate~$c_{ij}$, and we get $\alpha \beta$
$y$-exposed indeterminates in this fashion.
\end{enumerate}
Since there are no other exposed indeterminates, we have $\# C^{\rm ex} = 2\alpha\beta$.
We know that we can eliminate all non-exposed indeterminates, and the exposed
indeterminates generate~$B_\OO$ (see Corollary~\ref{cor-exposedgenerate}.a). Therefore
$\dim(\BO) = 2 \mu = 2 \alpha \beta$ shows that $\BO \cong \mathbb{A}^{2\alpha\beta}_K$
and that~$B_\OO$ is a polynomial ring over~$K$ generated by the residue classes
of the indeterminates in $C^{\rm ex}$ described in~(a) and~(b).
\end{remark}

Let us see an example which illustrates the remark.

\begin{example}{\bf (The (2,3)-Box)}\label{ex-23Exposed}\\
In $P=\mathbb{Q}[x,y]$, and consider the order ideal 
$\OO = \{t_1,t_2,t_3,t_4, t_5, t_6\}$, where 
$t_1=1$, $t_2=y$, $t_3=x$, $t_4=y^2$, $t_5=xy$, and $t_6=xy^2$. Then the border of~$\OO$
is $\partial\OO = \{b_1, b_2, b_3, b_4, b_5\}$ with $b_1=x^2$, $b_2=y^3$, $b_3 = x^2y$,
$b_4 = xy^3$, and $b_5 = x^2y^2$.  
\begin{figure}[ht]
	\makebox{\beginpicture
		\setcoordinatesystem units <1cm,1cm>
		\setplotarea x from 0 to 3, y from 0 to 3.5
		\axis left /
		\axis bottom /
		\arrow <2mm> [.2,.67] from  2.5 0  to 3 0
		\arrow <2mm> [.2,.67] from  0 3.5  to 0 4
		\put {$\scriptstyle x^i$} [lt] <0.5mm,0.8mm> at 3.1 0
		\put {$\scriptstyle y^j$} [rb] <1.7mm,0.7mm> at 0 4
		\put {$\bullet$} at 0 0
		\put {$\bullet$} at 1 0
		\put {$\bullet$} at 0 1
		\put {$\bullet$} at 1 1
		\put {$\bullet$} at 1 2
		\put {$\bullet$} at 0 2
		\put {$\scriptstyle 1$} [lt] <-1mm,-1mm> at 0 0
		\put {$\scriptstyle t_1$} [rb] <-1.3mm,0.4mm> at 0 0
		\put {$\scriptstyle t_3$} [rb] <-1.3mm,0.4mm> at 1 0
		\put {$\scriptstyle t_2$} [rb] <-1.3mm,0mm> at 0 1
		\put {$\scriptstyle t_4$} [rb] <-1.3mm,0mm> at 0 2
		\put {$\scriptstyle t_5$} [rb] <-1.3mm,0mm> at 1 1
		\put {$\scriptstyle t_6$} [lb] <-4mm,0.4mm> at 1 2
		\put {$\scriptstyle b_1$} [lb] <0.8mm,0.8mm> at 2 0
		\put {$\scriptstyle b_2$} [lb] <0.8mm,0.8mm> at  0 3
		\put {$\scriptstyle b_3$} [lb] <0.8mm,0.8mm> at 2 1
	        \put {$\scriptstyle b_4$} [lb] <0.8mm,0.8mm> at 1 3
	        \put {$\scriptstyle b_5$} [lb] <0.8mm,0.8mm> at 2 2
		\put {$\times$} at 0 0
	
		\put {$\circ$} at 2 0
		\put {$\circ$} at 2 1
		\put {$\circ$} at 2 2
		\put {$\circ$} at 1 3
		\put {$\circ$} at 0 3
		\endpicture}
\end{figure}
The elements of the upper border are~$b_2$ and~$b_4$, and the elements of the upper rim
are~$t_4$ and~$t_6$. The elements of the right border are $b_1$, $b_3$, and~$b_5$. 
The elements of the right rim are~$t_3$, $t_5$, and~$t_6$.

The next-door neighbour pair $x b_2 = b_4$ yields the $x$-exposed indeterminates 
$c_{32}$, $c_{52}$, and $c_{62}$. The across-the-street neighbour pair $x b_4 = y b_5$ 
yields the $x$-exposed indeterminates $c_{34}$, $c_{54}$, and~$c_{64}$.
Hence the set of $x$-exposed indeterminates is $\{ c_{32}, c_{34}, c_{52}, c_{54},
c_{62}, c_{64} \}$. Analogously, the set of $y$-exposed indeterminates is 
$\{ c_{41}, c_{43}, c_{45}, c_{61}, c_{63}, c_{65}\}$. Altogether, we have found 12 exposed
indeterminates and we conclude that there is an isomorphism 
$$
B_\OO \cong \QQ[c_{32}, c_{34}, c_{41}, c_{43}, c_{45}, c_{52}, c_{54}, 
c_{61}, c_{62}, c_{63}, c_{64}, c_{65}]
$$
The weights provided by Algorithm~\ref{alg-WeightAssign} are
$$
\begin{array}{ccccccccccccccc}
c_{11}&  c_{12}&  c_{13}&  c_{14}&  c_{15}&  c_{21}&  c_{22}&  c_{23}&  c_{24}&  c_{25}&  
   c_{31}&  c_{32}&  c_{33}&  c_{34}&  c_{35}  \cr
13& 15& 13& 20& 19& 3& 5& 3& 4& 3& 2& 0& 3& 0& 9 \cr
c_{41}&  c_{42}&  c_{43}&  c_{44}&  c_{45}&  c_{51}&  c_{52}&  c_{53}&  c_{54}&  c_{55}&  
   c_{61}&  c_{62}&  c_{63}&  c_{64}&  c_{65} \cr
  0&1& 0& 1& 0& 1& 0& 1& 0& 2& 0& 0& 0& 0& 0
\end{array}
$$
and the algorithm computes the 18 natural generators whose leading terms are the indeterminates 
with positive weight, i.e., the non-exposed indeterminates.
\end{example}

Let us point out a special merit of Algorithm~\ref{alg-WeightAssign} connected to 
Remark~\ref{rem:differentF}.

\begin{remark}\label{rem:goodweights}
Not all orderings which eliminate non-exposed indeterminates find suitable natural generators
having such indeterminates as leading terms.
For instance, in the above example the polynomial $f\!=\! -c_{22}c_{41} -c_{24}c_{61} -c_{11}+c_{23}$
is a natural generator of the ideal~$I(\BO)$. The weights of the indeterminates 
appearing in~$\Supp(f)$ are 
$$
\begin{matrix}
c_{11}&    c_{22}&  c_{23}&  c_{24}& c_{41}&  c_{61}\cr
   13&  5& 3& 4& 0 & 0
\end{matrix}
$$ 
Hence we have $\LT_\sigma(f) = c_{11}$. Instead, if we use an elimination ordering~$\tau$ 
represented by a matrix whose first row assigns~1 to the non-exposed indeterminates
and 0 to the exposed indeterminates, we get $\LT_\tau(f) = c_{22}c_{41}$. 
Therefore~$\tau$ produces the same reduced Gr\"obner basis as~$\sigma$,
but does not find 18 natural generators having the correct leading terms.
\end{remark}

Next we make good use of exposed indeterminates in the process of 
searching for an optimal embedding.  This goal is achieved by the following variant 
of Algorithm~\ref{alg-compEmb}.

\begin{proposition}{\bf (Optimal Re-embeddings Using $C^{\rm ex}$)}\label{prop-E0plusNonexp}\\
Assume that we are in the setting of Algorithm~\ref{alg-compEmb} with $n=2$. In other words, let
$\OO = \{ t_1,\dots, t_\mu\}$ be an order ideal in~$\mathbb{T}^2$ with border
$\partial \OO = \{ b_1,\dots,b_\nu\}$, and let~$C$ be the set of indeterminates which are the
coefficients of the generic $\OO$-preborder basis. Replace Step~(4) by the following steps.
\begin{enumerate}
\item[(4.1)] Compute the set $C^{\rm ex}$  
of exposed indeterminates in~$C$.

\item[(4.2)] For $i=1,\dots,q$, compute $\widetilde{E}_i =  E_i \cap C^{\rm ex}$.

\item[(4.3)] From each set $E_0\cup (C\setminus C^{\rm ex}) \cup \widetilde{E}_1^\ast \cup \cdots 
\cup \widetilde{E}_q^\ast$, where $\widetilde{E}_i^\ast$ is obtained from~$\widetilde{E}_i$ 
by deleting one element, form a tuple~$Z$ and perform the following steps. 
\end{enumerate}
Then the result is an algorithm which computes the set~$S$ of all tuples~$Z$ 
of distinct indeterminates in~$C$ such that there exists an optimal
$Z$-separating re-embedding of~$I(\BO)$ which eliminates the non-exposed indeterminates.
\end{proposition}

\begin{proof}
From Algorithm~\ref{alg-compEmb} we know that the indeterminates 
in~$E_0$ must be eliminated. From Theorem~\ref{thm-elimNex} we know that non-exposed 
indeterminates can be eliminated. Therefore, if we want to eliminate them, 
we do not need to use the sets~$E_i$, as we do in Algorithm~\ref{alg-compEmb}.
Instead, it suffices to use the sets~$\widetilde{E}_i$.
\end{proof}

Let us use Example~\ref{ex-Section5} to show this proposition at work.

\begin{example}\label{ex-Section5-continued}
In Example~\ref{ex-Section5} we have
\begin{align*}
E_0 = & \{ c_{11},  c_{12},  c_{13},  c_{14},  c_{15},  c_{21},  c_{22},  c_{23},  
c_{24},  c_{25},  c_{31},  c_{32},  c_{33},  c_{34},  c_{35}, \\ 
& \; c_{42},  c_{44},  c_{45},  c_{55},  c_{65} \}\\
E_1 = & \{c_{51},  c_{85}\}, \qquad E_2 = \{c_{43},  c_{54}\}, \qquad
E_3 = \{c_{41},  c_{52},  c_{75}\}
\end{align*}
The exposed indeterminates are
$C^{\rm ex} = \{ c_{51},  c_{53},  c_{54},  c_{61},  c_{62},  c_{63},  c_{64},  
c_{65},_{71},  c_{72},  \allowbreak   c_{73},  c_{74},  c_{75},  
c_{81},  c_{82},  c_{83},  c_{84},  c_{85} \}$.

Notice that $c_{52}$ and~$c_{55}$ are rim indeterminates which are not exposed.
Moreover, as~$c_{65}$ is exposed and trivial, it will be eliminated. 
At this point we already know that $C^{\rm ex} $ does not generate $B_\OO$ minimally.
Furthermore, the indeterminates $c_{41}$, $c_{43}$, and $c_{52}$ are non-exposed. 
Therefore we have $\widetilde{E}_1 = E_1$,  $\widetilde{E}_2 =\{c_{54} \}$, 
and $\widetilde{E}_3 = \{c_{75} \}$, so that $\widetilde{E}_2^\ast =\emptyset$ 
and $\widetilde{E}_3^\ast =\emptyset$.
Consequently, the remaining tuples to be examined are only the two tuples~$Z_1$ and~$Z_2$
formed from $E_0 \cup (C \setminus C^{\rm ex}) \cup \{ c_{51}\}$ and 
$E_0 \cup (C\setminus C^{\rm ex}) \cup \{ c_{85}\}$, respectively.

It turns out that both provide optimal re-embeddings of $I(\BO)$.
In the first case, the exposed indeterminates $c_{65}$ and $c_{51}$ are eliminated, 
and in the second case, the exposed indeterminates $c_{65}$ and $c_{85}$ are eliminated. 

As already observed in Remark~\ref{rem-optZsepPlanar},
this is an example of an order ideal~$\OO$ which does not have the ``sawtooth'' form.
Nevertheless, we get an optimal re-embedding of $I(\BO)$, and $\BO$ turns out to be an affine cell.
Notice also that the example mentioned in~\cite{Hui1}, Remark 7.5.3, is similar to this example.  
\end{example}

\bigskip\bigbreak
%
%

\section{Simplicial Order Ideals}
\label{Simplicial Order Ideals}

In this section we continue to examine re-embeddings of border basis schemes
in special cases. We use the notation introduced previously.
Let~$K$ be a field, let $P=K[x_1,\dots,x_n]$, and let
$\M= \langle x_1,\dots,x_n\rangle$ be the homogeneous maximal ideal of~$P$.
The order ideals we are looking at in this section are defined as follows.

\begin{definition}
Let $d\ge 1$. The order ideal 
$$
\OO =\mathbb{T}^n \setminus \M^{d+1} =
\{ t\in\mathbb{T}^n \mid \deg(t) \le d\}
$$ 
is called the {\bf simplicial order ideal} of type~$d$ in~$\mathbb{T}^n$.
\end{definition}

Let us collect some basic numerical information about simplicial order ideals.

\begin{remark}\label{rem-formulas}
Let $\OO = \{t_1,\dots,t_\mu\}$ be the simplicial order ideal of type~$d$
in~$\mathbb{T}^n$, and let $\partial\OO = \{ b_1,\dots,b_\nu\}$ be its border.
Then the following formulas hold:
\begin{enumerate}
\item[(a)] $\mu = \binom{d+n}{n}$

\item[(b)] $\nu = \binom{d+n}{n-1}$

\item[(c)] $\# \OO^\circ = \binom{d+n-1}{n}$

\item[(d)] $\# \nabla\OO = \binom{d+n-1}{n-1}$

\item[(e)] $\# C = \binom{d+n}{n} \cdot \binom{d+n}{n-1}$

\item[(f)] $\# C^\circ = \binom{d+n-1}{n} \cdot \binom{d+n}{n-1}$

\item[(g)] $\# C^\nabla = \binom{d+n-1}{n-1} \cdot \binom{d+n}{n-1}$
\end{enumerate}
These formulas follow immediately from the 
affine Hilbert function of a polynomial ring (cf.~\cite{KR2}, Propositions 5.1.13
and 5.6.3.c).
\end{remark}

The following proposition allows us to find optimal re-embeddings of border basis
schemes in the case of simplicial order ideals.

\begin{proposition}\label{prop-simplicialZ}
Let $d\ge 1$, and let $\OO= \{t_1,\dots,t_\mu\}$ be the simplicial order ideal of type~$d$
in~$\mathbb{T}^n$, and let  $\partial\OO = \{b_1,\dots, b_\nu\}$ be its border.
Furthermore, let $s = \#C^\circ = \#\OO^\circ \cdot \#(\partial\OO)=
\binom{d+n-1}{n} \cdot \binom{d+n}{n-1}$.
Then there exists a term ordering~$\sigma$ such that for each indeterminate $c_{ij} \in C^\circ$ there
is a natural generator~$f$ of~$I(\BO)$ with $c_{ij}=\LT_\sigma(f)$.
Consequently, $\sigma$ is a $Z$-separating term ordering for the ideal ~$I(\BO)$.
\end{proposition}

\begin{proof}
Let~$Z$ be the tuple of indeterminates in~$C^\circ$, ordered arbitrarily. 
We define an elimination ordering~$\sigma$ for~$Z$ by constructing a suitable matrix~$M$.
In the first row of~$M$ and the column corresponding to an indeterminate~$c_{ij}$, 
we put the number $2\delta -1$, where~$\delta$ is the
total arrow degree $\deg(c_{ij}) = \deg(b_j) - \deg(t_i) =  d+1 -\deg(t_i)$.

Next we define the second row of the matrix~$M$ as follows. For each arrow degree~$\delta$,
let $S_\delta$ be the set of all indeterminates $c_{ij}\in C^\circ$ with
$\deg_W(c_{ij})=\delta$. Let $D$ be the set of all arrow degrees~$\delta$ such that
$S_\delta\ne \emptyset$. For each $\delta = (\delta_1,\dots,\delta_n)$ in~$D$, 
we choose a fixed,  for instance the first, number $k\in \{1,\dots,n\}$ such that $\delta_k > 0$, 
and note  that this is possible by Proposition~\ref{prop-poscomp}. 
For $c_{ij} \in S_\delta$, we write
$t_i = x_1^{\alpha_1} \cdots x_n^{\alpha_n}$ with $\alpha_1,\dots,\alpha_n\in\mathbb{N}$.
Then we put the number $\alpha_k$ into the column corresponding to~$c_{ij}$
of the second row of the matrix~$M$. In the columns corresponding to $c_{ij}\in C^\nabla$
we put zeros into the second row of~$M$. Finally, we complete the matrix~$M$ to
a matrix defining an elimination ordering $\sigma = {\rm ord}(M)$ for~$Z$.

Our goal is to show that~$\sigma$ is a $Z$-separating term ordering for~$I(\BO)$.
For this purpose we prove that for each~$c_{ij} \in C^\circ$ there exists a
generator $f_{ij}\in \AR_{\OO}$ such that the following conditions are satisfied:
\begin{enumerate} 
\item[(1)] The linear part of~$f_{ij}$ is (up to sign) $c_{ij}$ or of the form
$c_{ij} - c_{i'j'}$ with $c_{ij} >_\sigma c_{i'j'}$.

\item[(2)] For every term $c_{\alpha \beta} c_{\gamma\delta}$ of degree two in 
$\Supp(f_{ij})$, we have $c_{ij} >_\sigma c_{\alpha\beta} c_{\gamma\delta}$.
\end{enumerate}
Then it follows that $\LT_\sigma(f_{ij}) = c_{ij}$ and hence~$\sigma$ is a
$Z$-separating term ordering for the ideal $I(\BO)$.

To prove the existence of~$f_{ij}$ with the desired properties, 
we let $c_{ij} \in C^\circ$, and we let $\delta = (\delta_1,\dots,\delta_n)$
be the arrow degree of~$c_{ij}$. Let $k\in \{1,\dots,n\}$ be the component of~$\delta$
with $\delta_k >0$ chosen in the above construction, and let $t_i = x_1^{\alpha_1}
\cdots x_n^{\alpha_n}$. Now we distinguish two cases:

{\it Case 1:} $\alpha_k=0$. By assumption, the exponent of~$b_j$ with respect to~$x_k$
is positive. Let $\ell\in \{1,\dots,n\} {\setminus} \{k\}$. Then $x_\ell \, b_j / x_k$ 
is an element of $\partial\OO$, and hence equal to $b_{j'}$ for some $j'\in \{1,\dots,\nu\}$.
Notice that $(b_j,b_{j'})$ is an across-the-rim neighbour pair.
Moreover, the term $t_m = x_\ell t_i$ is contained in~$\OO$, since $t_i\in \OO^\circ$.
However $t_m$ is not divisible by~$x_k$. Next, consider the polynomial
$f_{ij} = \AR(j,j')_m$. The linear part of~$f_{ij}$ is (up to sign) $c_{ij}$. Furthermore, the
terms in the quadratic part are of the form $c_{m\kappa} c_{\rho j''}$ where $j''\in \{j,j'\}$
and $c_{\rho j''} \in C^\nabla$. Here $\deg(t_m)= \deg(t_i)+1$ implies that the entry 
in the first row of~$M$ corresponding to $c_{m\kappa}$ is two less than the entry 
corresponding to~$c_{ij''}$.
Since the entry corresponding to $c_{\rho j''}$ is one, it follows that we have
$c_{ij} >_\sigma c_{m\kappa} c_{\rho j''}$. Thus we get $\LT_\sigma(f_{ij}) = c_{ij}$.

{\it Case 2:} $\alpha_k >0$. In this case, we choose an index $\ell \in \{1,\dots,n\}
{\setminus} \{k\}$ and note that $x_\ell b_j / x_k$ is in~$\partial\OO$. Let $j'\in\{1,\dots,\nu\}$
such that $x_k b_{j'} = x_\ell b_j$. Since $t_i \in \OO^\circ$, we know that $t_m = x_\ell t_i$
is in~$\OO$. Since $\alpha_k>0$, it follows that $x_k$ divides $t_m$ and there exists an
index $i'\in \{1,\dots,\mu\}$ with $x_k t_{i'} = t_m = x_\ell t_i$.
Then $(b_j, b_{j'})$ is an across-the-rim neighbour pair, and we consider the polynomial
$f_{ij} = \AR(j,j')_m$. The linear part of~$f_{ij}$ is (up to sign) given by
$c_{ij} -c_{i' j'}$. Here the facts that the total arrow degree of $c_{ij}$ and $c_{i' j'}$
are the same and that the $x_k$-exponent of~$t_{i'}$ is $\alpha_k -1$ imply
$c_{ij} >_\sigma c_{i'j'}$. 

The terms of degree two in the support of~$f_{ij}$ are seen to be smaller than~$c_{ij}$
with respect to~$\sigma$ in exactly the same way as in Case~1. Altogether, it follows
that $\LT_\sigma(f_{ij}) = c_{ij}$, and the proof of the proposition is complete.
\end{proof}

As mentioned above, this proposition has a number of useful consequences.

\begin{corollary}
Let $d\ge 1$, and let $\OO= \{t_1,\dots,t_\mu\}$ be the simplicial order ideal of type~$d$
in~$\mathbb{T}^n$ with border $\partial\OO = \{b_1,\dots, b_\nu\}$.
\begin{enumerate}
\item[(a)] The linear part of~$I(\BO)$ is given by
$\Lin_\M(I(\BO)) = \langle C^\circ \rangle_K$.

\item[(b)] We have  $\dim_K (\Cot_\m(B_\OO))= \# C^\nabla = \binom{d+n-1}{n-1} 
\cdot \binom{d+n}{n-1}$.

\item[(c)] We have $\edim(B_\OO) = \dim_K (\Cot_{\m}(B_\OO)) $.  
\end{enumerate}
In particular, for the term ordering~$\sigma$ constructed in the proposition, 
the $Z$-sep\-a\-ra\-ting re-embedding  
$\Phi: K[C] / I(\BO) \longrightarrow K[C^\nabla] / ( I\cap K[C^\nabla])$
defined by the reduced $\sigma$-Gr\"obner basis of~$I(\BO)$ is optimal.
\end{corollary}

\begin{proof}
To prove~(a), we note that the indeterminates of~$C^\circ$ are in
$\Lin_\M(I(\BO))$ by the proposition. By~\cite{KSL}, Cor.~2.8.b, there is no
linear part of one of the natural generators of~$I(\BO)$ which contains
an indeterminate of~$C^\nabla$ in its support. Hence we get
$\Lin_\M(I(\BO)) = \langle C^\circ\rangle_K$.

The first equality in~(b) follows directly from the proposition and~(a).
Since we have the equality $\#C = \#\OO\cdot \#(\partial\OO) =
\binom{d+n}{n}\cdot \binom{d+n}{n-1}$, we get
$$
\# C^\nabla = \#C - \#C^\circ =   \tbinom{d+n}{n}\cdot \tbinom{d+n}{n-1} 
- \tbinom{d+n-1}{n} \cdot \tbinom{d+n}{n-1}
= \tbinom{d+n-1}{n-1}\cdot \tbinom{d+n}{n-1}
$$
Finally, we  prove~(c). As an application of~\cite{KLR2}, Cor.~4.2, we know that
the embedding dimension of~$B_\OO$ is $\#C - \#C^\circ =\# C^\nabla$, 
and this number is equal to $\dim_K (\Cot_\m(B_\OO))$ by~(b).
\end{proof}

This proposition allows us to check when a simplicial order ideal 
has a border basis scheme which is an affine cell.

\begin{corollary}\label{cor-affand singular}
In the setting of the proposition, let $p=(0,\dots,0)$ be the monomial point of~$\BO$.
\begin{enumerate}
\item [(a)] If $n=2$ then $\BO$ is isomorphic to an affine space.

\item[(b)] If $n\ge 3$ then~$p$ is a singular point of~$\BO$.
\end{enumerate}
\end{corollary}

\begin{proof}
For $n = 2$, it follows from the proposition that
$\BO$ can be embedded into an affine space of dimension
$\dim_K (\Cot_\m(B_\OO)) = \binom{d+1}{1}\cdot \binom{d+2}{1} =   (d+1)(d+2)$.
It~is also known that the  scheme $\BO$ is irreducible and its dimension is given by
the formula
$\dim(\BO)=2\cdot  \#(\OO) = 2\frac{(d+1)(d+2)}{2} = (d+1)(d+2)$.
The conclusion follows.

To prove~(b), we observe that the dimension of the principal component of~$\BO$
is $n\cdot \#\OO = n\cdot \binom{d+n}{n}$.
Moreover, it follows from part~(c) of the proposition that
$\dim_K (\Cot_\m(B_\OO)) = \binom{d+n-1}{n-1} \cdot \binom{d+n}{n-1}$.
Hence it suffices to show that for $n\ge 3$ we have
$$
n\cdot \tbinom{d+n}{n} < \tbinom{d+n-1}{n-1} \cdot \tbinom{d+n}{n-1}  \eqno{(1)}
$$
To prove (1), we note that $n\cdot \binom{d+n}{n} {=}(d+n)\cdot \binom{d+n-1}{n-1}$.
Hence~(1) is equivalent to
$$
d+n <  \tbinom{d+n}{n-1}   \eqno{(2)}
$$
Therefore we need to show that for $n\ge 3$ we have
$$
1 <  \tfrac{ (d+n-1)\cdot(d+n-2)\cdots(d+2)}{(n-1)!}=
\tfrac{d+n-1}{n-1}\cdot \tfrac{d+n-2}{n-2}\cdots \tfrac{d+2}{2}
$$
which is obvious, since $d\ge 1$.

Notice that for $n=2$ we have $d+n = \binom{d+n}{n-1}$. Thus we have equality in~(2),
and hence also equality in~(1). This shows $\dim_K(\Cot_\m(B_\OO)) = \dim(\BO)$ 
and implies that~$p$ is a smooth point of $\BO$, in agreement with~(a).
\end{proof}

\begin{remark}\label{rem-rim=exposed}
It is easy to see that for planar simplicial order ideals the sets of exposed and rim 
indeterminates coincide. Therefore in this case it is also possible to use 
Algorithm~\ref{alg-WeightAssign} to construct an optimal re-embedding of~$\BO$
(see also our next Example~\ref{ex-123}).
\end{remark}

The following examples illustrate the main results of this section.

\begin{example}\label{ex-123}
In $P=\mathbb{Q}[x,y]$, and consider the simplicial order ideal 
$\OO = \{t_1, t_2, t_3, t_4, t_5, t_6\}$, 
where $t_1=1$, $t_2=y$, $t_3=x$, $t_4=y^2$, $t_5=xy$, and $t_6=x^2$. 
Then the border of~$\OO$ is $\partial\OO = \{b_1, b_2, b_3, b_4\}$ with 
$b_1=y^3$, $b_2=xy^2$, $b_3 = x^2y$, and $b_4 = x^3$.

Thus we have $d=2$, $n=2$, $\mu = \binom{d+n}{n} = 6$, $\nu = \binom{d+n}{n-1} = 4$,
$\# \OO^\circ {=} \binom{d+n-1}{n} {=} 3$,
$\# \nabla\OO = \binom{d+n-1}{n-1} = 3$,
$\# C = \binom{d+n}{n} \cdot \binom{d+n}{n-1} = 24$,
$\# C^\circ = \binom{d+n-1}{n} \cdot \binom{d+n}{n-1} = 12$, and
$\# C^\nabla = \binom{d+n-1}{n-1} \cdot \binom{d+n}{n-1} = 12$
(see Remark~\ref{rem-formulas}).

By following the proof of Proposition~\ref{prop-simplicialZ} 
we construct a matrix whose first two rows are
$$
\begin{array}{ll}
5,\  5,\  5,\  5,\     3,\  3,\  3,\  3,\     3,\  3,\  3,\  3,\     1,\  1,\  1,\  1,\     
   1,\  1,\  1,\  1,\     1,\  1,\  1,\  1 \cr
0,\  0,\  0,\  0,\     1,\  0,\  0,\  0,\     0,\  0,\  1,\  1,\     0,\  0,\  0,\  0,\     
   0,\  0,\  0,\  0,\     0,\  0,\  0,\  0
\end{array}
$$
and then complete it to a matrix~$M$ which defines an elimination ordering for
$Z = (c_{11},c_{12}, c_{13}, c_{14},  c_{21},c_{22}, c_{23}, c_{24}, 
c_{31},c_{32}, c_{33}, c_{34})$. 

Using the term ordering $\sigma = \ord(M)$, we check that the natural generators of~$I(\BO)$
given by
$$
\begin{array}{l}
  c_{32}c_{41} -c_{31}c_{42} +c_{33}c_{51} -c_{32}c_{52} +c_{34}c_{61} -c_{33}c_{62} +c_{11}, \cr
  c_{22}c_{41} -c_{21}c_{42} +c_{23}c_{51} -c_{22}c_{52} +c_{24}c_{61} -c_{23}c_{62} -c_{12}, \cr
  c_{22}c_{42} -c_{21}c_{43} +c_{23}c_{52} -c_{22}c_{53} +c_{24}c_{62} -c_{23}c_{63} -c_{13}, \cr
  c_{22}c_{43} -c_{21}c_{44} +c_{23}c_{53} -c_{22}c_{54} +c_{24}c_{63} -c_{23}c_{64} -c_{14}, \cr
  -c_{42}c_{51} +c_{41}c_{52} -c_{52}^2 +c_{51}c_{53} +c_{54}c_{61} -c_{53}c_{62} +c_{21} -c_{32}, \cr
  c_{43}c_{51} -c_{42}c_{52} +c_{44}c_{61} -c_{43}c_{62} -c_{22}, \cr
  c_{42}^2 -c_{41}c_{43} +c_{43}c_{52} -c_{42}c_{53} +c_{44}c_{62} -c_{43}c_{63} -c_{23}, \cr
  c_{42}c_{43} -c_{41}c_{44} +c_{43}c_{53} -c_{42}c_{54} +c_{44}c_{63} -c_{43}c_{64} -c_{24}, \cr
  -c_{42}c_{61} +c_{41}c_{62} -c_{52}c_{62} +c_{51}c_{63} -c_{62}c_{63} +c_{61}c_{64} +c_{31}, \cr
  -c_{43}c_{61} +c_{42}c_{62} -c_{53}c_{62} +c_{52}c_{63} -c_{63}^2 +c_{62}c_{64} +c_{32}, \cr
  -c_{44}c_{61} +c_{43}c_{62} -c_{54}c_{62} +c_{53}c_{63} +c_{33}, \cr
  -c_{44}c_{51} +c_{43}c_{52} +c_{53}^2 -c_{52}c_{54} +c_{54}c_{63} -c_{53}c_{64} +c_{23} -c_{34}
\end{array}
$$
define a tuple~$F$ whose tuple of $\sigma$-leading terms is~$Z$. Moreover, the entries of~$F$ form a 
minimal set of generators of~$I(\BO)$, and the reduced $\sigma$-Gr\"obner basis of 
$I(\BO) =\langle F \rangle$ defines an isomorphism between $B_\OO$ and  $\QQ[C {\setminus}Z]$.

As mentioned in Remark~\ref{rem-rim=exposed}, we can also use Algorithm~\ref{alg-WeightAssign} 
to construct a suitable term ordering. It is given by a matrix~$N$ whose first row is
$$
\begin{array}{ll}
7, \ 6, \ 6, \ 6, \   2, \ 1, \ 1, \ 1,  \  1, \ 1, \ 2, \ 2, \  0, \ 0, \ 0, \ 0,\  
0, \ 0, \ 0,\  0, \  0, \ 0, \ 0,\ 0
\end{array}
$$
and then~$N$ is completed to a matrix which defines an elimination ordering~$\sigma'$ for 
$Z = (c_{11},c_{12}, c_{13}, c_{14},  c_{21},c_{22}, c_{23}, c_{24}, 
c_{31},c_{32}, c_{33}, c_{34})$. 
Using the term ordering $\sigma' =\ord(N)$, we get exactly the same result which we got with~$\sigma$.

As in Remark~\ref{rem:goodweights}, if we use instead an elimination ordering~$\tau$ represented 
by a matrix whose first row assigns~1 to the indeterminates in~$Z$ and~0 to 
the remaining indeterminates, then~$\tau$ yields the same reduced Gr\"obner basis as~$\sigma$,
but does not find the 12 natural generators with the correct leading terms.
\end{example}

In the final example of this section we construct an optimal embedding of a border basis scheme
which is not an isomorphism with an affine space.

\begin{example}\label{ex-13}
In $P=\mathbb{Q}[x,y,z]$, and consider the simplicial order ideal $\OO = \{t_1,t_2,t_3,t_4\}$, 
where $t_1=1$, $t_2=z$, $t_3=y$, and $t_4=z$. Then the border of~$\OO$
is $\partial\OO = \{b_1, b_2, b_3, b_4,b_5,b_6\}$ with $b_1=z^2$, $b_2=yz$, $b_3 = xz$, and
$b_4 = y^2$, $b_5 = xy$, and $b_6 = x^2$.

Here we have $d=1$, $n=3$, $\mu = 4$, $\nu = 6$, and thus
$\# \OO^\circ = \binom{d+n-1}{n} = 1$,
$\# \nabla\OO = \binom{d+n-1}{n-1} = 3$,
$\# C = \binom{d+n}{n} \cdot \binom{d+n}{n-1} = 24$,
$\# C^\circ = \binom{d+n-1}{n} \cdot \binom{d+n}{n-1} = 6$, and
$\# C^\nabla = \binom{d+n-1}{n-1} \cdot \binom{d+n}{n-1} = 18$.

As in the proof of Proposition~\ref{prop-simplicialZ},
we construct a matrix whose first two rows are
$$
\begin{array}{ll}
3,\  3,\  3,\  3,\  3,\  3,  \     1,\  1,\  1,\  1,\ 1,\ 1, \     1,\  1,\  1,\  1,\ 1,\ 1, \    
   1,\  1,\  1,\  1,\ 1, \ 1     \cr
0,\  0,\  0,\  0,\  0,\  0, \      0,\  0,\  0,\  0,\  0,\  0, \     0,\  0,\  0,\  0,\  0,\  0, \   
   0,\  0,\  0,\  0,\  0,\  0
\end{array}
$$
and then, after deleting the second useless row, we complete it to a matrix~$M$ which defines 
an elimination ordering $\sigma = \ord(M)$ for
$Z = (c_{11},c_{12}, c_{13}, c_{14}, c_{15}, c_{16})$.
Then the polynomials
$$
\begin{array}{l}
  -c_{22}c_{31} +c_{21}c_{32} -c_{32}^2 +c_{31}c_{34} +c_{35}c_{41} -c_{33}c_{42} +c_{11}, \cr
  c_{24}c_{31} -c_{22}c_{32} +c_{25}c_{41} -c_{23}c_{42} -c_{12}, \cr
  c_{25}c_{31} -c_{22}c_{33} +c_{26}c_{41} -c_{23}c_{43} -c_{13}, \cr
  c_{22}^2 -c_{21}c_{24} +c_{24}c_{32} -c_{22}c_{34} +c_{25}c_{42} -c_{23}c_{44} -c_{14}, \cr
  c_{22}c_{23} -c_{21}c_{25} +c_{25}c_{32} -c_{22}c_{35} +c_{26}c_{42} -c_{23}c_{45} -c_{15}, \cr
  c_{23}^2 -c_{21}c_{26} +c_{25}c_{33} -c_{22}c_{36} +c_{26}c_{43} -c_{23}c_{46} -c_{16}, \cr
    \end{array}
$$
yield a tuple~$F$ such that~$Z$ is exactly the tuple of $\sigma$-leading terms of~$F$. 
The reduced $\sigma$-Gr\"obner basis of $I(\BO)=\langle F \rangle$ 
defines an optimal re-embedding of~$B_\OO$. 

More specifically, it defines an isomorphism between $B_\OO$ and a ring of the form
$\QQ[C{\setminus }Z]/ J$, where~$J$ is a homogenous ideal which is minimally generated 
by 15 quadratic polynomials (see Example~4.6 of~\cite{KLR2}).
In accordance with Corollary~\ref{cor-affand singular}, we see that the point~$p$ is singular.
\end{example}

\bigskip\bigbreak
%
%

\section{Further Examples}
\label{Further Examples}

In this section we apply the methods developed earlier to various examples
in order to show how they can be applied to find good or even optimal re-embeddings
of border basis schemes and to discover some of their properties.

\begin{example}{\bf (The L-Shape Reconsidered)}\label{ex-theLshapeCont}\\
As in Example~\ref{ex-theLshape}, we use the order ideal 
$\OO = \{1,\, y,\, x,\, y^2,\, x^2\}$ in  $P=\mathbb{Q}[x,y]$ which was 
called the {\bf L-shape}. Its border is 
$\partial\OO = \{ xy,\, y^3,\, xy^2,\, x^2y, x^3\}$.
We have $\mu = \nu =5$, and $\QQ[C] = \QQ[c_{11}, \dots, c_{55}]$ 
is a polynomial ring in 25 indeterminates.
The cotangent equivalence classes are
\begin{align*}
E_0 = & \{c_{11},\ c_{12},\; c_{13},\; c_{14},\; c_{15},\; c_{23},\; c_{24},\;
          c_{25},\; c_{32},\; c_{33},\; c_{34},\; c_{44},\; c_{53} \}\\
E_1 = & \{c_{31},\ c_{54}\}, \qquad E_2 = \{ c_{21},\ c_{43}\}
\end{align*}
and $C\setminus C^{\rm ex} = \{ c_{11},\ c_{12},\ c_{13},\ c_{14},\ c_{15},\ c_{23},
c_{25},\ c_{32},\ c_{34}\}$ is the set of non-exposed indeterminates. 
It is contained in $E_0$.

According to Definition~\ref{def-genGFan} and Theorem~\ref{thm-shapeofSsigma}.c,
we have 
\begin{align*}
\LTGFan(I(\BO)_{\rm lin})  \;=\; \{
&E_0\cup \{c_{21},c_{31} \},\; E_0\cup \{c_{21}, c_{54} \},\\ 
&E_0\cup \{c_{31}, c_{43} \},\;  E_0\cup \{c_{43},c_{54} \} \}
\end{align*}
Using the methods explained in Subsection~1.B, in particular the characterization 
given in Proposition~\ref{prop-CheckZ}, we find that no subset of cardinality 
at least~14 of the elements of $\LTGFan(I(\BO)_{\rm lin})$ yields a separating re-embedding 
of~$I(\BO)$. However, many subsets of cardinality~13 do, although~$E_0$ is not among them. 
One of the best re-embeddings of~$I(\BO)$ is given by 
$$
Z \;=\; (c_{11},\; c_{12},\; c_{13},\; c_{14},\; c_{15},\; c_{21},\; c_{23},\; 
     c_{24},\; c_{25},\; c_{32},\; c_{34},\; c_{44},\; c_{53})
$$
It gives rise to a $\QQ$-algebra isomorphism 
$$
B_\OO \cong \QQ[c_{22}, c_{31}, c_{33}, c_{35},c_{41}, c_{42}, c_{43}, c_{45}, 
c_{51}, c_{52}, c_{54}, c_{55}]/\langle f_1,\, f_2\rangle$$
where
\begin{align*}
f_1 \;=\; & c_{31}c_{41}^2c_{51}^2 +c_{41}^3c_{51}c_{52} +c_{41}^2c_{51}^2c_{54} -c_{41}c_{43}c_{51}^2 
       -c_{45}c_{51}^3 +c_{41}c_{51}^2c_{55} \\
     & +c_{31}c_{41}c_{51} -c_{41}c_{42}c_{51} +c_{41}^2c_{52} +c_{43}c_{51} +c_{51}c_{55} +c_{31} -c_{54},\\
f_2 \;=\; & 2c_{31}c_{41}^3c_{51}c_{52} +2c_{41}^4c_{52}^2 +2c_{41}^3c_{51}c_{52}c_{54} 
            +c_{31}c_{41}c_{43}c_{51}^2 -c_{31}c_{45}c_{51}^3\\ 
     & -c_{41}^2c_{43}c_{51}c_{52} -3c_{41}c_{45}c_{51}^2c_{52} +c_{41}c_{43}c_{51}^2c_{54} 
       -c_{45}c_{51}^3c_{54} \\
     & +c_{31}c_{41}c_{51}^2c_{55} +3c_{41}^2c_{51}c_{52}c_{55} +c_{41}c_{51}^2c_{54}c_{55} 
        -c_{22}c_{41}c_{51}\\
     & +2c_{33}c_{41}c_{51} -c_{41}c_{42}^2c_{51} +c_{35}c_{51}^2 -c_{43}^2c_{51}^2 +c_{31}c_{41}^2c_{52}
        -c_{41}^2c_{52}c_{54}\\ 
     &   -c_{31}c_{43}c_{51} 
       +c_{42}c_{43}c_{51} +c_{41}c_{43}c_{52} +c_{45}c_{51}c_{52} -c_{31}c_{51}c_{55}\\ 
     &  +c_{42}c_{51}c_{55} -c_{41}c_{52}c_{55} -c_{33}
       + c_{31}c_{42} +  c_{31}c_{54} - c_{42}c_{54}
\end{align*}
In other words,  we have an embedding of~$\BO$ as a codimension two, smooth complete 
intersection in $\AA_\QQ^{12}$.
\end{example}

For the next examples, we let $P=\QQ[x,y,z]$
and $\OO = \{t_1,\dots,t_\mu\}$ an order ideal in~$\mathbb{T}^3$. 
In particular, with the next example we reprove a result of~\cite{Rob}, Corollary 3.13, using
the techniques of this paper.

\begin{example}\label{ex-1x3indets}
Let $P=\QQ[x,y,z]$ and let $\OO =\{ t_1,t_2,t_3 \}$, where $t_1=1$, $t_2=z$, and $t_3=z^2$.
Then the border of~$\OO$ is given by $\partial\OO = \{ b_1,\dots,b_7\}$,
where $b_1 = y$, $b_2 = x$, $b_3 = yz$, $b_4 = xz$, $b_5 = z^3$, $b_6 = yz^2$,
and $b_7 = xz^2$. Consequently, the ring $K[C] = K[c_{11}, \dots,c_{37}]$ is a polynomial ring 
in~21 indeterminates.

Since all~21 indeterminates $c_{ij}$ are exposed, this information is not useful for finding 
an optimal re-embedding of~$I(\BO)$. Instead, let us look at the cotangent equivalence classes.
The trivial equivalence class is $E_0=\{ c_{13}, c_{14}, c_{16}, c_{17}, c_{26}, c_{27}\}$.
There are four proper equivalence classes, namely $\{ c_{21}, c_{33} \}$,
$\{ c_{22}, c_{34}\}$, $\{c_{11}, c_{23}, c_{36}\}$, and $\{c_{12}, c_{24}, c_{37}\}$. 
Using Definition~\ref{def-genGFan} and Theorem~\ref{thm-shapeofSsigma}.c,
we see that the first element of $\LTGFan(I(\BO)_{\rm lin})$ is 
$$
E_0\cup \{ c_{21}\} \cup \{c_{22}\}  \cup  \{  c_{11}, c_{23}   \} \cup \{c_{12}, c_{24}\}
$$
It corresponds to  $Z = (c_{11}, c_{12}, c_{13}, c_{14},  
c_{16}, c_{17}, c_{21}, c_{22}, c_{23}, c_{24}, c_{26}, c_{27})$, and 
it turns out that the ideal~$I(\BO)$ is $Z$-separating.
Consequently, every elimination ordering $\sigma$ for~$Z$ yields a $Z$-separating
re-embedding of~$I(\BO)$. 
The reduced $\sigma$-Gr\"obner basis of $I(\BO)$ is 
$$
\begin{matrix}
\{ & c_{11} +c_{25}c_{31} +c_{33}c_{35} -c_{36},  &\qquad &  c_{21} +c_{31}c_{35} -c_{33},\\
   & c_{12} +c_{25}c_{32} +c_{34}c_{35} -c_{37},  &\qquad &  c_{22} +c_{32}c_{35} -c_{34},\\
   & c_{13} -c_{15}c_{31},  & \qquad &  c_{23} +c_{33}c_{35} -c_{36},\\
   & c_{14} -c_{15}c_{32},  & \qquad &  c_{24} +c_{34}c_{35} -c_{37},\\
   & c_{16} -c_{15}c_{33},  & \qquad &  c_{26} -c_{15}c_{31} -c_{25}c_{33},\;\;\\
   & c_{17} -c_{15}c_{34},  & \qquad &  c_{27} -c_{15}c_{32} -c_{25}c_{34} \; \}
\end{matrix}
$$
It defines an isomorphism $B_\OO \cong \QQ[c_{15}, c_{25}, c_{31}, c_{32},
c_{33},c_{34},c_{35},c_{36},c_{37}]$.
Consequently, the border basis scheme $\BO \cong \AA_\QQ^9$ is an affine cell.

More generally, it is shown in~\cite{Rob}, Corollary 3.13, that~$\BO$ is an affine
cell for every order ideal of the form $\OO = \{1,\, x_n,\, \dots,\, 
x_n^d\}$ in $K[x_1,\dots,x_n]$.
\end{example}

In our next example we find an optimal $Z$-separating re-embedding of a singular border basis scheme.

\begin{example}\label{ex-13z^2}
Let $P =\QQ[x,y,z]$, and let $\OO = \{t_1, \dots,t_5\}$, where $t_1=1$, $t_2=z$, $t_3=y$, $t_4=x$, and $t_5=z^2$. 
The border of~$\OO$ is $\partial\OO = \{ b_1,\dots,b_8\}$, where
$b_1=yz$, $b_2=xz$, $b_3=y^2$, $b_4=xy$, $b_5=x^2$, $b_6=z^3$, $b_7=yz^2$, and $b_8=xz^2$.
Consequently, we have $\mu=5$, $\nu=8$, and the ring $K[C] = K[c_{11}, \dots, c_{58}]$ 
is a polynomial ring in~$40$ indeterminates. The cotangent equivalence classes are
\begin{align*}
E_0 &= \{ c_{11}, c_{12}, c_{13},  c_{14},  c_{15},  c_{16},  c_{17},  c_{18},  c_{23},  
c_{24},  c_{25},  c_{27},  c_{28},  c_{37},  c_{38},  c_{47},  c_{48}\hskip-1pt \}\\
E_1 &= \{c_{22},\ c_{58}\} \\
E_2 &= \{ c_{21},\ c_{57}\}
\end{align*}
Using Definition~\ref{def-genGFan} and Theorem~\ref{thm-shapeofSsigma}.c,
we see that the first element of the set $\LTGFan(I(\BO)_{\rm lin})$ is 
$E_0\cup \{ c_{22}\} \cup \{c_{21}\}$. The corresponding tuple is  
$$
Z \;=\; ( c_{11}, c_{12}, c_{13}, c_{14}, c_{15},  c_{16}, c_{17}, c_{18},  
c_{21}, c_{22}, c_{23}, c_{24}, c_{25}, c_{27}, c_{28}, c_{37}, c_{38}, c_{47}, c_{48})
$$
The dimension of the principal component of~$\BO$ is~15, while the dimension of the 
cotangent space at the origin is $40-\#(Z) = 40 -19 = 21$. Therefore the origin is a singular 
point of~$\BO$. As suggested in Remark~\ref{rem-randomFforZ}, we may use an LPF-solver  
to detect and verify an explicit tuple~$F$ of natural generators of $I(\BO)$ which is $Z$-separating.
It follows that the reduced Gr\"obner basis of $\langle F \rangle$ yields a $Z$-separating re-embedding  
$B_\OO \cong \QQ[C{\setminus} Z]/J$ which is optimal.
\end{example}

At first glance, the order ideal~$\OO$ in the next example is similar to the preceding one.
But as soon as we start calculating $Z$-separating re-embeddings of~$\BO$, the
results turn out to be quite different.

\begin{example}\label{ex-13xy}
Let $P =\QQ[x,y,z]$, and let  $\OO = \{ t_1,\dots,t_5\}$, 
where $t_1=1$, $t_2=z$, $t_3=y$, $t_4=x$, and $t_5=xy$. 
Then we have $\partial\OO = \{ b_1,\dots,b_8\}$, 
where $b_1=z^2$, $b_2=yz$, $b_3=xz$, $b_4=y^2$, $b_5=x^2$, $b_6=xyz$, $b_7=xy^2$, 
and $b_8=x^2y$. As in the previous example, the ring $K[C] = K[c_{11},\dots,c_{58}]$ 
has 40 indeterminates and the set of non-exposed indeterminates turns out to be equal to
$C^\circ = \{ c_{11}, \dots, c_{18} \}$.

When we compute the linear part of~$I(\BO)$ in this setting, we find
that the trivial cotangent equivalence class has 21 elements, there are 12
basic equivalence classes, and there are three proper equivalence classes, namely
$\{ c_{21}, c_{56} \}$, $\{ c_{3,3}, c_{57} \}$, and $\{ c_{22}, c_{34}, c_{58}\}$.
In particular, it follows that $\dim_K(\m/\m^2) = 15$ which equals the dimension of
the principal component of~$\BO$. Hence the origin is a smooth point of~$\BO$.
By Corollay~\ref{cor-affinescheme}, the existence of an optimal 
re-embedding of~$I(\BO)$ would imply that $\BO$ is an affine cell.
In this case there is no optimal $Z$-separating re-embedding of~$I(\BO)$.
However, as already pointed out in Remark~\ref{rem-onlysufficient}, this does not 
exclude that $\BO$ is an affine cell.
\end{example}

%
%
\section*{Acknowledgements}
The third author thanks the University of Passau for its hospitality
during part of the preparation of this paper. The authors thank Mark
Huibregtse for useful discussions about Example~\ref{ex-theLshapeCont}.

%
%

\end{document}